\documentclass[8pt,a4paper]{article}

 \usepackage{amssymb,amsmath,amsfonts,amssymb}
 \usepackage{graphics,graphicx,color}
 -\marginparwidth \oddsidemargin -4mm \evensidemargin
 \oddsidemargin%
 \advance\hoffset by 5mm
 \usepackage{amsmath}    
 \advance\hoffset by 5mm

\newcommand\grad{{\bf \nabla}}

\newcommand\la{{\lambda}}
\newcommand\ga{{\gamma}}

\newcommand\eps{{\epsilon}}

\newcommand{\Acal}{{\cal A}}
\newcommand{\Ecal}{{\cal E}}
\newcommand{\Rr}{{\mathbb R}}

 \title{Landau-De Gennes theory of nematic liquid
crystals: \\ the Oseen-Frank limit and beyond}
 \date{\today}
\author{Apala Majumdar and Arghir Zarnescu\footnote{ Mathematical Insitute, University of Oxford, 24-29 St. Giles', OX1 3LB, U.K.}}

\newtheorem{remark}{Remark}
 \newtheorem{lemma} {Lemma}
 \newtheorem{proposition}{Proposition}
 
 \newtheorem{corollary}{Corollary}

\begin{document}
 \maketitle
\begin{abstract}
We study global minimizers of a continuum Landau-De Gennes
energy functional for nematic liquid crystals, in 
three-dimensional domains, subject to uniaxial boundary conditions. We analyze the physically relevant
limit of small elastic constant and show that global minimizers
converge strongly, in $W^{1,2}$, to a global minimizer predicted
by the Oseen-Frank theory for uniaxial nematic liquid crystals
with constant order parameter. Moreover, the convergence is
uniform in the interior of the domain, away from the singularities of the
limiting Oseen-Frank global minimizer. We obtain results on the rate of convergence
of the eigenvalues and the
regularity of the eigenvectors of the Landau-De Gennes global
minimizer. 
\par We also study the interplay between biaxiality
and uniaxiality in Landau-De Gennes global energy minimizers and obtain estimates
for various related quantities such as the biaxiality parameter
and the size of admissible strongly biaxial regions.
\end{abstract}
\section{Introduction}

Nematic liquid crystals are an intermediate phase of matter
between the commonly observed solid and liquid states of matter
\cite{dg}. The constituent nematic molecules translate freely as
in a conventional liquid but whilst flowing, tend to align along
certain locally preferred directions i.e. exhibit a certain degree
of long-range orientational order. Nematic liquid crystals break
the rotational symmetry of isotropic liquids; the resulting
anisotropic properties make liquid crystals suitable for a wide
range of physical applications and the subject of very interesting
mathematical modelling \cite{lin}.

\par There are three main continuum theories for nematic liquid
crystals \cite{lin}. The simplest mathematical theory for nematic
liquid crystals is the \emph{Oseen-Frank} theory
\cite{oseenfrank}. The Oseen-Frank theory is restricted to
uniaxial nematic liquid crystal materials (liquid crystal
materials with a single preferred direction of molecular
alignment) with constant degree of orientational order. The state
of a uniaxial nematic liquid crystal is described by a unit-vector
field, $n(x)\in S^2$, which  represents the preferred direction of molecular alignment. In
the simplest setting, the liquid crystal energy reduces to:

\begin{equation}
\mathcal{F}_{OF}[n]=\int_{\Omega} n_{i,k}(x)n_{i,k}(x)~dx,
\label{of}
\end{equation} the standard Dirichlet energy  for vector-valued maps into the unit sphere. The equilibrium configurations
(the physically observable configurations) correspond to minimizers of
the $\mathcal{F}_{OF}$-energy, subject to the imposed boundary conditions.
 In particular, the minimizers of $\mathcal{F}_{OF}$
are examples of $S^2$-valued harmonic maps \cite{lin,virga}. The
Oseen-Frank theory has been extensively studied in the literature,
see the review \cite{brezis}, and there are
rigorous results on the existence, regularity and singularities of
Oseen-Frank minimizers.

The Oseen-Frank theory is limited in the sense that it can only
account for point defects in liquid crystal systems but not
the more complicated line and surface defects that are observed experimentally. A second,
more comprehensive theory is the continuum \emph{Ericksen theory}
\cite{ericksen}. The Ericksen theory is also restricted to
uniaxial liquid crystal materials but can account for spatially
varying orientational order i.e. the state of the liquid crystal
is described by a pair, $\left(s,n\right)\in \Rr \times S^2$,
where $s\in \Rr$ is a real scalar order parameter that measures
the degree of orientational ordering and $n$ represents the
direction of preferred molecular alignment. In the simplest
setting, the corresponding energy functional is given by
\begin{equation}\mathcal{F}_{E}[s,n]=\int_\Omega s(x)^2|\nabla n(x)|^2+k|\nabla
s(x)|^2+W_0(s)\ dx \label{e}
\end{equation} where $k$ is a material-dependent elastic constant
and $W_0(s)$ is a bulk potential. The Ericksen theory is based on
the premise that $s$ vanishes wherever $n$ has a singularity and
this theory can account for all physically observable
defects.
\par
\par However Ericksen recognizes that his theory is but a simplified description  of
a possibly  more complex
situation (see \cite{ericksen}):
\par {\it ``There is the third possibility, that the three eigenvalues
of Q are all distinct, giving what are called biaxial nematic
configurations. Theories fitting MACMILLAN'S [11] format permit any
of the three types of configurations to occur. Certainly it is not
unreasonable to think that flows or other influences could convert a
rather stable nematic configuration to one of the biaxial type, etc.
I [19] am one of those who have argued that, near isotropic-nematic
phase transitions, it should be quite easy to induce such changes.
Accounting for such possibilities does add significant complications
to the equations and the problems of analyzing them. Experimental
information concerning the biaxial configurations is still quite
slim and, for me, it is too early to think seriously about them. So,
I will develop a theory representing a kind of compromise.''}

\par The most general continuum theory for nematic liquid crystals
is the \emph{Landau-De Gennes} theory \cite{dg,newtonmottram} which can account for uniaxial and biaxial
phases (biaxiality implies the existence of more than one
preferred direction of molecular alignment). Indeed, this theory was one of the major reasons for awarding P.G. De Gennes a Nobel prize for physics in 1991.
In the Landau-De Gennes framework, the
state of a nematic liquid crystal is modelled by a symmetric,
traceless $3\times 3$ matrix $Q\in M^{3\times 3}$, known as the
\emph{$Q$-tensor order parameter}. A nematic liquid crystal is
said to be (a) \emph{isotropic} when $Q=0$, (b) \emph{uniaxial} when the
$Q$-tensor has two equal non-zero eigenvalues; a uniaxial $Q$-tensor can be
written in the special form \begin{equation} \label{eq:uniaxial}
Q=s\left(n\otimes n-\frac{1}{3}Id\right);~ s\in \Rr\setminus\{0\},~n\in S^2
\end{equation} and (c) \emph{biaxial} when $Q$ has three distinct
eigenvalues; a biaxial $Q$-tensor can always be represented as follows
(see Proposition~\ref{prop:rep1})
\begin{equation}
\label{eq:biaxial} Q =s\left(n\otimes n-\frac{1}{3}Id\right) +
r\left(m\otimes m-\frac{1}{3}Id\right) \quad s,r \in \Rr;~ n,m \in
S^2.
\end{equation}

The Landau-De Gennes energy functional, $\mathcal{F}_{LG}[Q]$, is
a nonlinear integral functional of $Q$ and its spatial
derivatives. We work with the simplest form of
$\mathcal{F}_{LG}[Q]$, with Dirichlet boundary conditions, $Q_b$
(refer to (\ref{eq:Qb})), on  three-dimensional domains
$\Omega \subset \Rr^3$. We take $\mathcal{F}_{LG}[Q]$ to be
\cite{smallelastic}
\begin{equation}
\label{energy}
 \mathcal{F}_{LG}[Q]=\int_{\Omega} \frac{L}{2}|\nabla Q|^2(x)+f_B(Q(x))\,dx
\end{equation}
where $f_B(Q)$ is the bulk energy density that accounts for
bulk effects, $\left| \grad Q\right|^2 =\sum_{i,j,k=1}^3 Q_{ij,k}Q_{ij,k}$ 
is the elastic energy density that penalizes spatial
inhomogeneities and $L>0$ is a material-dependent elastic
constant. We take $f_B(Q)$ to be a quartic polynomial in the
$Q$-tensor components, since this is the simplest form of $f_B(Q)$ that allows for multiple local minima and a first-order nematic-isotropic phase transition \cite{dg,virga}.
This form of $f_B(Q)$ has been widely-used in the  literature and is defined as follows
$$f_B(Q) = \frac{\alpha(T- T^*)}{2}\textrm{tr}\left(Q^2\right) -
\frac{b}{3}\textrm{tr}\left(Q^3\right) +
\frac{c}{4}\left(\textrm{tr}Q^2\right)^2$$ where $\alpha,b,c \in
\Rr$ are material-dependent positive constants, $T$ is the
absolute temperature and $T^*$ is a characteristic liquid crystal
temperature. We work in the low-temperature regime $T < T^*$ for
which $\alpha (T - T^*) < 0$. Keeping this in mind, we recast the
bulk energy density as follows:
\begin{equation}
f_B(Q) = -\frac{a^2}{2}\textrm{tr}\left(Q^2\right) -
\frac{b^2}{3}\textrm{tr}\left(Q^3\right) +
\frac{c^2}{4}\left(\textrm{tr}Q^2\right)^2, \label{eq:10}
\end{equation} where $a^2,b^2,c^2 \in \Rr^+$ are material-dependent and
temperature-dependent
positive constants. The equilibrium configurations (the physically observable
configurations) then correspond to minimizers of
$\mathcal{F}_{LG}[Q]$, subject to the imposed boundary conditions.

In the first part of the paper, we study the  the limit of vanishing elastic
constant $L\to 0$  for global minimizers,
$Q^{(L)}$, of $\mathcal{F}_{LG}[Q]$. This study is in the spirit of the asymptotics
for minimizers of Ginzburg-Landau functionals for superconductors
\cite{bbh}. The limit $L\to 0$ is a physically relevant limit
since the elastic constant $L$ is typically very small, of the order of $10^{-11}$ Joule/metre.
\cite{smallelastic}.

 We define a  \emph{limiting harmonic map} $Q^{(0)}$ as
follows $$ Q^{(0)} = s_+\left(n^{(0)} \otimes n^{(0)} -
\frac{1}{3}Id\right)$$ where $s_+$ is defined in
(\ref{s+}), $n^{(0)}$ is a minimizer of the Oseen-Frank energy,
$\mathcal{F}_{OF}[n]$ in (\ref{of}), subject to the fixed boundary
condition $n=n_b\in C^\infty(\partial\Omega,\mathbb{S}^2)$  and $Q_b$ and $n_b$ are related as in (\ref{eq:Qb}).
Our main results are:

$\bullet$ There exists a sequence of global minimizers $\left\{
Q^{(L_k)}\right\}$ such that $Q^{(L_k)} \stackrel{L_k \to
0^+}\longrightarrow Q^{(0)}$ strongly in the Sobolev space $W^{1,2}(\Omega,\mathbb{R}^9)$.

$\bullet$ The sequence $\left\{Q^{(L_k)}\right\}$ as above converges uniformly to $Q^{(0)}$ as $L_k\to 0$, in the
interior of $\Omega$, away from the (possible) singularities of $Q^{(0)}$. 

$\bullet$ The bulk energy density, $f_B\left(Q^{(L_k)}\right)$, converges
uniformly to its minimum value away from the (possible) singularities of
$Q^{(0)}$; the uniform convergence of the bulk energy density holds in
the interior and up to the boundary, away from the (possible) singularities
of $Q^{(0)}$.

These results show that the predictions of the Oseen-Frank theory
(described by the limiting map $Q^{(0)}$) and the Landau-De Gennes
theory agree away from the singularities of $Q^{(0)}$. The
global minimizers, $Q^{(L)}$, are real analytic (see
Proposition~\ref{analyticity}) and have no singularities as such. However, one of the most intriguing features of  nematic liquid crystals are the optical `defects' that appear in the Schlieren textures \cite{dg}. From a physical point of view, these defects are regions of
rapid changes in the configurational properties of a nematic
liquid crystal \cite{dg}. We conjecture that certain types of optical defects in $Q^{(L_k)}$ (for small $L_k$),
when they exist, may be localized near the analytic singularities of the limiting
map $Q^{(0)}$, since $Q^{(L_k)}$ can have strong variations only near  the  singularities of $Q^{(0)}$ (more precisely,  the gradient, $\grad Q^{(L_k)}$, cannot be bounded independently of $L_k$ on any set containing  a singularity of $Q^{(0)}$).   
 There is existing literature on the location of singularities in harmonic maps \cite{lieb} and this may allow one to predict the location of (optical) defects in
a global Landau-De Gennes minimizer. 


Our convergence results analyze the limit of vanishing elastic constant $L \to 0$. 
Physical situations are modelled by small {\it but non-zero} values of the elastic constant $L$. Thus our convergence results show that for $L$ sufficiently small, the limiting harmonic map $Q^{(0)}$ provides but a `rough' description of  $Q^{(L)}$ i.e. $Q^{(L)}$ can be thought of as having  a `leading'
uniaxial part plus a small biaxial perturbation, away from the
singularities of $Q^{(0)}$.   This small biaxial perturbation is  of order $O(\sqrt{L})$ where $L << 1$
(see Section~\ref{sec:consequence} for details). However, numerical simulations show that biaxiality may become prominent 
in the vicinity of defects \cite{mg,rosso&virga}.  In the second part of our paper, we
study biaxiality and their role in global minimizers $Q^{(L)}$, noting that biaxiality (if it exists) is one of the main differences between $Q^{(L)}$ and the limiting approximation $Q^{(0)}$. More precisely, in  Propositions~\ref{prop:isotropic}
and \ref{prop:nondefect}, we obtain estimates  for the size of the regions where $Q^{(L_k)}$ can deviate significantly from $Q^{(0)}$ and on the size of admissible
strongly biaxial regions in $Q^{(L)}$, in terms of the biaxiality parameter
$\beta$ (defined in (\ref{eq:am2})) and the
material-dependent constants. While Proposition~\ref{prop:isotropic} may be relevant to the properties of $Q^{(L)}$ near the singular set of $Q^{(0)}$, Proposition~\ref{prop:nondefect} is relevant to the equilibrium properties away from the singular set of $Q^{(0)}$.

Using a simple nearest-neighbour projection
argument (see Corollary~\ref{eigenvector}), we show that the
`leading eigendirection', corresponding to
the leading uniaxial part (see Section \ref{sec:consequence} for definitions)
is smooth on any compact set $K$ not containing  any singularity of $Q^{(0)}$. Further, in Proposition~\ref{zeroL}, we also
show that $Q^{(L)}$ is either (a) uniaxial everywhere (except for possibly a set of measure zero where $Q$ can be isotropic) or (b) $Q^{(L)}$
is biaxial everywhere and can be uniaxial or isotropic only on sets of measure zero. It is known that as long as the number
of distinct eigenvalues does not change, the eigenvectors of $Q^{(L)}$
enjoy the same degree of regularity as $Q^{(L)}$ itself \cite{nomizu}. In Corollary~\ref{eigenvectornew}, we show that the eigenvectors are necessarily smooth everywhere except for possibly a zero-measure set where the number of distinct eigenvalues changes and therefore, if the
eigenvectors of $Q^{(L)}$ suffer any discontinuities, these
discontinuities must be localized on the uniaxial-biaxial, uniaxial-isotropic or biaxial-isotropic
interfaces. This result may be relevant to the
interpretation of optical data from experiments and we hope to explore this connection in future work.

\par Finally, we note that the Landau-De Gennes theory for uniaxial liquid crystal materials has strong analogies
with the $3D$ version of the  Ginzburg-Landau theory for superconductors \cite{bbh} . The Ginzburg-Landau energy functional for a three-dimensional vector field,
$u:\Omega \to \Rr^3$, is typically of the form
\begin{equation}
\label{eq:gl}
\mathcal{F}_{GL}[u] = \int_{\Omega} \frac{1}{2}\left|\grad u\right|^2 + \frac{1}{4\eps^2}\left(1 - |u|^2\right)^2~dx
\end{equation}
where $\eps>0$ is a very small parameter. The functional $\mathcal{F}_{GL}$ has been rigorously studied in the limit $\eps\to 0$ which is analogous to the limit $L\to 0$ in our problem.
 The new mathematical complexities in the Landau-De Gennes theory for nematic liquid crystals come from the high dimensionality of the target space and also
 from the  possibility of biaxiality in global energy minimizers. Future challenges include a better understanding of the qualitative properties of global minimizers for small but non-vanishing values of $L$, a better description of $Q^{(L)}$ near the singularities of the limiting harmonic map $Q^{(0)}$, the regularity of the eigenvectors and eigenvalues, along with a  deeper understanding of the appearance and role of biaxiality in global minimizers.

The paper is organized as follows - in Section~\ref{sec:conventions}, we introduce the conventions and notations that are used in the rest of the paper. In
Section~\ref{sec:preliminaries}, we state two representation
formulae for $Q$-tensors that are useful for subsequent computations in later sections. In
Section~\ref{sec:main}, we study the properties of global energy
minimizers in the limit $L \to 0$ and
prove the convergence results. In Section~\ref{sec:consequence},
we discuss the consequences of our convergence results and their relevance to
the bulk energy density, the biaxiality parameter, the eigenvalues and the eigenvectors of
a global Landau-De Gennes minimizer. In Section~\ref{sec:6}, we derive estimates for the bulk energy density, obtain bounds for the size of admissible strongly biaxial regions  and discuss the interplay between biaxiality and uniaxiality in a global energy minimizer.

\section{Preliminaries}
\label{sec:conventions}

\par We take our domain, $\Omega \subset \Rr^3$, to be bounded and
simply-connected with smooth boundary, $\partial\Omega$. Let
$S_0\subset \mathbb{M}^{3\times 3}$ denote the space of Q-tensors,  i.e.
\begin{displaymath}
 S_0\stackrel{def}{=} \left\{Q \in \mathbb{M}^{3\times 3};
Q_{ij}=Q_{ji},~Q_{ii} = 0 \right\}
\end{displaymath}
where we have used the Einstein summation convention; the Einstein convention will be assumed in the rest of the paper. The
corresponding matrix norm is defined to be
\begin{displaymath}
 \left| Q \right|\stackrel{def}{=}\sqrt{\textrm{tr}Q^2} =\sqrt{ Q_{ij}
Q_{ij}}.
\end{displaymath}

As stated in the introduction, we take the bulk energy density term to be   
$$f_B(Q) = -\frac{a^2}{2}\textrm{tr}\left(Q^2\right) -
\frac{b^2}{3}\textrm{tr}\left(Q^3\right) +
\frac{c^2}{4}\left(\textrm{tr}(Q^2)\right)^2$$ where $a^2,b^2,c^2\in \Rr$ are material-dependent and temperature-dependent
positive constants. One can readily verify that $f_B(Q)$ is
bounded from below (see  Proposition~\ref{prop:bulk}, \cite{bm}), and we
define a non-negative bulk energy density, $\tilde{f}_B(Q)$, that differs from
$f_B(Q)$ by an additive constant as follows:
\begin{equation}
\label{eq:bulknew}\tilde{f}_B(Q) = f_B(Q) - \min_{Q\in
S_0}f_B\left(Q\right).
\end{equation} It is clear that
$\tilde{f}_B(Q)\geq 0$ for all $Q\in S_0$ and the set of
minimizers of $\tilde{f}_B(Q)$ coincides with the set of
minimizers for $f_B(Q)$. In Proposition~\ref{prop:bulk}, we show
that the function $\tilde{f}_B(Q)$ attains its minimum on the set of uniaxial
$Q$-tensors with constant order parameter $s_+$ as shown below
\begin{eqnarray}
\label{eq:Qmin} && \tilde f_B(Q)=0\Leftrightarrow Q\in Q_{min}~ where \nonumber \\ &&
 Q_{\min} = \left\{Q\in S_0, Q= s_+ \left( n\otimes n -
\frac{1}{3}Id \right), n\in\mathbb{S}^2 \right\}\end{eqnarray}with \begin{equation}
s_+=\frac{b^2+\sqrt{b^4+24a^2c^2}}{4c^2}. \label{s+}
\end{equation}

We work with Dirichlet boundary conditions, referred to as
\emph{strong anchoring} in the liquid crystal literature
\cite{dg}. The boundary condition $Q_b\in Q_{\min}$ is smooth and is given by
\begin{equation}
\label{eq:Qb} Q_b = s_+\left(n_b \otimes n_b -
\frac{1}{3}Id\right),\, n_b \in C^\infty\left(\partial\Omega;S^2\right).
\end{equation}
 We define our admissible space to be
\begin{eqnarray}
&& \Acal_{Q} = \left\{Q\in W^{1,2}\left(\Omega;S_0\right);
\textrm{$Q=Q_b$ on $\partial\Omega$, \textrm{ with $Q_b$ as in }(\ref{eq:Qb})}\right\}\label{eq:am41},
\end{eqnarray} where $W^{1,2}\left(\Omega;S_0\right)$ is the
Sobolev space of square-integrable $Q$-tensors with
square-integrable first derivatives \cite{evans}. 
The corresponding
$W^{1,2}$-norm is given by $\| Q
\|_{W^{1,2}(\Omega)} =\left( \int_{\Omega}
|Q|^2 + |\grad Q|^2~dx\right)^{1/2}.$ In addition to the $W^{1,2}$-norm, we also
use the $L^{\infty}$-norm in this paper, defined to be $\|Q\|_{L^{\infty}(\Omega)} =
\textrm{ess sup}_{x\in\Omega}|Q(x)|$ .

We study global minimizers of a modified Landau-De Gennes energy functional, $\tilde{F}_{LG}[Q]$, 
in the admissible space $\Acal_{Q}$. The functional $\tilde{F}_{LG}[Q]$ differs from $\mathcal{F}_{LG}[Q]$ in (\ref{energy})
by an additive constant and is defined to be  
\begin{equation}\tilde F_{LG}[Q]=\int_{\Omega}
\frac{L}{2}Q_{ij,k}(x)Q_{ij,k}(x)+\tilde f_B(Q(x))\,dx
\label{LDGfunctional}.
\end{equation} 
For a fixed $L>0$, let $Q^{(L)}$ denote a global minimizer of
$\tilde{F}_{LG}[Q]$ in the admissible class, $\Acal_Q$. 
The existence of $Q^{(L)}$ is immediate from the
direct methods in the calculus of variations \cite{evans}.  The bulk energy
density, $\tilde f_B(Q)$, is bounded from below, the energy density is convex in
$\grad Q$ and therefore, $\tilde F_{LG}[Q]$ is weakly
sequentially lower semi-continuous. Moreover, it is clear that $\tilde{F}_{LG}[Q]$ and $\mathcal{F}_{LG}[Q]$ have the same
set of global minimizers for a fixed set of material-dependent and temperature-dependent 
constants $\left\{a^2, b^2, c^2, L\right\}$.

The global minimizer $Q^{(L)}$ is a weak solution of the corresponding Euler-Lagrange
equations \cite{bm}
\begin{equation}
L\Delta
Q_{ij}=-a^2Q_{ij}-b^2\left(Q_{ik}Q_{kj}-\frac{\delta_{ij}}{3}\textrm{tr}(Q^2)\right)
+c^2Q_{ij}\textrm{tr}(Q^2)\,~ ~i,j=1,2,3 \label{ELeq}.
\end{equation} where the term $b^2\frac{\delta_{ij}}{3}\textrm{tr}(Q^2)$ is a Lagrange multiplier 
that enforces the tracelessness constraint. It follows from standard arguments in elliptic
regularity that $Q^{(L)}$ is actually a classical solution of
(\ref{ELeq}) and $Q^{(L})$ is smooth and real analytic (see also Section \ref{section:analyticitybiaxiality}).

Finally, we introduce a \textit{``limiting uniaxial harmonic map''}
$Q^{(0)}:\Omega \to Q_{min}$; 
$Q^{(0)}$ is defined to be a global minimizer (not necessarily unique) of
$\tilde F_{LG}[Q]$ in the restricted class, 
$\Acal_{Q}\cap \{Q:\Omega\to S_0,\,Q(x)\in Q_{min}\,\textrm{a.e.}\,x\in\Omega\}$. Then
 $Q^{(0)}$ is necessarily of the form
\begin{equation}
\label{eq:Q0} Q^{(0)} = s_+\left(n^{(0)} \otimes n^{(0)} -
\frac{1}{3}Id\right),\end{equation} where $n^{(0)}$ is a global
minimizer of $\mathcal{F}_{OF}[n]$ (see \cite{bz}, \cite{bc}),
\begin{equation}
\label{eq:n0} \int_{\Omega}|\grad n^0(x)|^2\,dx = \min_{n \in
\Acal_n}\int_{\Omega}|\grad n(x)|^2\,dx
\end{equation}
in the admissible class $\Acal_n = \left\{ n \in
W^{1,2}\left(\Omega;S^2\right);~n = n_b
~on~\partial\Omega\right\}$ and $n_b$ and $Q_b$ are related as in (\ref{eq:Qb}). 
This ``limiting harmonic'' map $Q^{(0)}$ is 
therefore obtained from an energy minimizer, $n^0$, (not necessarily unique) within the Oseen-Frank theory for
uniaxial nematic liquid crystals with constant order parameter (for more results about the relation between $n^{(0)}$ and $Q^{(0)}$ see \cite{bz}). It
follows from standard results in harmonic maps \cite{virga} that
$Q^{(0)}$ has at most a finite number of isolated point singularities 
(points where $n^{(0)}$ has singularities). In the following sections we will elaborate on the relation between $Q^{(L)}$ and $Q^{(0)}$.

\bigskip\section{Representation formulae for  $Q$-tensors}
\label{sec:preliminaries}

\par We have:

\begin{proposition}\label{prop:rep1}
A  matrix $ Q\in S_0$ can be represented in the form
\begin{equation}
Q=s(n\otimes n-\frac{1}{3}Id)+r(m\otimes m-\frac{1}{3}Id)
\label{Qrep}
\end{equation} with $n$ and $m$ unit-length eigenvectors of $Q$, $n\cdot
m=0$ and
\begin{equation}
0\le r\le \frac{s}{2}\,\,\textrm{ or } \frac{s}{2}\le r\le
0\label{rs}
\end{equation}
\par The scalar order parameters $r$ and $s$ are piecewise linear combinations of the
eigenvalues of  $Q$. \label{prop:Qrep}
\end{proposition}
\par{\bf Proof.} From the spectral decomposition theorem we have
 \begin{equation}
Q = \lambda_1 n_1\otimes n_1 + \lambda_2 n_2\otimes n_2
+\lambda_3n_3 \otimes n_3
\end{equation} where $\lambda_1, \lambda_2,\lambda_3$ are  eigenvalues of $Q$
and $n_1,n_2,n_3$ are the corresponding unit eigenvectors,
pairwise perpendicular. We have $I=\sum_{i=1}^{3} n_i\otimes n_i$
and the tracelessness condition implies that
$\lambda_1+\lambda_2+\lambda_3=0$. Thus
$$Q=\lambda_1 n_1\otimes n_1+\lambda_2n_2\otimes
n_2-(\lambda_1+\lambda_2)(I-n_1\otimes n_1-n_2\otimes n_2)$$

\par We consider six regions $R_i^+,i=1,\dots,6$ in the
$(\lambda_1,\lambda_2)$ - plane
which cover exactly half of the whole plane. This corresponds to
the representation (\ref{Qrep}) with $0\leq r \leq \frac{s}{2}$.
The other half of the plane is covered by the regions
$R_i^{-},i=1,\dots,6$ , (which are obtained by reflecting $R_i^+$
through the origin $(0,0)$) and the regions $R_i^{-}$ correspond
to the representation (\ref{Qrep}),with $r,s\le 0$.

\par We let
$R_1^+=\{(\lambda_1,\lambda_2)\in\mathbb{R}^2,-2\lambda_1\le
\lambda_2,\lambda_1\le 0\}$. In this case
$r\stackrel{def}{=}2\lambda_1+\lambda_2$ and
$s\stackrel{def}{=}2\lambda_2+\lambda_1$ with
 $n\stackrel{def}{=}n_2,m\stackrel{def}{=}n_1$.  One can directly verify
that for
$r,s$ thus defined, we have
 $$ r = 2\la_1 + \la_2 \leq \frac{s}{2} = \la_2 + \frac{\la_1}{2}.$$
 Interchanging $\lambda_1$
 with $\lambda_2$ in the definition of $r$ and $s$ and $m$ with $n$, we
 obtain the region $R_2^+=\{(\lambda_1,\lambda_2);\lambda_2\ge
 -\lambda_1/2;\lambda_2\le 0\}$.

 \par Let $R_3^+=\{(\lambda_1,\lambda_2)\in\mathbb{R}^2,\lambda_2\le
 0,\lambda_2\ge \lambda_1\}$. Taking
 $r\stackrel{def}{=}\lambda_2-\lambda_1,s\stackrel{def}{=}-2\lambda_1-\lambda_2$,
 $n\stackrel{def}{=}n_3,m\stackrel{def}{=}n_2$, one can check that
 $$ r = \la_2 - \la_1 \leq \frac{s}{2}=-\la_1 - \frac{\la_2}{2}.$$ The
 region $R_4^+$ is obtained from interchanging $\lambda_1$ and
 $\lambda_2$.

 \par We have
 $R_5^+=\{(\lambda_1,\lambda_2)\in\mathbb{R}^2,\lambda_1\le
 0,-2\lambda_1\ge\lambda_2\ge -\lambda_1\}$ with
 $r\stackrel{def}{=}-2\lambda_1-\lambda_2,s\stackrel{def}{=}\lambda_2-\lambda_1$,
 $n=n_2$ and $m\stackrel{def}{=}n_3$. Again, it is straightforward to
check that
 $$ r = -2\la_1 - \la_2 \leq \frac{\la_2}{2} - \frac{\la_1}{2}.$$
 Interchanging $\lambda_1$ with $\lambda_2$, we obtain the region $R_6^+$.

 \par
 Finally the remaining half of the $\left(\la_1,\la_2\right)$-plane is
covered by
 the regions $R_i^{-}$ (obtained from $R_i^{+}$ by changing the signs of the
inequalities and keeping
 the definitions of $r$ and $s$ unchanged). For example, $R_1^{-}$
 is defined to be
 $$ R_1^{-} = \left\{\left(\la_1,\la_2\right) \in \Rr^2; \la_1
 \geq 0,~2\la_1 \leq - \la_2 \right\}$$
 with $r=2\la_1 + \la_2$ and $s=2\la_2 + \la_1$. One can then
 directly check that
 $$ \frac{s}{2}\leq r \leq 0.$$ The remaining five regions $R_i^-$
 for $i=2\ldots 6$ can be defined analogously.
 $\Box$
 
 \begin{remark}
 The representation formula (\ref{Qrep}) is known in the literature
\cite{newtonmottram}. In Proposition~\ref{prop:rep1}, we state
that it suffices to consider the two cases given by (\ref{rs}); we
have not found references for this fact.
\end{remark}

In Proposition~\ref{prop:secondrep}, we state a second
representation formula for admissible $Q\in S_0$ and its relation
to the representation formula (\ref{Qrep}). The representation
formula (\ref{eq:second1}) is  known in the literature
\cite{virgadematteis} and will be used in
Section~\ref{sec:consequence}. For reader's convenience we provide a quick proof.

\begin{proposition} \label{prop:secondrep} (A second representation formula) 
A matrix $Q\in S_0$ can be represented  as:
\begin{equation}
\label{eq:second1} Q = S\left(n\otimes  n -
\frac{1}{3}Id\right) + R \left(m\otimes m - p\otimes p\right)
\end{equation}

\par The vectors $n, m$ and $p$ are unit-length and pairwise perpendicular eigenvectors of $Q$ with corresponding eigenvalues
$\lambda_1,\lambda_2,\lambda_3$.  The scalar order parameters $S$
and $R$ are given by
\begin{equation}
\label{eq:second2} S = 3 \frac{\la_1}{2} \qquad R =
\frac{1}{2}\left(2\la_2 + \la_1 \right).
\end{equation}
\end{proposition}

\textbf{Proof.}  We have the spectral decomposition of $Q$, namely
$$Q=\lambda_1 n\otimes n+\lambda_2 m\otimes m+\lambda_3 p\otimes p$$ with $n,m,p$ pairwise perpendicular unit-length eigenvectors of $Q$ and
$$Id=n\otimes n+m\otimes m+p\otimes p.$$ 

Combining the last two relations and taking  $S = 3 \frac{\la_1}{2}, R =
\frac{1}{2}\left(2\la_2 + \la_1 \right)$ we obtain the claim.$\Box$

\section{The limiting harmonic map}
\label{sec:main}
\subsection{ The uniform convergence in the interior}

Firstly, we recall that for a $Q\in S_0$  the biaxiality parameter $\beta(Q)$ (see for instance \cite{mg}) is defined to be
\begin{equation}
\label{eq:am2} \beta(Q) = 1 -
6\frac{\left(\textrm{tr}Q^3\right)^2}{\left(\textrm{tr}Q^2\right)^3}
\end{equation}

\par The  significance of  $\beta(Q)$ as a measure of biaxiality is due to the following
\begin{lemma}
\label{betalemma}
 (i) The biaxiality parameter $\beta(Q)\in\left[0,1\right]$ and
$\beta(Q)=0$ if and only if $Q$ is purely uniaxial i.e. if $Q$ is
of the form, $Q=s\left(n\otimes n-\frac{1}{3}Id\right)$ for some
$s\in\mathbb{R},n\in\mathbb{S}^2$. (ii) The biaxiality parameter,
$\beta(Q)$, can be bounded in terms of the ratio $\frac{r}{s}$,
where $(s,~r)$ are the scalar order parameters in
Proposition~\ref{prop:rep1} . These bounds are given by
\begin{equation}
\label{eq:am3} \frac{1}{2}\left(1 - \sqrt{1 - \sqrt{\beta}}\right)
\leq \frac{r}{s} \leq \frac{1}{2}\left(1 + \sqrt{1 -
\sqrt{\beta}}\right).
\end{equation} Equivalently,
\begin{equation}
\label{eq:am3new} \frac{1 -\sqrt{1 - \sqrt{\beta}}}{3 + \sqrt{1 -
\sqrt{\beta}}} \leq \frac{R}{S}\leq \frac{1 +\sqrt{1 -
\sqrt{\beta}}}{3 - \sqrt{1 - \sqrt{\beta}}} \end{equation} where
$(S,~R)$ are the order parameters in
Proposition~\ref{prop:secondrep}. Further $\beta(Q)=1$ if and only
if $r=\frac{s}{2}$ or if and only if $\frac{R}{S}=\frac{1}{3}$.
(iii)For an arbitrary $Q\in S_0$, we have that
\begin{equation}
-\frac{|Q|^3}{\sqrt{6}}\left(1 - \frac{\beta}{2}\right) \leq
\textrm{tr} Q^3 \leq \frac{|Q|^3}{\sqrt{6}}\left(1 -
\frac{\beta}{2}\right).\label{eq:am9}
\end{equation}
\end{lemma}

\textbf{Proof:} The proof of Lemma~\ref{betalemma} is deferred to
the Appendix. $\Box$

\par The next proposition gives us apriori $L^\infty$ bounds, independent of $L$.
\begin{proposition}
\label{prop:max} Let $\Omega \subset \mathbb{R}^3$ be a
bounded and simply-connected  open set with smooth boundary. Let $Q^{(L)}$
be a global minimizer of the Landau-De Gennes energy functional (\ref{LDGfunctional})  , in the space
(\ref{eq:am41}).

 Then
\begin{equation}
\|Q^{(L)} \|_{L^{\infty}(\Omega)} \leq
\sqrt{\frac{2}{3}}s_+ \label{eq:am42}
\end{equation}where $s_+$ is defined in (\ref{s+}).
\end{proposition}

\textbf{Proof.} Proposition~\ref{prop:max} has been proven in
\cite{bm}; we reproduce the proof here for completeness.

The proof  proceeds by contradiction. In the following we drop the superscript $L$ for convenience.
We assume that there exists a point $x^*\in\bar\Omega$
where $|Q|$ attains its maximum and $|Q(x^*)| >\sqrt{\frac{2}{3}}s_+$. On
$\partial \Omega$, $|Q| = \sqrt{\frac{2}{3}}s_+$ by our choice of the
boundary condition $Q_b$ (note that if $Q\in Q_{min}$ then $|Q|=\sqrt{\frac{2}{3}}s_+$). If $Q$ is a global
minimizer of $\tilde F[Q]$  then $Q$ is a classical
solution (see Section ~\ref{section:analyticitybiaxiality} for regularity)  of the Euler-Lagrange equations 
\begin{equation} \label{eq:am43} L\Delta
Q_{ij} = -a^2 Q_{ij} - b^2\left(Q_{ip}Q_{pj} -
\frac{1}{3}\textrm{tr}Q^2\delta_{ij}\right) + c^2
\left(\textrm{tr}Q^2\right) Q_{ij}.
\end{equation}

Since the function $|Q|^2:\bar{\Omega}\to \Rr$ must attain its
maximum at $x^*\in\Omega$,  we necessarily have that

\begin{equation} \label{eq:am44} \Delta
\left(\frac{1}{2}|Q|^2\right)(x^*)\leq 0 
\end{equation}

We multiply both sides of (\ref{eq:am43}) by $Q_{ij}$ and obtain
\begin{equation}
\label{eq:am46} L~\Delta\left(\frac{1}{2}|Q|^2\right) = -a^2
\textrm{tr}Q^2 - b^2\textrm{tr}Q^3 + c^2
\left(\textrm{tr}Q^2\right)^2 + L|\nabla Q|^2.
\end{equation} We  note that
\begin{equation}
\label{eq:am47} -a^2 \textrm{tr}(Q^2) - b^2\textrm{tr}(Q^3) +
c^2 \left(\textrm{tr}(Q^2)\right)^2 \geq f(|Q|)
\end{equation}
where
\begin{equation}
\label{eq:am48} f(|Q|) = -a^2 |Q|^2 -
\frac{b^2}{\sqrt{6}}|Q|^3 + c^2 |Q|^4,
\end{equation}
since $\textrm{tr}(Q^3) \leq \frac{|Q|^3}{\sqrt{6}}$ from
(\ref{eq:am9}). One can readily verify that
\begin{equation}
\label{eq:am49} f(|Q|) > 0 \quad \textrm{for $|Q| >
\sqrt{\frac{2}{3}}s_+$}
\end{equation}
which together with (\ref{eq:am46}) and (\ref{eq:am47}) imply that
\begin{equation}
\label{eq:am50} \Delta\left(\frac{1}{2}|Q|^2\right) (x)> 0
\end{equation}
for all interior points $x\in\Omega$, where $|Q(x)|> \sqrt{\frac{2}{3}}s_+$.
This contradicts  (\ref{eq:am44}) and thus gives the conclusion. $\Box$

\par In what follows, let $e_L(Q(x))$ denote the energy density $e_L(Q(x))\stackrel{def}{=}\frac{1}{2}|\nabla Q|^2+\frac{\tilde f_B(Q(x))}{L}$. 
We consider the normalized energy on balls $B(x,r)\subset \Omega
= \left\{y\in\Omega; \left| x - y \right| \leq r \right\}$

\begin{equation}
\label{eq:normalized}
\mathcal{F}(Q,x,r)\stackrel{def}{=}\frac{1}{r}\int_{B(x,r)} e_L(Q(x))\,dx=\frac{1}{r}\int_{B(x,r)}\frac{\tilde{f}_B(Q)}{L}
+ \frac{1}{2}\left|\grad Q \right|^2\,dx.
\end{equation}

\par We have:

\begin{lemma} (Monotonicity lemma)
\label{lem:mon} Let $Q^{(L)}$ be a global minimizer of
$\tilde F_{LG}[Q]$ in (\ref{LDGfunctional}). Then
\begin{equation}
\mathcal{F}(Q^{(L)},x,r)\le\mathcal{F}(Q^{(L)},x,R), \forall
x\in\Omega,r\le R, \textrm{ so that } B(x,R)\subset\Omega
\end{equation}
\label{lemma:mon1}
\end{lemma}

\smallskip\par{\bf Proof.} The proof follows a standard pattern (see for instance \cite{linriviere}) and is a consequence of the Pohozaev identity. 
We assume, without loss of generality, that $x=0$ and
$0<R<d(0,\partial\Omega)$, where $d$ denotes the Euclidean
distance. Since $Q^{(L)}$ is a global energy minimizer, it is a
classical solution (see Section ~\ref{section:analyticitybiaxiality} for regularity) of the system (\ref{ELeq}):
\begin{equation}
\Delta Q_{ij}=\frac{1}{L}\left[\frac{\partial\tilde
f_B(Q)}{\partial
Q_{ij}}+b^2\frac{\delta_{ij}}{3}\textrm{tr}(Q^2)\right]
\label{Eqgrad}
\end{equation} In (\ref{Eqgrad}) and in what follows, we drop the
superscript $L$ for convenience.

\par We multiply (\ref{Eqgrad}) by $x_k\cdot\partial_k Q_{ij}$,
sum over repeated indices and integrate over $B(0,R)$ to obtain
the following 

\begin{eqnarray}
0=\int_{B(0,R)}Q_{ij,ll}(x)\cdot x_k\cdot\partial_k
Q_{ij}(x)-\frac{1}{L}\frac{\partial \tilde {f_B}(Q(x))}{\partial
Q_{ij}}\cdot x_k\cdot\partial_k
Q_{ij}(x)-\frac{1}{L}b^2\frac{\delta_{ij}}{3}\textrm{tr}(Q^2(x))\cdot
x_k\cdot\partial_kQ_{ij}(x)\,dx\nonumber\\
=\underbrace{\int_{B(0,R)}Q_{ij,ll}(x)\cdot x_k\cdot\partial_k
Q_{ij}(x)\,dx}_{I}-\underbrace{\int_{B(0,R)}\frac{1}{L}\frac{\partial
\tilde f_B(Q(x))}{\partial Q_{ij}}\cdot x_k\cdot\partial_k
Q_{ij}\,dx}_{II}\label{mon1}
\end{eqnarray} where we have used the tracelessness condition $Q_{ii}=0$.
\par Integrating by parts, we have that:
\begin{eqnarray}
I=\int_{B(0,R)}Q_{ij,ll}(x) x_k\partial_k
Q_{ij}(x)\,dx\nonumber\\=-\int_{B(0,R)}Q_{ij,l}(\delta_{lk}Q_{ij,k}(x)+x_k
Q_{ij,kl}(x))dx+
\int_{\partial B(0,R)}Q_{ij,l}x_kQ_{ij,k}\frac{x_l}{R}\,dx\nonumber\\
=-\int_{B(0,R)}Q_{ij,l}(x)Q_{ij,l}(x)\,dx+3\int_{B(0,R)}\frac{1}{2}Q_{ij,l}(x)Q_{ij,l}(x)\,dx\nonumber\\-\int_{\partial
B(0,R)}\frac{Q_{ij,l}(x)Q_{ij,l}(x)}{2}\frac{x_k\cdot
x_k}{R}\,dx+\int_{\partial B(0,R)}\frac{(Q_{ij,k}(x)\cdot
x_k)^2}{R}\,dx
\end{eqnarray}

\begin{eqnarray}
II=\int_{B(0,R)}\frac{1}{L}\frac{\partial\tilde
f_B(Q(x))}{\partial Q_{ij}}\cdot x_k\cdot \partial_k
Q_{ij}(x)\,dx=\frac{1}{L}\int_{B(0,R)}\partial_k \tilde
f_B(Q(x))\cdot
x_k\,dx\nonumber\\
=-\frac{3}{L}\int_{B(0,R)}\tilde
f_B(Q(x))\,dx+\frac{1}{L}\int_{\partial B(0,R)}\tilde
f_B(Q(x))\cdot\frac{x_k\cdot x_k}{R}\,dx
\end{eqnarray} Hence (\ref{mon1}) becomes:

\begin{eqnarray}
-\int_{B(0,R)}\frac{Q_{ij,l}(x)Q_{ij,l}(x)}{2}+\frac{\tilde
f_B(Q(x))}{L}\,dx+R\int_{\partial
B(0,R)}\frac{Q_{ij,k}(x)Q_{ij,k}(x)}{2}+\frac{\tilde
f_B(Q(x))}{L}\,dx\nonumber\\=\frac{1}{R}\int_{\partial B(0,R)}
(Q_{ij,k}(x)\cdot x_k)^2\,dx+2\int_{B(0,R)}\frac{\tilde
f_B(Q(x))}{L}\,dx \label{mon2}
\end{eqnarray}

\par We have
\begin{eqnarray}\label{mon3}
\frac{\partial}{\partial
R}\left(\frac{1}{R}\int_{B(0,R)}\frac{Q_{ij,l}(x)Q_{ij,l}(x)}{2}+\frac{\tilde
f_B(Q(x))}{L}\,dx\right)
\nonumber=-\frac{1}{R^2}\int_{B(0,R)}\frac{Q_{ij,l}(x)\cdot
Q_{ij,l}(x)}{2}+\frac{\tilde
f_B(Q(x))}{L}\,dx\nonumber\\+\frac{1}{R}\int_{\partial
B(0,R)}\frac{Q_{ij,l}(x)\cdot Q_{ij,l}(x)}{2}+\frac{\tilde
f_B(Q(x))}{L}\,dx.
\end{eqnarray} The right-hand side of (\ref{mon3}) is positive from
(\ref{mon2}) and
hence the conclusion.
$\Box$

\medskip
\begin{lemma}( $W^{1,2}$ convergence to harmonic maps) Let
$\Omega\subset \Rr^3$ be a simply-connected bounded open set with
smooth boundary. Let $Q^{(L)}$ be a global minimizer of
$\tilde F_{LG}[Q]$ in the admissible class $\mathcal{A}_Q$
defined in (\ref{eq:am41}). Then there exists a sequence
$L_{k}\to 0$ so that $Q^{(L_{k})}\to Q^{(0)}$ strongly in
$W^{1,2}\left(\Omega;S_0\right)$, where $Q^{(0)}$ is the limiting
harmonic map defined in (\ref{eq:Q0}). \label{lemma:W^{1,2}}
\end{lemma}

\smallskip\par{\bf Proof.} Our proof follows closely, up to a point, the ideas of Proposition $1$ in \cite{bbh}.   Firstly, we note that the limiting
harmonic map $Q^{(0)}$ belongs to our admissible space $\Acal_{Q}$ and
since $Q^{(0)}(x)\in Q_{\min}$, a.e. $x\in\Omega$ (see Section ~\ref{sec:conventions}) we have that
$\tilde{f}_B\left(Q^{(0)}(x)\right) = 0$ a.e. $x\in\Omega$. Therefore

\begin{equation}
\int_\Omega \frac{1}{2}Q^{(L)}_{ij,k}(x)Q^{(L)}_{ij,k}(x)~dx\le
\int_\Omega\frac{1}{2}
Q^{(L)}_{ij,k}(x)Q^{(L)}_{ij,k}(x)+\frac{1}{L}\tilde f_B(Q^{(L)}(x))\,dx
\le \int_\Omega\frac{1}{2}Q^{(0)}_{ij,k}(x)Q^{(0)}_{ij,k}(x)\,dx
\label{upperbdh1}
\end{equation}

\par The $Q^{(L)}$'s are subject to the same boundary condition, $Q_b$, for
all $L$.
Therefore (\ref{upperbdh1}) shows that the $W^{1,2}$-norms of the
$Q^{(L)}$'s are bounded uniformly in $L$. Hence there exists a
weakly-convergent subsequence $Q^{(L_{k})}$ such that
$Q^{(L_{k})}\rightharpoonup Q^{(1)}$ in $W^{1,2}$, for some $Q^{(1)}\in
\Acal_Q$ as $L_k\to 0$. Using the lower semicontinuity of the
$W^{1,2}$ norm with respect to the weak convergence, we have that
\begin{equation}
\int_\Omega |\nabla Q^{(1)}(x)|^2~dx\le \int_\Omega |\nabla
Q^{(0)}(x)|^2~dx \label{q0q1}
\end{equation}

\par Relation (\ref{upperbdh1}) shows that
$\int_\Omega \tilde f_B(Q^{(L_k)}(x))\,dx\le L_k\int_\Omega
Q^{(0)}_{ij,k}(x)Q^{(0)}_{ij,k}(x)\,dx$ and
 hence
$\int\tilde f_B(Q^{(L_k)}(x))\,dx\to 0$ as $L_k\to 0$. Taking into
account that $\tilde f_B(Q)\ge 0, \forall Q\in S_0$ we have
that, on a subsequence  ${L_{k_j}}$,  $\tilde
f_B(Q^{(L_{k_j})}(x))\to 0$ for almost all $x\in\Omega$. From
Proposition~\ref{prop:bulk}, we know that $\tilde{f}_B(Q)=0$ if
and only if $Q\in Q_{\min}$ i.e. if $Q=s_+\left(n\otimes
n-\frac{1}{3}Id\right)$ for $n\in S^2$. On the other hand, the sequence $Q^{(L_k)}$ converges
weakly in $W^{1,2}$ and, on a subsequence, strongly in $L^2$ to
$Q^{(1)}$. Therefore, the weak limit $Q^{(1)}$ is of the form
\begin{equation}
Q^{(1)}(x)=s_+\left(n^{(1)}(x)\otimes
n^{(1)}(x)-\frac{1}{3}Id\right),\,n^{(1)}(x)\in\mathbb{S}^2,\,a.e.\,x\in\Omega
\label{eq:Q1}
\end{equation}

\par It was proved in \cite{bz} (see also \cite{bc}) that if $Q^{(1)}\in W^{1,2}$ and the domain $\Omega$ is simply-connected, we can assume, without loss of generality, that $n^{(1)}\in W^{1,2}(\Omega,\mathbb{S}^2)$ and its trace is $n_b$.
Then (\ref{eq:Q1}) implies $|\grad Q^{(1)}(x)|^2 = 2s_+^2|\grad n^{(1)}(x)|^2$ for a.e. $x\in\Omega$. Also, recalling the definition of $Q^{(0)}$ from Section ~\ref{sec:conventions} we have
 $|\grad Q^{(0)}(x)|^2 = 2s_+^2 |\grad n^{(0)}(x)|^2$  for a.e. $x\in\Omega$.
 
  Combining  (\ref{q0q1}) with (\ref{eq:n0}) and the observations in the previous paragraph,
we obtain
 $\int_\Omega |\nabla
n^{(1)}(x)|^2\,dx=\int_\Omega |\nabla n^{(0)}(x)|^2\,dx$ and
$\int_\Omega |\nabla Q^{(1)}(x)|^2\,dx= \int_\Omega |\nabla
Q^{(0)}(x)|^2\,dx$. Then:

$$\int_\Omega |\nabla Q^{(0)}(x)|^2\,dx\le
\liminf_{L_{k_j}\to 0}\int_\Omega |\nabla Q^{(L_{k_j})}(x)|~dx\le
\limsup_{L_{k_j}\to 0}\int_\Omega |\nabla Q^{(L_{k_j})}(x)|^2~dx\le
\int_\Omega |\nabla Q^{(0)}(x)|^2 ~dx,$$ which demonstrates that
$\lim_{L_{k_j}\to 0} \|\nabla Q^{(L_{k_j})}\|_{L^2}=\|\nabla
Q^{(0)}\|_{L^2}$. This together with the weak convergence
$Q^{(L_{k_j})}\to Q^{(0)}$ suffices to show the strong convergence
$Q^{(L_{k_j})}\to Q^{(0)}$ in $W^{1,2}$. $\Box$

\par The following has an elementary proof, that will be omitted:

\begin{lemma} The function $\tilde f_B:S_0\to \mathbb{R}_+$
is locally Lipschitz. \label{lemma:LLip}
\end{lemma}

\par We can now prove the uniform convergence of the
bulk energy density in the interior, away from the singularities of the limiting harmonic map $Q^{(0)}$.
\begin{proposition} Let $\Omega\subset\mathbb{R}^3$ be a simply-connected bounded
open set
with smooth boundary.
 Let $Q^{(L)}\in W^{1,2}(\Omega,S_0)$ denote a global minimizer of $\tilde F_{LG}[Q]$ in the admissible class $\Acal_Q$. Assume that we have a sequence $\{Q^{(L_k)}\}_{k\in\mathbb{N}}$ so that  $Q^{(L_k)}\to Q^{(0)}$ in $W^{1,2}(\Omega, S_0)$ as $L_k\to 0$.
\par  For any compact $K\subset \Omega$ such that $Q^{(0)}$ has no
singularity in $K$ we have

\begin{equation}
\lim_{L_k\to 0}\tilde f_B(Q^{(L_k)}(x))=0\quad x\in K
\end{equation} and the limit is uniform on $K$.
\label{uniformbulk}
\end{proposition}

\smallskip\par{\bf Proof.} {\it Lemma} ~\ref{lemma:W^{1,2}} shows that the strong limit $Q^{(0)}$ is a limiting harmonic map, as
defined in Section ~\ref{sec:conventions},
$Q^{(0)}=s_+(n^{(0)}(x)\otimes n^{(0)}(x)-\frac{1}{3}Id)$ where
$n^{(0)}(x)\in W^{1,2}(\Omega,\mathbb{S}^2)$ a  global energy
minimizer of the harmonic map problem, subject to the boundary
condition $n=n_b$ on $\partial\Omega$.
\par Let $\alpha_{L^k}=\tilde f_B(Q^{(L_k)}(x_0))$, for $x_0\in K$ an
arbitrary
point. {\it Proposition} ~\ref{prop:max} and {\it Lemma}
~\ref{lemma:LLip} imply that there exists a constant $\beta$
(independent of $x_0$)  so that

\begin{equation}
|\tilde f_B(Q^{(L)}(x))-\tilde f_B(Q^{(L)}(y))|\le \beta
|Q^{(L)}(x)-Q^{(L)}(y)|\label{beta}
\end{equation} for any $x,y\in\Omega,L>0$.
\par We  then have
\begin{eqnarray}
\alpha_{L^k}\le\tilde
f_B(Q^{(L_k)}(x))+\beta|Q^{(L_k)}(x)-Q^{(L_k)}(x_0)|\nonumber\\
\le\tilde f_B(Q^{(L_k)}(x))+\beta\|\nabla
Q^{(L_k)}\|_{L^\infty(K')}|x-x_0|\le \tilde
f_B(Q^{(L_k)}(x))+\frac{\tilde C\beta}{\sqrt{L_k}}|x-x_0|,\forall
x\in K'
\label{seqalphabound}
\end{eqnarray} where  $K'\subset \Omega$ is a compact  neighborhood of $K$
to be precisely defined  later. In the last relation above we use {\it Lemma}
$A.1$ from \cite{bbh} and the apriori bound given by {\it
Proposition} ~\ref{prop:max}. For reader's convenience we recall
that {\it Lemma} $A.1$ in \cite{bbh} states that if $u$ is a
scalar-valued function such that $-\Delta u=f$ on
$\Omega\subset\mathbb{R}^n$ then $|\nabla u(x)|^2\le
C\left(\|f\|_{L^\infty(\Omega)}\|u\|_{L^\infty(\Omega)}+\frac{1}{\textrm{dist}^2(x,\partial\Omega)}\|u\|^2_{L^\infty(\Omega)}\right)$
where $C$ is a constant that depends on $n$ only. In our case the
constant $\tilde C$ depends on the dimension, $n=3$, on
$a^2,b^2,c^2$ and on the distance $\sup_{y\in
K}d(y,\partial\Omega)$ only.
\par From (\ref{seqalphabound}) we have that
\begin{equation}
\alpha_{L^k}-\frac{\tilde C\beta\rho_k}{\sqrt{L_k}}\le \tilde
f_B(Q^{(L_k)}(x)),\forall x\in K',|x-x_0|<\rho_k
\end{equation}
 We argue similarly as in \cite{bbh} and divide by $L_k$ and
 integrate over $B_{\rho_k}(x_0)$ to obtain:

 \begin{equation}
 \frac{\rho_k^3}{L_k}(\alpha_{L_k}-\frac{\tilde C\beta\rho_k}{\sqrt{L_k}})\le
 \int_{B_{\rho_k}(x_0)}\frac{\tilde f_B(Q^{L_k}(x))}{L_k}\,dx
\label{half1}
 \end{equation}

\par Take an arbitrary $\varepsilon>0$. Recall that $K$ is a compact set that
does not contain singularities of $Q^{(0)}$. Then there exists a larger
compact set $K'$, so that $K\subset K'$, that does not contain
singularities either, and a constant $C_{K'}$ such that $|\nabla
Q^{(0)}(x)|^2<C_{K'}, \,\forall x\in K'$. For $R_0$ small enough,
with $R_0<\textrm{dist}(K,\partial\Omega)$ and such that
$B(x_0,R_0)\subset K',\forall x_0\in K$ we have
\begin{equation}
\frac{1}{R_0}\int_{B_{R_0}(x_0)}\frac{|\nabla
Q^{(0)}(x)|^2}{2}\,dx\le\frac{4\pi}{6} C_{K'} R_0^2\le
\frac{\varepsilon}{3}, \forall x_0\in K
\end{equation}

\par We fix an  $R_0$  as before. As $Q^{(L_k)}\to Q^{(0)}$ in
$W^{1,2}$, we have that there exists an $\bar L_0>0$ so that:

\begin{equation}
\frac{1}{R_0}\int_{B_{R_0}(x_0)}\frac{|\nabla
Q^{(L_k)}(x)|^2}{2}\,dx<\frac{1}{R_0}\int_{B_{R_0}(x_0)}\frac{|\nabla
Q^{(0)}(x)|^2}{2}\,dx+\frac{\varepsilon}{3}, \textrm{ for }L_k<\bar
L_0, \forall x_0\in K
\end{equation}

\par The arguments in \cite{bbh}  fail to work in our case as we have
a three
 dimensional domain, unlike in the quoted paper, where the
domain is two dimensional. In our case, using the monotonicity
 formula from {\it Lemma} $2$ and taking $\rho_k<R_0$ we obtain:

\begin{equation}
\int_{B_{\rho_k}(x_0)}\frac{\tilde f_B(Q^{(L_k)}(x))}{L_k}\,dx\le
\frac{\rho_k}{R_0}\int_{B_{R_0}(x_0)}\frac{ |\nabla Q^{(L_k)}(x)|^2}{2}
+\frac{\tilde f_B(Q^{(L_k)}(x))}{L_k}\,dx\le
\rho_k\left(\frac{2\varepsilon}{3}+\frac{\varepsilon}{3}\right)
\label{half2}
\end{equation} for $L_k<\bar L_1$ with $\bar L_1$ small enough so that
$\frac{1}{R_0}\int_{B_{R_0}(x_0)} \frac{\tilde f_B(Q^{(L_k)}(x))}{L_k}\,dx<
\frac{\varepsilon}{3}$ (note that there exists such an $\bar L_1$
as  the proof of {\it Lemma} $3$ shows that $\int_\Omega
\frac{\tilde f_B(Q^{(L_k)}(x))}{L_k}\,dx=o(1)$ as $L_k\to 0$).

\par We take $\rho_k=\frac{\alpha_{L_k}\sqrt{L_k}}{2\tilde C\beta}$ . Then, from  (\ref{half1}) and (\ref{half2}) we obtain

$$\alpha_{L_k}^3<8(\tilde C\beta)^2\varepsilon$$ for $L_k<\min\{\bar
L_0,\bar L_1\}$. As $\varepsilon>0$  is arbitrary and the
estimate on $\alpha_{L_k}=\tilde f_B(Q^{(L_k)}(x_0)),\,x_0\in K$
is obtained in a manner independent of $x_0$, we have the claimed
result. $\Box$

\par We also need the following

\smallskip
\begin{lemma} There exists $\varepsilon_0>0$  so that:
\begin{eqnarray}
\frac{1}{\tilde C}\tilde f_B(Q)\le \sum_{i,j=1}^3
\left(\frac{\partial \tilde f_B(Q)}{\partial Q_{ij}}
+b^2\frac{\delta_{ij}}{3}\textrm{tr}(Q^2)\right)^2\le\tilde
C\tilde
f_B(Q)\nonumber\\
\forall Q\in S_0\textrm{ such that }|Q-s_+(n\otimes n-\frac{1}{3}Id)|\le
\varepsilon_0, \textrm{ for some }n\in\mathbb{S}^2 \label{2prime}
\end{eqnarray} where $s_+=\frac{b^2+\sqrt{b^4+24a^2c^2}}{4c^2}$ and the
constant
$\tilde C$ is independent of $Q$, but depends on
$a^2,b^2,c^2$. \label{lemma:comparison}
\end{lemma}

\smallskip\par{\bf Proof.} Recall from Proposition \ref{prop:bulk}, \cite{bm} that  $\tilde f_B(Q)\ge 0$ and $\tilde
f_B(Q)=0\leftrightarrow Q=s_+(n\otimes n-\frac{1}{3}Id)$ with
 $s_+=\frac{b^2+\sqrt{b^4+24a^2c^2}}{4c^2}$ and $n\in\mathbb{S}^2$.
 \par Let the eigenvalues of $Q$ be $x,y,-x-y$. We define
$F(x,y)\stackrel{def}{=}-a^2(x^2+y^2+xy)+b^2xy(x+y)+c^2(x^2+y^2+xy)^2$
and $D\stackrel{def}{=}\min_{(x,y)\in\mathbb{R}^2}F(x,y)$. Then
$\tilde F(x,y)\stackrel{def}{=}F(x,y)-D=\tilde f_B(Q)$.
\par Then $\tilde F=0$ only at three pairs $(x,y)$ namely
$(-\frac{s_+}{3},-\frac{s_+}{3}),(-\frac{s_+}{3},2\frac{s_+}{3})$
and $(2\frac{s_+}{3},-\frac{s_+}{3})$.
\par   On the other hand we have
\begin{eqnarray}
\sum_{i,j=1}^3 \left(\frac{\partial\tilde f_B}{\partial
Q_{ij}}+\frac{b^2\delta_{ij}}{3}\textrm{tr}(Q^2)\right)^2=a^4\textrm{tr}(Q^2)+(\frac{b^4}{6}-2a^2c^2)(\textrm{tr}(Q^2))^2\nonumber\\+
c^4(\textrm{tr}(Q^2))^3+2a^2b^2\textrm{tr}(Q^3)-2b^2c^2\textrm{tr}(Q^2)\textrm{tr}(Q^3)
\end{eqnarray} (where we used the identity
$\textrm{tr}(Q^4)=\frac{(\textrm{tr}(Q)^2)^2}{2}$, valid for a traceless
symmetric
$3\times 3$ matrix)
\par If we denote $h(Q)= \sum_{i,j=1}^3 \left(\frac{\partial \tilde
f_B(Q)}{\partial Q_{ij}}+b^2\frac{\delta_{ij}}{3}\textrm{tr}(Q^2)\right)^2$
we have $h(Q)=H(x,y)$ where $H:\mathbb{R}^2\to\mathbb{R}$ is given by

\begin{eqnarray}
H(x,y)\stackrel{def}{=}2a^4(x^2+y^2+xy)+4(\frac{b^4}{6}-2a^2c^2)(x^2+y^2+xy)^2+8c^4(x^2+y^2+xy)^3
\nonumber\\+12b^2c^2xy(x+y)(x^2+y^2+xy)-6a^2b^2xy(x+y)\nonumber
\end{eqnarray}

\par We claim that there exist
$\varepsilon_1,\varepsilon_2,\varepsilon_3>0$ so that

\begin{equation}
\frac{1}{\tilde C}\tilde F(x,y)\le H(x,y)\le\tilde C\tilde
F(x,y),\forall (x,y)\in
B_{\varepsilon_1}(-\frac{s_+}{3},-\frac{s_+}{3}),B_{\varepsilon_2}(-\frac{s_+}{3},2\frac{s_+}{3}),B_{\varepsilon_3}
(2\frac{s_+}{3},-\frac{s_+}{3}) \label{singularbounds}
\end{equation} which gives the conclusion.
\par We prove the inequality (\ref{singularbounds}) only for $(x,y)\in
B_{\varepsilon_1}(-\frac{s_+}{3},-\frac{s_+}{3})$; the other two
cases can be dealt with similarly.

\par Careful computations show:

\begin{displaymath}
H(-\frac{s_+}{3},-\frac{s_+}{3})=\frac{\partial H}{\partial
x}(-\frac{s_+}{3},-\frac{s_+}{3})=\frac{\partial H}{\partial
y}(-\frac{s_+}{3},-\frac{s_+}{3})=0
\end{displaymath}

\begin{displaymath}
\frac{\partial^2 H}{\partial
y^2}(-\frac{s_+}{3},-\frac{s_+}{3})=\frac{\partial^2 H}{\partial
x^2}(-\frac{s_+}{3},-\frac{s_+}{3})=4(b^4+6a^2c^2)\frac{b^4+12a^2c^2+b^2\sqrt{b^4+2a^2c^2}}{24c^4}
\end{displaymath}

\begin{displaymath}
\frac{\partial ^2 H}{\partial x\partial
y}(-\frac{s_+}{3},-\frac{s_+}{3})=-2(b^4-12a^2c^2)\frac{b^4+12a^2c^2+b^2\sqrt{b^4+24a^2c^2}}{24c^4}
\end{displaymath}

\begin{displaymath}
\tilde F(-\frac{s_+}{3},-\frac{s_+}{3})=\frac{\partial \tilde
F}{\partial x}(-\frac{s_+}{3},-\frac{s_+}{3})=\frac{\partial
\tilde F}{\partial y}(-\frac{s_+}{3},-\frac{s_+}{3})=0
\end{displaymath}

\begin{displaymath}
\frac{\partial^2 \tilde F}{\partial
y^2}(-\frac{s_+}{3},-\frac{s_+}{3})=\frac{\partial^2 \tilde
F}{\partial
x^2}(-\frac{s_+}{3},-\frac{s_+}{3})=\frac{1}{4c^2}(b^4+12a^2c^2+
b^2\sqrt{b^4+24a^2c^2})
\end{displaymath}

\begin{displaymath}
\frac{\partial ^2 \tilde F}{\partial x\partial
y}(-\frac{s_+}{3},-\frac{s_+}{3})=3a^2
\end{displaymath}

\par Let $(x_0,y_0)=(-\frac{s_+}{3},-\frac{s_+}{3})$. We have

\begin{displaymath}
\frac{H(x,y)}{\tilde F(x,y)}=\frac{H_1(x,y)+R_H(x,y)}{\tilde
F_1(x,y)+R_{\tilde F}(x,y)}
\end{displaymath} where $H_1(x,y)=(x-x_0)^2\frac{\partial^2 H}{\partial
x^2}(x_0,y_0)+2(x-x_0)(y-y_0)\frac{\partial^2 H}{\partial x\partial
y}(x_0,y_0)+
(y-y_0)^2 \frac{\partial^2 H}{\partial y^2}(x_0,y_0)$  and $\tilde
F_1(x,y)=(x-x_0)^2\frac{\partial^2 \tilde F}{\partial
x^2}(x_0,y_0)+2(x-x_0)(y-y_0)\frac{\partial^2 \tilde F}{\partial
x\partial y}(x_0,y_0)+ (y-y_0)^2 \frac{\partial^2 \tilde
F}{\partial y^2}(x_0,y_0)$ with $R_H$, $R_{\tilde F}$ the
remainders in the Taylor expansions around $(x_0,y_0)$.

\par From the definition of Taylor expansions, we have that there exists
$\varepsilon_0>0$ so that on $B_{\varepsilon_1}(x_0,y_0)$ we have
\begin{equation}
|R_H(x,y)|\le\frac{1}{2} H_1(x,y)\textrm{     and }|R_{\tilde
F}(x,y)|\le\frac{1}{2} \tilde F_1(x,y), \forall (x,y)\in
B_{\varepsilon_1}(x_0,y_0) \label{rhfh}
\end{equation}
\par  On the other hand we have
\begin{equation}
\tilde F_1(x,y)\frac{1}{8(b^4+6a^2c^2)}\le H_1(x,y)\le \tilde
F_1(x,y)8(b^4+6a^2c^2)\,\forall (x,y)\in\mathbb{R}^2 \label{f1h1}
\end{equation} hence, combining (\ref{rhfh}) and (\ref{f1h1}), we get:

\begin{equation}
\tilde F(x,y)\frac{1}{24(b^4+6a^2c^2)}\le H(x,y)\le \tilde
F(x,y)24(b^4+6a^2c^2),\forall (x,y)\in
B_{\varepsilon_1}(-\frac{s_+}{3},-\frac{s_+}{3})
\end{equation} which yields  claim (\ref{singularbounds}) for
$(x,y)\in B_{\varepsilon_1}(-\frac{s_+}{3},-\frac{s_+}{3})$. The other two cases can be analyzed analogously. $\Box$

\smallskip

We continue by proving a  Bochner-type inequality that is crucial for the derivation of uniform (in $L$) Lipschitz bounds, away from the singularities  of the limiting harmonic map.
This type of inequalities were first used (to the best of our knowledge) in the context of 
harmonic maps (see \cite{schoen} and the references there) and later
adapted to other, more complicated contexts (see for instance \cite{remarks}). The main difficulty in the proof 
of Proposition ~\ref{prop:unifconvinterior} (to follow) is the derivation of the next lemma.

\begin{lemma}  There exists $\varepsilon_0>0$ and  a constant $C>0$,
 independent of $L$, so that for $Q^{(L)}$ a global minimizer of $\tilde F_{LG}[Q]$ in the admissible space $\Acal_Q$, we have

\begin{equation}
-\Delta e_L(Q^{(L)})(x)\le Ce^2_L(Q^{(L)}(x)) \label{bochner}
\end{equation} provided there exists a ball $B_{\rho(x)}(x)$ for some $\rho(x)>0$ such that  
$|Q^{(L)}(y)-s_+\left(m(y)\otimes
m(y)-\frac{1}{3}Id\right)|<\varepsilon_0$ with $m(y)\in\mathbb{S}^2$, for all $y\in B_{\rho(x)}(x)$.
\label{lemma:bochner}
\end{lemma}

\smallskip\par{\bf Proof.}  In the following we drop the superscript $L$ for convenience. We have:

\begin{eqnarray}
-\Delta\left(\frac{Q_{ij,k}Q_{ij,k}}{2}\right)=-\Delta(Q_{ij,k})Q_{ij,k}-Q_{ij,kl}Q_{ij,kl}\le\nonumber\\
\le-\partial_k\left[\frac{1}{L}\frac{\partial\tilde f_B}{\partial
Q_{ij}}(Q(x))+\frac{b^2\delta_{ij}}{3L}\textrm{tr}(Q^2)\right]Q_{ij,k}=-
\partial_k\left[\frac{1}{L}\frac{\partial\tilde f_B}{\partial Q_{ij}}(Q(x))\right]Q_{ij,k}
\label{bochp1}
\end{eqnarray}

\par On the other hand:
\begin{eqnarray}
-\Delta\left[\frac{1}{L}\tilde
f_B(Q(x)\right]=-\partial_k\left(\frac{1}{L}\frac{\partial\tilde
f_B}{\partial Q_{ij}}(Q(x))\partial_k Q_{ij}\right)=
-\partial_k\left[\left[\frac{1}{L}\frac{\partial\tilde f_B}{\partial
Q_{ij}}(Q(x))+
\frac{b^2\delta_{ij}}{3L}\textrm{tr}(Q^2)\right]\partial_k Q_{ij}\right]\nonumber\\
=-\underbrace{\left(\frac{1}{L}\frac{\partial\tilde f_B}{\partial
Q_{ij}}(Q(x))+\frac{b^2\delta_{ij}}{3L}\textrm{tr}(Q^2)\right)}_{\stackrel{def}{=}Z}\times\underbrace{\Delta
Q_{ij}}_{=Z}-\partial_k\left(\frac{1}{L}\frac{\partial\tilde
f_B}{\partial Q_{ij}}(Q(x))\right)Q_{ij,k}\nonumber\\
\le -\partial_k\left(\frac{1}{L}\frac{\partial\tilde f_B}{\partial
Q_{ij}}(Q(x))\right)Q_{ij,k} \label{bochp2}
\end{eqnarray}

\par We take $\varepsilon_1>0$  a small number, to be made precise later.
For any such $\varepsilon_1$ we can pick  $\varepsilon_0>0$ small
enough so that if the eigenvalues of $Q(x)$ are
$(\lambda,\mu,-\lambda-\mu)$ then one of the three numbers
$(\lambda+\frac{s_+}{3})^2+(\mu+\frac{s_+}{3})^2+(\lambda+\mu+2\frac{s_+}{3})^2$,
$(\lambda+\frac{s_+}{3})^2+(\mu-2\frac{s_+}{3})^2+(\lambda+\mu-\frac{s_+}{3})^2$,
$(\lambda-2\frac{s_+}{3})^2+(\mu+\frac{s_+}{3})^2+(\lambda+\mu-\frac{s_+}{3})^2$
is less than or equal to $\varepsilon_1$ (this can be done because the
eigenvalues are continuous functions of matrices, \cite{kato}, and
the matrix $s_+(n\otimes n-\frac{1}{3}Id)$ has eigenvalues
$-\frac{s_+}{3},-\frac{s_+}{3}$ and $2\frac{s_+}{3}$).  Note moreover
that we need to choose $\varepsilon_0$ to be smaller than the
choice (of $\varepsilon_0$) in {\it Lemma} ~\ref{lemma:comparison}
as we will need to use that lemma in the remainder of this proof.

\par For the matrix $Q(x)$, let us denote its eigenvectors by $n_1(x),n_2(x),n_3(x)$  and let
$\lambda_1(x),\lambda_2(x)$ ,
$\lambda_3(x)=-\lambda_1(x)-\lambda_2(x)$ denote the corresponding
eigenvalues. From the preceeding discussion,  we can, without loss
of generality, assume that
\begin{equation}
(\lambda_1+\frac{s_+}{3})^2+(\lambda_2+\frac{s_+}{3})^2+(\lambda_1+\lambda_2+2\frac{s_+}{3})^2<\varepsilon_1
\label{eigenvaluesclose}
\end{equation}

\par We define the matrix
$$Q^x\stackrel{def}{=}-\frac{s_+}{3} n_1(x)\otimes
n_1(x)-\frac{s_+}{3} n_2(x)\otimes n_2(x)+\frac{2s_+}{3}
n_3(x)\otimes n_3(x)$$ (Note that  there exists a
$m(x)\in\mathbb{S}^2$ so that $Q^x=s_+(m(x)\otimes
m(x)-\frac{1}{3}Id)$). 

Taking into account
(\ref{eigenvaluesclose}) and the fact that $Q(x)$ and $Q^x$ have the same
eigenvectors,   we have :
\begin{equation}
\textrm{tr}(Q(x)-Q^x)^2=(\lambda_1+\frac{s_+}{3})^2+(\lambda_2+\frac{s_+}{3})^2+(\lambda_1+\lambda_2+2\frac{s_+}{3})^2<\varepsilon_1
 \label{diagtrick}
\end{equation}

\par Using the of  Taylor expansion of
$\frac{1}{2}\frac{\partial^2\tilde f_B}{\partial Q_{ij}\partial
Q_{mn}}(Q(x))$ around $Q^x$ we obtain:

\begin{equation}
\frac{1}{2}\frac{\partial^2\tilde f_B}{\partial Q_{ij}\partial
Q_{mn}}(Q(x))=\frac{1}{2}\frac{\partial^2\tilde f_B}{\partial
Q_{ij}\partial Q_{mn}}(Q^x)+ \frac{1}{2}\frac{\partial^3\tilde
f_B}{\partial Q_{ij}\partial Q_{mn}\partial
Q_{pq}}(Q^x)(Q_{pq}(x)-Q^x_{pq})+\mathcal{R}^{ijmn}(Q^x,Q(x))
\label{taylor2}
\end{equation} where $\mathcal{R}^{ijmn}(Q^x,Q(x))$ is the remainder.

\par From (\ref{taylor2}) we have:

\begin{eqnarray}
-\partial_k\left(\frac{1}{L}\frac{\partial \tilde f_B}{\partial
Q_{ij}}(Q(x)\right)Q_{ij,k}=-\frac{1}{L}\frac{\partial^2 \tilde f_B}{\partial
Q_{ij}\partial Q_{mn}}Q_{mn,k}Q_{ij,k}=\nonumber\\
=\underbrace{-\frac{1}{L}\frac{\partial^2 \tilde f_B}{\partial
Q_{ij}\partial Q_{mn}}(Q^x)Q_{mn,k}Q_{ij,k}}_{\le 0}-\frac{1}{L}
\frac{\partial^3\tilde f_B}{\partial Q_{ij}\partial Q_{mn}\partial
Q_{pq}}(Q^x)(Q_{pq}(x)-Q^x_{pq})Q_{ij,k}Q_{mn,k}-\nonumber\\-\frac{1}{L}\mathcal{R}^{ijmn}(Q(x),Q^x)Q_{ij,k}Q_{mn,k}\le\nonumber\\
\le\frac{C_0\delta}{L^2}\Sigma_{i,j,m,n=1}^3\left(\frac{\partial^3\tilde
f_B}{\partial Q_{ij}\partial Q_{mn}\partial
Q_{pq}}(Q^x)\right)^2(Q_{pq}(x)-Q^x_{pq})^2\nonumber\\+\frac{C_0\delta}
{L^2}\sum_{i,j,m,n=1}^3\left(\mathcal{R}^{ijmn}\right)^2(Q(x),Q^x)+\frac{1}{\delta}|\nabla
Q|^4\le\nonumber\\
\le \frac{\delta}{L^2}\sum_{i,j,m,n=1}^3\left[\bar
C_0\left(\frac{\partial^3 f}{\partial Q_{ij}\partial
Q_{mn}\partial Q_{pq}}(Q^x)\right)^2+1\right](Q_{pq}(x)-
Q_{pq}^x)^2\frac{1}{\delta}|\nabla Q|^4\nonumber\\
\le \frac{C_1\delta}{L^2}\textrm{tr}(Q(x)-Q^x)^2+\frac{1}{\delta}|\nabla
Q|^4 \label{estimateboch}
\end{eqnarray} where $0<\delta<1$ and  $C_0,\bar C_0,C_1$ are independent
of $L$ and
$x$. For the first term in the second line  above we use the fact that
the Hessian
matrix of a function $\tilde f_B(Q)$ is non-negative definite at a global
minimum
(which holds true in our case as well, as one can easily check,  even
though we have $\tilde f_B(Q)$ restricted to the linear space $S_0$).

\par Let us recall (from the proof of the previous lemma) the definitions of $F$ and $\tilde F$. Then, for  a matrix $Q\in S_0$ with eigenvalues 
$(\lambda_1,\lambda_2,-\lambda_1-\lambda_2)$ we have

\begin{equation}
\tilde f_B(Q)=\tilde F(\lambda_1,\lambda_2) \label{spectralred}
\end{equation}

\par We claim that for $\varepsilon_1>0$ small enough there exists $C_2$
independent
of $L, \lambda_1,\lambda_2$ so that

\begin{eqnarray}
C_2\left((\lambda_1+\frac{s_+}{3})^2+(\lambda_2+\frac{s_+}{3})^2+(\lambda_1+\lambda_2+2\frac{s_+}{3})^2\right)\le
\tilde F(\lambda_1,\lambda_2)\,\nonumber\\
\textrm{ for all } (\lambda_1,\lambda_2) \,\textrm{ so
that}(\lambda_1+\frac{s_+}{3})^2+(\lambda_2+\frac{s_+}{3})^2+(\lambda_1+\lambda_2+2\frac{s_+}{3})^2<\varepsilon_1.
\label{spectrallowerbd}
\end{eqnarray}

\par Careful computations show:

\begin{displaymath}
\tilde F(-\frac{s_+}{3},-\frac{s_+}{3})=\frac{\partial
\tilde F}{\partial
\lambda_1}(-\frac{s_+}{3},-\frac{s_+}{3})=\frac{\partial
\tilde F}{\partial \lambda_2}(-\frac{s_+}{3},-\frac{s_+}{3})=0
\end{displaymath}

\begin{displaymath}
\frac{\partial^2 \tilde F}{\partial
\lambda_2^2}(-\frac{s_+}{3},-\frac{s_+}{3})=\frac{\partial^2
\tilde F}{\partial
\lambda_1^2}(-\frac{s_+}{3},-\frac{s_+}{3})=\frac{1}{4c^2}(b^4+12a^2c^2+
b^2\sqrt{b^4+24a^2c^2})
\end{displaymath}

\begin{displaymath}
\frac{\partial ^2 \tilde F}{\partial \lambda_1\partial
\lambda_2}(-\frac{s_+}{3},-\frac{s_+}{3})=3a^2
\end{displaymath}

\par Using a Taylor expansion around
$(\lambda_1,\lambda_2)=(-\frac{s_+}{3},-\frac{s_+}{3})$ we have

\begin{eqnarray}
\tilde F(\lambda_1,\lambda_2)=\frac{1}{8c^2}\left(b^4+12a^2c^2+b^2\sqrt{b^4+24a^2c^2}\right)[(\lambda_1+\frac{s_+}{3})^2+
(\lambda_2+\frac{s_+}{3})^2]+\nonumber\\
+3a^2(\lambda_1+\frac{s_+}{3})(\lambda_2+\frac{s_+}{3})+R(\lambda_1,\lambda_2)\ge\nonumber\\
\ge
\frac{1}{2}\left\{\frac{1}{8c^2}\left(b^4+12a^2c^2+b^2\sqrt{b^4+24a^2c^2}\right)[(\lambda_1+\frac{s_+}{3})^2+
(\lambda_2+\frac{s_+}{3})^2]+3a^2(\lambda_1+\frac{s_+}{3})(\lambda_2+\frac{s_+}{3})\right\}
\label{comb1}
\end{eqnarray} where $R(\lambda_1,\lambda_2)$ is the remainder in the Taylor
expansion, and the inequality holds {\it provided} that the remainder $R$ is small enough. We choose 
$\varepsilon_1>0$ to be small enough so that if 
$(\lambda_1+\frac{s_+}{3})^2+(\lambda_2+\frac{s_+}{3})^2+(\lambda_1+\lambda_2+2\frac{s_+}{3})^2<\varepsilon_1$
then $R$ is small enough and the inequality  above holds.

\par As the quadratic form
$\frac{1}{16c^2}\left(b^4+12a^2c^2+b^2\sqrt{b^4+24a^2c^2}\right)[(\lambda_1+\frac{s_+}{3})^2+
(\lambda_2+\frac{s_+}{3})^2]+\frac{3}{2}a^2(\lambda_1+\frac{s_+}{3})(\lambda_2+\frac{s_+}{3})$
is positive definite, there exists a $C_2>0$, depending only on
$a^2$,$b^2$ and $c^2$ such that

\begin{eqnarray}
\frac{1}{2}\left\{\frac{1}{8c^2}\left(b^4+12a^2c^2+b^2\sqrt{b^4+24a^2c^2}\right)[(\lambda_1+\frac{s_+}{3})^2+
(\lambda_2+\frac{s_+}{3})^2]+3a^2(\lambda_1+\frac{s_+}{3})(\lambda_2+\frac{s_+}{3})\right\}\ge\nonumber\\
\ge
C_2\left((\lambda_1+\frac{s_+}{3})^2+(\lambda_2+\frac{s_+}{3})^2+(\lambda_1+\lambda_2+\frac{2s_+}{3})^2\right)\,\forall
(\lambda_1,\lambda_2)\in\mathbb{R}^2\nonumber
\end{eqnarray}
\par Combining this last inequality with (\ref{comb1}) we obtain the claim
(\ref{spectrallowerbd}).
\par The relation (\ref{spectrallowerbd}) together with
(\ref{spectralred}) and
(\ref{diagtrick}) show that
$\textrm{tr}(Q(x)-Q^x)^2\le C_3 \tilde f_B(Q(x))$  for some $C_3$
independent of $L$ and $x$, which combined with
(\ref{estimateboch}) shows
$$-\partial_k\left(\frac{1}{L}\frac{\partial \tilde f_B}{\partial
Q_{ij}}(Q(x)\right)Q_{ij,k}
\le \frac{\delta C_4}{L^2}\tilde f_B(Q(x))+\frac{1}{\delta}|\nabla
Q(x)|^4$$ with $C_4$ a constant independent of $L$ and $x$ and any
$\delta>0$. This last inequality together with (\ref{bochp1}) and
(\ref{bochp2}) show:

$$-\Delta e_L+\frac{1}{L^2}\sum_{i,j=1}^3
\left(\frac{\partial\tilde f_B}{\partial
Q_{ij}}+\frac{b^2\delta_{ij}}{3}\textrm{tr}(Q^2) \right)^2 \le
\frac{\delta C_4}{L^2}\tilde f_B(Q)+\frac{1}{\delta}|\nabla Q|^4$$

\par Taking into account {\it Lemma}~\ref{lemma:comparison} and choosing
$\delta$
small enough (depending only on $C_4$ and the constant $\tilde C$
from {\it Lemma}~\ref{lemma:comparison}) we can absorb the term
$\frac{\delta C_4}{L^2}\tilde f_B(Q)$ on the right hand side into
the left hand side and obtain

$$-\Delta e_L\le \frac{1}{\delta}|\nabla Q|^4,$$ giving  the desired conclusion.
$\Box$

\begin{lemma} Let $\Omega\subset\mathbb{R}^3$ be a simply-connected bounded open set
with smooth boundary.
 Let $Q^{(L_k)}\in W^{1,2}(\Omega,S_0)$ be a
sequence of global minimizers for the energy $\tilde F_{LG}[Q]$ in the admissible space $\Acal_Q$. Assume that as $L_k\to 0$
we have $Q^{(L_k)}\to Q^{(0)}$ in $W^{1,2}(\Omega,S_0)$.
\par Let $K\subset\Omega$ be a compact set which contains no singularity
of $Q^{(0)}$.
There exists $C_1>0,C_2>0,\bar L_0>0$( all constants independent
of $L_k$) so that if for
$a\in K,0<r<d(a,\partial K)$ we
have
$$\frac{1}{r}\int_{B_r(a)}e_{L_k}(Q^{(L_k})(x))\,dx\le C_1$$
then
$$r^2\sup_{B_{\frac{r}{2}}(a)}e_{L_k}(Q^{(L_k)})\le C_2.$$ for all
$L_k<\bar L_0$.
\label{lemma:uniformreg}
\end{lemma}

\smallskip\par{\bf Proof.} Taking into account our assumptions on the
sequence
$\left(Q^{(L_k)}\right)_{k\in\mathbb{N}}$, {\it Proposition}
~\ref{uniformbulk} shows that for any given $\tilde \varepsilon_0$
smaller than $\varepsilon_0$ in {\it Lemma}~\ref{lemma:comparison}
and also smaller than the $\varepsilon_0$ in {\it
Lemma}~\ref{lemma:bochner}, we have that there exists a $\bar L_0$
so that for $L_k<\bar L_0$ we have

\begin{equation}
\|Q^{(L_k)}(x)-s_+\left(n(x)\otimes n(x)-\frac{1}{3}Id\right)\|\le
\tilde\varepsilon_0,\,\forall x\in K,\, \textrm{for some
}n(x)\in\mathbb{S}^2 \label{closeQ}
\end{equation}

\par  We continue reasoning similarly as in \cite{schoen}. We fix an
arbitrary
$L_k<\bar L_0$ and
an $a\in\Omega$ and take a $r>0$ so that
$0<r<\min\{d(a,\partial\Omega),d(a,K)\}$ . We let $r_1>0$ and
$x_1\in B_{r_1}(a)$ be such that

$$\max_{0\le s\le \frac{2}{3}r}(\frac{2}{3} r-s)^2\max_{B_s(a)}
e_{L_k}(Q^{(L_k)})=(\frac{2}{3} r-r_1)^2
\max_{B_{r_1}(a)}e_{L_k}(Q^{(L_k)})=(\frac{2}{3}r-r_1)^2e_{L_k}(Q^{(L_k)}(x_1))$$
\par Define
$e_1^{(L_k)}\stackrel{def}{=}\max_{B_{r_1}(a)}e_{L_k}(Q^{(L_k}))$. Then:

\begin{eqnarray}
\max_{B_{\frac{2/3\cdot r-r_1}{2}}(x_1)}e_{L_k}(Q^{(L_k)})\le
\max_{B_{\frac{2/3\cdot r+r_1}{2}}(a)}e_{L_k}(Q^{(L_k)})\nonumber\\
\le \frac{(2/3\cdot
r-r_1)^2\max_{B_{r_1}(a)}e_{L_k}(Q^{(L_k)})}{(2/3\cdot r-(2/3\cdot
r+r_1)/2)^2}= 4\max_{B_{r_1}(a)}e_{L_k}(Q^{(L_k)})=4e_1^{(L_k)}
\label{rel57}
\end{eqnarray} where for the first inequality we use the fact that
$B_{\frac{(2/3r-r_1)}{2}}(x_1)\subset B_{\frac{(2/3\cdot
r+r_1)}{2}}(a)$ and for the second inequality, we use the
definition of $r_1$.
\par Let $r_2=\frac{(2/3\cdot r-r_1)\sqrt{e_1^{(L_k)}}}{2}$ and define
$R^{(L_k)}(x)=Q^{(L_k)}\left(x_1+\frac{x}{\sqrt{e_1^{(L_k)}}}\right)$. We let
$\bar L_k=e_1^{(L_k)}L_k$ and then 

\begin{eqnarray}e_{\bar L_k}(R^{(L_k)})=\frac{1}{2}|\nabla
R^{(L_k)}|^2+\frac{\tilde
f_B(R^{(L_k)})} {\bar L_k}=\frac{1}{2} \frac{|\nabla
Q^{(L_k)}|^2}{e_1^{(L_k)}}+\frac{\tilde f_B(Q^{(L_k)})}
{e_1^{(L_k)} L_k}=\frac{1}{e_1^{(L_k)}}e_{L_k}(Q^{(L_k)})\nonumber
\end{eqnarray}

Equation (\ref{rel57}) then implies  $$\max_{x\in B_{r_2}(0)}e_{\bar
L_k}(R^{(L_k)})=\max_{x\in B_{\frac{(2/3r-r_1)}{2}}(x_1)}
\frac{e_{L_k}(Q^{(L_k)}(x))} {e_1^{L_k}}\le 4$$

where  the equality above follows from the definition of $r_2$ and  $R^{(L_k)}$
and  the inequality  above follows from  equation
(\ref{rel57}). Thus, we have 

\begin{equation}
\max_{B_{r_2(0)}}e_{\bar L_k}(R^{(L_k)})\le 4,\,e_{\bar
L_k}(R^{(L_k)})(0)=1 \label{restr}
\end{equation} where $R^{(L_k)}$ satisfies the following system of elliptic PDEs

\begin{equation}
\bar L_k
R^{(L_k)}_{ij,kk}=-a^2R^{(L_k)}_{ij}-b^2\left(R^{(L_k)}_{ik}R^{(L_k)}_{kj}-\frac{\delta_{ij}}{3}\textrm{tr}((R^{(L_k)})^2)\right)+c^2R^{(L_k)}_{ij}\textrm{tr}((R^{(L_k)})^2)
\label{rescaledeq}
\end{equation}

\par We now  claim that
\begin{equation}r_2\le 1
\label{r2claim}
\end{equation}

\par It is clear  that $r_2\le 1$ implies the conclusion. Let us assume for
contradiction  that $r_2>1$. Then we claim that there
exists a  constant $C>0$ , independent of $L_k$, so that

\begin{equation}
1\le C\int_{B_1} e_{\bar L_k}(R^{(L_k)})(x)\,dx \label{lowerbd}
\end{equation}

\par The matrix $R^{(L_k)}$ satisfies the system
(\ref{rescaledeq}) (which
is the rescaled version of
(\ref{ELeq}) ); using  relation (\ref{closeQ}) and the
definition of $R^{(L_k)}$ as well as the fact that $r_2>1$, we can
apply {\it Lemma} ~\ref{lemma:bochner} to $e_{\bar
L_k}(R^{(L_k)})$ and obtain $$-\Delta e_{\bar L_k}(R^{(L_k)}(x))\le
Ce^2_{\bar L_k}(R^{(L_k)}(x)) \stackrel{(\ref{restr})}{\le}
4Ce_{\bar L_k}(R^{(L_k)}(x)),\forall x\in B_1(0)$$

\par Combining (\ref{restr}) and the Harnack inequality (see for instance
\cite{taylor}, Ch.$14$, Thm. $9.3$) along with
the above relation we obtain (\ref{lowerbd}).

\smallskip\par  We have
\begin{eqnarray}
\int_{B_1}e_{\bar L_k}(R^{(L_k)}(x))\,dx\le
\frac{1}{r_2}\int_{B_{r_2}(0)}\frac{|\nabla
R^{(L_k)}(x)|^2}{2}+\frac{\tilde f_B(R^{(L_k)}(x))}
{L_ke_1^{(L_k)}}\,dx=\nonumber\\
=\frac{2}{2/3\cdot r-r_1}\int_{B_{(2/3\cdot
r-r_1)/2}(x_1)}e_{L_k}(Q^{(L_k)}(x))\,dx\le
\frac{3}{r}\int_{B_{r/3}(x_1)}e_{L_k}(Q^{(L_k)}(x))\,dx\nonumber\\
\le \frac{3}{r}\int_{B_r(a)}e_{L_k}(Q^{(L_k)}(x))\,dx\le 3C_1
\label{rel62}
\end{eqnarray} where for the first inequality we use the monotonicity
inequality (Lemma
~\ref{lemma:mon1}) and the assumption that $r_2\ge 1$ (note that
the equation satisfied by $R^{(L_k)}$ , equation
(\ref{rescaledeq}) is the same as the equation satisfied by
$Q^{(L_k)}$, up to a different elastic constant, hence the use of
Lemma ~\ref{lemma:mon1} here is justified). For the
equality in relation (\ref{rel62}) we use the change of
variables $y=x_1+\frac{x} {\sqrt{e_1^{L_k}}}$ and use the
relation: $e_{\bar
L_k}(R^{(L_k)})=\frac{1}{e_1^{(L_k)}}e_{L_k}(Q^{(L_k)})$. For the
second inequality in (\ref{rel62}) we use the monotonicity
inequality and the fact that
$\frac{2/3r-r_1}{2}\le \frac{r}{3}$. For the third inequality in
(\ref{rel62}) we use the fact that $B_{r/3}(x_1)\subset B_r(a)$
since $|x_1-a|<r_1<\frac{2}{3}r$. The last step in (\ref{rel62}) follows from 
the hypothesis of the Lemma.

\smallskip
\par  Choosing $C_1$ small enough we reach a contradiction with
(\ref{lowerbd}) which in turn implies that $r_2\le 1$   and hence the
conclusion.$\Box$

\par We can now prove the uniform convergence of $Q^{(L_k)}$ away from
singularities of the limiting harmonic map $Q^{(0)}$:

\begin{proposition}
Let $\Omega\subset\mathbb{R}^3$ be a simply-connected bounded open set
with smooth boundary.
 Let $Q^{(L_k)}\in W^{1,2}(\Omega,S_0)$ be a
sequence of global minimizers for the energy $\tilde F_{LG}[Q]$ in the admissible space $\Acal_Q$. Assume that as $L_k\to 0$
we have $Q^{(L_k)}\to Q^{(0)}$ in $W^{1,2}(\Omega,S_0)$.
\par Let $K\subset\Omega$ be a compact set which contains no singularity of
$Q^{(0)}$. Then
\begin{equation}
\lim_{k\to \infty}Q^{(L_k)}(x)=Q^{(0)}(x), \textrm{ uniformly for
} x\in K \label{uniformconv}
\end{equation}
\label{prop:unifconvinterior}
\end{proposition}

\smallskip\par{\bf Proof.} From the hypothesis and {\it Proposition} ~\ref{uniformbulk} we have that
 $\tilde f_B(Q^{(L_k)}(x))\to 0$ uniformly in $K$.
 Thus for any
$\varepsilon_0>0$ there exists a  $\bar L_0>0$ such that for
$L_k<\bar L_0$ we have that $|Q^{(L_k)}(x)-s_+\left(n(x)\otimes
n(x)-\frac{1}{3}Id\right)|<\varepsilon_0$ for all $x\in K$ (and
for each $x\in K$, we have $n(x)\in\mathbb{S}^2$). Thus we can
apply Lemmas ~\ref{lemma:comparison}, ~\ref{lemma:bochner} and
~\ref{lemma:uniformreg}.
\par In order to show the uniform convergence it suffices to show
that we have uniform (independent of $L_k$) Lipschitz bounds on $Q^{(L_k)}(x)$ for
$x\in K$. We reason similarily to the proof in Proposition ~\ref{uniformbulk} (see also \cite{remarks}).  We first claim that there exists an
$\varepsilon_1>0$ so that

\begin{eqnarray}
\forall \varepsilon\in(0,\varepsilon_1),\textrm{ there exists }
r_0(\varepsilon)\textrm{ depending only on
}\varepsilon,\,\Omega, K,\textrm{ and boundary data }Q_b\textrm{ so
that}\nonumber\\
\frac{1}{r_0}\int_{K\cap B_{r_0}(x)}\frac{1}{2}|\nabla
Q^{(L_k)}(x)|^2+\frac{\tilde f_B(Q^{(L_k)}(x))}{L_k}\,dx\le
\varepsilon\,,\forall x\in K,\textrm{ provided that
}L_k<L_*(\varepsilon,r_0(\varepsilon)) \label{smallclaim}
\end{eqnarray}
\par In order to prove the claim let us first recall that  $Q^{(0)}$
has no singularities on the compact set $K$. Thus there exists a larger
compact set $K'$ with $K\subset K'$  and a constant $C>0$ so that
$|\nabla Q^{(0)}(x)|\le C,\forall x\in K'$. We choose
$\varepsilon_1>0$ so that $B(x,\varepsilon_1)\cap K\subset K'$
hence for an arbitrary $\varepsilon\in (0,\varepsilon_1)$ there
exists $r_0(\varepsilon)>0$ so that
$$
\frac{1}{r_0}\int_{K\cap B_{r_0}(x)}\frac{1}{2}|\nabla
Q^{(0)}(x)|^2\,dx<\frac{\varepsilon}{3}
$$ provided that $x\in K$  and $r_0(\varepsilon)$ is chosen small enough. We also have, from the
$W^{1,2}(\Omega,S_0)$
convergence of $Q^{(L_k)}$ to $Q^{(0)}$, that there exists $\bar
L(\varepsilon)$ so that
$$
\frac{1}{r_0}\int_{K\cap B_{r_0}(x)}\frac{1}{2} |\nabla Q^{(L_k)}(x)|^2(x)\,dx\le
\frac{1}{r_0}\int_{K\cap B_{r_0}(x)}\frac{1}{2}|\nabla
Q^{(0)}(x)|^2\,dx+\frac{\varepsilon}{3},\,\forall L_k<\bar
L(\varepsilon) \label{smalllimit+}
$$
\par Recall from the proof of {\it Lemma} ~\ref{lemma:W^{1,2}} that 
 $\lim_{L_k\to 0}\int_\Omega \frac{\tilde
f_B(Q^{(L_k)}(x))}{L_k}\,dx=0$. Hence there exists $\tilde
L(\varepsilon)$ so that $\frac{1}{r_0} \int_\Omega \frac{\tilde
f_B(Q^{(L_k)}(x))}{L_k}\,dx<\frac{\varepsilon}{3}, \forall
L<\tilde L(\varepsilon)$. Letting
$L_*(\varepsilon,r_0(\varepsilon))=\min\{\bar L,\tilde L\}$ and
combining the two relations above we obtain the claim
(\ref{smallclaim}).
\par Choosing $\varepsilon>0$ smaller than the constant $C_1$
from {\it Lemma}~\ref{lemma:uniformreg}, we  apply {\it
Lemma}~\ref{lemma:uniformreg} to conclude that $|\nabla Q^{(L_k)}(x)|$ can be bounded, independently of $L_k$, on the set $K$. The uniform convergence result now follows.
$\Box$

\subsection{The analysis near the boundary}
\label{sec:boundary} In this section we consider the behaviour of
a global minimizer $Q^{(L)}$ near the boundary, $\partial\Omega$, in the
limit $L\to 0$.
For $x^0 \in \partial\Omega$ we define the region $\Omega_r$ to
be:
\begin{equation}
\label{eq:omegar} \Omega_r \stackrel{def}{=} \Omega \cap
B_r(x^0),\, r > 0.
\end{equation}


\begin{lemma} Let $\Omega$ be a simply-connected, bounded open set with
Lipschitz
boundary. There exists a constant $D>0$, depending only on
$\Omega$, and a constant $r_0>0$ such that for all $r<r_0$ and for
any $x^0\in\partial\Omega$, we have:
\begin{equation}
\mathcal{H}^2(\partial \Omega\cap B_r(x^0))\le Dr^2.
\end{equation}
\label{lemma:boundarysize}
\end{lemma}

\smallskip\par{\bf Proof.} Since $\Omega$ has Lipschitz boundary, we have
that
 for any $x^0\in\partial\Omega$, there exists a $\lambda(x^0)>0$ and
an orthonormal coordinate system $X=(x_1,x_2,x_3)$ such
that $x^0=(0,0,0)$ and there exists a Lipschitz function,
$f_{x^0}:\mathbb{R}^2\to\mathbb{R}$, with the property

\begin{displaymath}
\mathcal{U}_{x^0}\stackrel{def}{=}\{x\in\Omega,|x_i|<\lambda(x^0),
i=1,2,3\}=\{x\in\mathbb{R}^3,x_3<f_{x^0}(x_1,x_2),|x_i|<\lambda(x^0),
i=1,2,3\}.
\end{displaymath}
\par As $\Omega$ is bounded, it is necessarily uniformly Lipschitz (see
for instance
\cite{fraenkel}). Hence, for each $x^0\in\partial\Omega$, we can choose the system of coordinates as before
such that there exists a constant
$\bar c>0$, independent of $x^0$, so that $\|\grad f_{x^0}\|\le\bar c,\forall
x^0\in\partial\Omega$.
\par Letting $r_0\stackrel{def}{=}\lambda$ we have:

\begin{equation}
\mathcal{H}^2(\partial\Omega\cap B_r(x^0))\le \int_{[-r,r]^2}
\sqrt{1+|\nabla f_{x^0}(x_1,x_2)|^2}\,dx_1~dx_2\le
\int_{[-r,r]^2}\sqrt{1+\bar c^2}~dx_1~dx_2=4\sqrt{1+\bar c^2}r^2
\nonumber
\end{equation}  $\forall r<r_0$. $\Box$

\par We have a boundary analogue of the interior mononicity lemma, {\it
Lemma} ~\ref{lemma:mon1}, namely :

\begin{lemma}(boundary monotonicity)
\label{boundarymonotonicity} Let $\Omega$ be a simply-connected bounded open set with smooth boundary. Let $Q^{(L)}$ be a global minimizer of
$\tilde F_{LG}[Q]$ in the admissible class $\Acal_Q$.
 Let
\begin{equation}
\Ecal_r = \frac{1}{r}\int_{\Omega_r}\frac{|\nabla Q^{(L)}|^2}{2} +
\frac{\tilde{f}_B(Q^{(L)})}{L}~dV \label{eq:normen}
\end{equation}
\par Then there exists $r_0>0$ so that 
\begin{equation}
\label{eq:bound2} \frac{d}{dr}\Ecal_r \geq - C\left(a^2, b^2,
c^2,Q_b,r_0,\Omega\right),\forall 0<r<r_0
\end{equation}
where the positive constant $C$ is independent of $L$.
\end{lemma}

\par{\bf Proof.} {\it Step 1} We assume that the  domain $\Omega$ is
star-shaped.
Then the proof of (\ref{eq:bound2}) closely follows the arguments in
\cite{linriviere} combined with an idea from \cite{bbh}.
\par Recall that $Q^{(L)}$ satisfies the equation:
\begin{equation}
\Delta Q^{(L)}_{ij}=\frac{1}{L}\left[\frac{\partial\tilde
f_B(Q^{(L)})}{\partial
Q^{(L)}_{ij}}+b^2\frac{\delta_{ij}}{3}\textrm{tr}(Q^{(L)})^2\right]
\label{Eqgrad+}
\end{equation} In what follows, we drop the superscript $L$ for convenience.
\par We multiply both sides of (\ref{Eqgrad+}) by $(x_p -
x^0_p)Q_{ij,p}$ and integrate over $\Omega_r$. Then
\begin{equation}
\label{eq:bound5} \int_{\Omega_r}Q_{ij,kk}(x_p - x^0_p)Q_{ij,p}~dx
= \int_{\partial\Omega_r}Q_{ij,k}Q_{ij,p} (x_p - x^0_p)n_k~d\sigma
- \int_{\Omega_r}|\nabla Q|^2 + Q_{ij,k}Q_{ij,kp} (x_p - x^0_p)~dx
\end{equation}
where $n$ is the unit outward normal to $\partial\Omega_r$ and
$d\sigma$ is the area element on $\partial\Omega_r$.

The integral $\int_{\partial\Omega_r}Q_{ij,k}Q_{ij,p} (x_p -
x^0_p) n_k~d\sigma$ is evaluated by considering the contributions
from $\partial\Omega \cap B_r(x^0) $ and
$\Omega \cap \partial B_r(x^0)$ separately. On
$\Omega \cap \partial B_r$, $n(x)=\frac{x-x^0}{|x-x^0|}$ so that
$$ \int_{\Omega \cap \partial B_r}Q_{ij,k}Q_{ij,p} (x_p - x^0_p) n_k\,d\sigma
 = \int_{\Omega \cap \partial B_r}r~\left|\frac{\partial Q}{\partial
n}\right|^2~d\sigma.$$ Similarly
$$ \int_{\partial\Omega \cap B_r}
Q_{ij,k}Q_{ij,p} (x_p - x^0_p) n_k~d\sigma = \int_{\partial\Omega
\cap B_r} (x - x^0)\cdot n \left|\frac{\partial Q}{\partial
n}\right|^2 + (x - x^0)\cdot \tau \frac{\partial Q_b}{\partial
\tau}\cdot \frac{\partial Q}{\partial n}~d\sigma$$ where
$\tau(x)\in\mathbb{S}^2$ is the tangential direction to the
boundary at $x\in\partial\Omega$.

In order to estimate $\int_{\Omega_r} Q_{ij,k}Q_{ij,kp}
(x_p - x^0_p)~dx$ we note that $$ Q_{ij,k}Q_{ij,kp} (x_p -
x^0_p) = \frac{\partial }{\partial x_p} \left[ (x_p - x^0_p)
\frac{1}{2}|\nabla Q|^2\right] - \frac{3}{2}|\nabla Q|^2$$ and
therefore
$$ \int_{\Omega_r}|\nabla Q|^2 + Q_{ij,k}Q_{ij,kp}
(x_p - x^0_p)~dx =  \int_{\partial\Omega_r}(x_p - x^0_p)
\frac{1}{2}|\nabla Q|^2n_p~d\sigma - \int_{\Omega_r}\frac{1}{2}\left|\grad
Q\right|^2~dx.$$ The surface integral over $\partial\Omega_r$ can
again be expressed in terms of separate contributions from
$\partial\Omega\cap B_r(x^0)$ and $\Omega\cap\partial B_r(x^0)$.

Combining the above, we have
\begin{eqnarray}
\label{eq:step1} && \int_{\Omega_r} Q_{ij,kk}(x_p -
x^0_p)Q_{ij,p}~dx = \int_{\Omega_r}\frac{|\nabla Q|^2}{2}~dx +
r\left( \int_{\Omega\cap\partial B_r}\left|\frac{\partial
Q}{\partial n}\right|^2 - \frac{|\nabla Q|^2}{2}~d\sigma \right)+
\nonumber \\ && + \int_{\partial\Omega\cap B_r}(x - x^0) \cdot n
\left[\frac{1}{2}\left|\frac{\partial Q}{\partial n}\right|^2 -
\frac{1}{2}\left|\frac{\partial
Q_b}{\partial\tau}\right|^2\right]~d\sigma +
\int_{\partial\Omega\cap B_r}(x - x^0) \cdot \tau \frac{\partial
Q_b}{\partial \tau}\cdot \frac{\partial Q}{\partial n}~d\sigma.
\end{eqnarray} In (\ref{eq:step1}), we use the fact that $\left|\grad
Q\right|^2 =
\left|\frac{\partial Q}{\partial n}\right|^2 + \left|\frac{\partial
Q_b}{\partial\tau}\right|^2$ on $\partial \Omega$.

Using the same sort of arguments as above, we compute
\begin{equation}
\label{eq:step2} \frac{1}{L}\int_{\Omega_r} \frac{\partial
\tilde{f}_B}{\partial Q_{ij}}(x_p - x^0_p)Q_{ij,p}~dx
 = \frac{1}{L} \int_{\Omega_r} \frac{\partial}{\partial
x_p}\left[\tilde{f}_B(Q)(x_p
- x^0_p)\right] - 3 \tilde{f}_B(Q)~dx
\end{equation} where $\tilde{f}_B(Q)=\tilde f_B(Q_b)=0$ on
$\partial\Omega$ (from
our choice of the boundary condition $Q_b$ in (\ref{eq:Qb})).

Equating (\ref{eq:step1}) and (\ref{eq:step2})  we obtain
\begin{eqnarray}
&& \label{eq:step3} \int_{\Omega_r}\frac{|\nabla Q|^2}{2}~dx +
r\left( \int_{\Omega\cap\partial B_r}\left|\frac{\partial
Q}{\partial n}\right|^2 - \frac{|\nabla Q|^2}{2}~d\sigma
\right)+\nonumber \\ && + \int_{\partial\Omega\cap B_r}(x - x^0)
\cdot n \left[\frac{1}{2}\left|\frac{\partial Q}{\partial
n}\right|^2 - \frac{1}{2}\left|\frac{\partial
Q_b}{\partial\tau}\right|^2\right]~d\sigma +
\int_{\partial\Omega\cap B_r}(x - x^0) \cdot \tau \frac{\partial
Q_b}{\partial \tau}\cdot \frac{\partial Q}{\partial n}~d\sigma  =
\nonumber \\ && = \int_{\Omega\cap\partial B_r} r \frac{
\tilde{f}_B(Q)}{L}~d\sigma - 3\int_{\Omega_r} \frac{
\tilde{f}_B(Q)}{L}~dx.
\end{eqnarray}

We multiply both sides of (\ref{eq:step3}) by $\frac{1}{r^2}$ and
after some re-arrangement, obtain
\begin{eqnarray}
\label{eq:step4} && -\frac{1}{r^2}\int_{\Omega_r} \frac{|\nabla
Q|^2}{2} + \frac{ \tilde{f}_B(Q)}{L}~dx + \frac{1}{r}\int_{\Omega
\cap \partial B_r}\frac{|\nabla Q|^2}{2} +
\frac{\tilde{f}_B(Q)}{L}~d\sigma = \nonumber \\ && =
\frac{1}{r}\int_{\Omega \cap \partial B_r} \left|\frac{\partial
Q}{\partial n}\right|^2~d\sigma + \frac{2}{r^2}\int_{\Omega_r}
\frac{\tilde{f}_B(Q)}{L}~dx + \nonumber
\\ && + \frac{1}{2r^2}\int_{\partial\Omega\cap B_r}(x - x^0)
\cdot n \left[\left|\frac{\partial Q}{\partial n}\right|^2 -
\left|\frac{\partial Q_b}{\partial\tau}\right|^2\right]~d\sigma +
\frac{1}{r^2}\int_{\partial\Omega\cap B_r}(x - x^0) \cdot \tau
\frac{\partial Q_b}{\partial \tau}\cdot \frac{\partial Q}{\partial
n}~d\sigma.
\end{eqnarray}
For a star-shaped domain $(x - x^0) \cdot n \geq 0$ on
$\partial\Omega$. Therefore, the negative contributions to the right hand side of
(\ref{eq:step4}) are $- \frac{1}{2r^2}\int_{\partial\Omega\cap
B_r}(x - x^0) \cdot n \left|\frac{\partial
Q_b}{\partial\tau}\right|^2~d\sigma$ and potentially
$\frac{1}{r^2}\int_{\partial\Omega\cap B_r}(x - x^0) \cdot \tau
\frac{\partial Q_b}{\partial \tau}\cdot \frac{\partial Q}{\partial
n}~d\sigma$. The first integral can be easily estimated since
$Q_b$ is known. Using the fact that $|\frac{\partial
Q_b}{\partial\tau}|^2 \leq C s_+^2$  for some $C>0$ (as $Q_b\in
C^\infty(\partial\Omega)$ by hypothesis) where $s_+$  is defined in
(\ref{s+}), we have that
\begin{equation}
\label{eq:step5}
 \frac{1}{2r^2}\int_{\partial\Omega\cap
B_r}|(x - x^0) \cdot n| \left|\frac{\partial
Q_b}{\partial\tau}\right|^2~d\sigma\leq C r s_+^2.
\end{equation}
Here we have used $\left|(x - x^0) \cdot n\right|\leq r$ and {\it
Lemma} ~\ref{lemma:boundarysize}.

Using Cauchy-Schwarz, we have
\begin{equation}
\label{eq:step6} \frac{1}{r^2}\int_{\partial\Omega\cap B_r}(x -
x^0) \cdot \tau \frac{\partial Q_b}{\partial \tau}\cdot
\frac{\partial Q}{\partial n}~d\sigma \leq \frac{1}{r}\left(
\int_{\partial\Omega\cap B_r}\left|\frac{\partial Q_b}{\partial
\tau}\right|^2~d\sigma\right)^{1/2}\left(\int_{\partial\Omega\cap
B_r}\left|\frac{\partial Q}{\partial
n}\right|^2~d\sigma\right)^{1/2}.
\end{equation}
The first integral on the right hand side is easily dealt
with i.e. $\int_{\partial\Omega\cap B_r}\left|\frac{\partial
Q_b}{\partial \tau}\right|^2~d\sigma \leq C s_+^2 r^2,$ from
Lemma~\ref{lemma:boundarysize}.

 The second integral involving $\left|\frac{\partial Q}{\partial
n}\right|^2$ is estimated using Lemma~\ref{lem:bbh}:

$$\label{eq:step7} \int_{\partial\Omega\cap B_r}\left|\frac{\partial
Q}{\partial n}\right|^2~d\sigma \leq G\left(Q_b,\Omega\right)$$
where $G>0$ is a constant independent of $L$.

\par Combining the above we have that
\begin{equation}
\label{eq:step8} -\frac{1}{r^2}\int_{\Omega_r} \frac{|\nabla
Q|^2}{2} + \frac{ \tilde{f}_B(Q)}{L}~dx + \frac{1}{r}\int_{\Omega
\cap
\partial B_r}\frac{|\nabla Q|^2}{2} +
\frac{\tilde{f}_B(Q)}{L}~d\sigma \geq - C r s_+^2 -
G^{'}\left(a^2,b^2,c^2,\Omega\right)
\end{equation}
where $C$ and $G^{'}$ are positive constants independent of $L$.
We note that
\begin{equation}
\label{eq:step9} \frac{d}{dr} \Ecal_r =
-\frac{1}{r^2}\int_{\Omega_r} \frac{|\nabla Q|^2}{2} + \frac{
\tilde{f}_B(Q)}{L}~dx + \frac{1}{r}\int_{\Omega \cap \partial
B_r}\frac{|\nabla Q|^2}{2} + \frac{\tilde{f}_B(Q)}{L}~d\sigma
\end{equation} and the above holds for any $0<r<r_0$ where $r_0$ is the constant from Lemma ~\ref{lemma:boundarysize}.
Therefore $$\frac{d}{dr} \Ecal_r \geq -
G^{''}\left(a^2, b^2,
c^2,Q_b,r_0, \Omega \right)
$$ where $G^{''} > 0$ is independent of $L$.

\par {\it Step 2:} General domain $\Omega$.

We do not assume that the domain $\Omega$ is star-shaped and take
into account the perturbation terms induced by omitting this
assumption. As in \cite{linriviere}, the boundary regularity of
the domain implies that
\begin{equation}
\label{eq:nonstarshaped} (x - x^0) \cdot n \geq \left|(x - x^0)
\cdot n\right| - c r^2
\end{equation}
where $c>0$ is independent of $r$ or $x^0\in\partial\Omega$. Then
\begin{equation}
\label{eq:ns2} \frac{1}{2r^2}\int_{\partial\Omega\cap B_r}(x -
x^0) \cdot n \left|\frac{\partial Q}{\partial n}\right|^2~d\sigma
\geq \frac{1}{2r^2}\int_{\partial\Omega\cap B_r}\left|(x - x^0)
\cdot n\right|\left|\frac{\partial Q}{\partial n}\right|^2~d\sigma
- \frac{c}{2}\int_{\partial\Omega\cap B_r}\left|\frac{\partial
Q}{\partial n}\right|^2~d\sigma.
\end{equation}
The inequality (\ref{eq:bound2}) now follows from
Lemma~\ref{lem:bbh}.$\Box$

\begin{lemma}
\label{lem:bbh} Let $Q^{(L)}$ be a minimizer of $\tilde F_{LG}[Q]$ in $\Acal_Q$ (see (\ref{eq:am41})) for a
fixed $L>0$. Then
\begin{equation}
\label{eq:bbh} \int_{\partial \Omega} \left|\frac{\partial
Q^{(L)}}{\partial n}\right|^2~d\sigma \leq G(Q_b,\Omega)
\end{equation}
where $G>0$ only depends on the boundary condition $Q_b$ and $\Omega$.
\end{lemma}

\par {\bf Proof. } The proof  follows closely the arguments of 
Proposition~3 in \cite{bbh}.
Let $V:\Omega \to \Rr^3$ be a smooth vector field on $\Omega$ such
that $V=n$ on $\partial \Omega$. We drop the superscript $L$ for convenience. We multiply (\ref{ELeq}) by $V_p
Q_{ij,p}$ and note that
\begin{equation}
\label{eq:bbh1} \int_{\Omega}Q_{ij,kk}Q_{ij,p}V_p~dx =
\int_{\partial\Omega}\left|\frac{\partial Q}{\partial
n}\right|^2~d\sigma -
\int_{\Omega}Q_{ij,k}\frac{\partial}{\partial x_k}
\left(Q_{ij,p}V_p \right)~dx.
\end{equation}

Proceeding similarly as in \cite{bbh}, we have that
\begin{eqnarray}
\label{eq:bbh2} \int_{\Omega}Q_{ij,k}\frac{\partial}{\partial x_k}
\left(Q_{ij,p}V_p \right)~dx = \int_\Omega Q_{ij,k}Q_{ij,kp}V_p+\int_\Omega Q_{ij,k}Q_{ij,k}\frac{\partial V_p}{\partial x_p}\,dx\nonumber\\
=\int_{\partial\Omega} \frac{|\nabla
Q|^2}{2}~d\sigma + O\left(s_+^2\right), \,\textrm{ as }L\to 0
\end{eqnarray}
Thus,
\begin{equation}
\label{eq:bbh3} \int_{\Omega}Q_{ij,kk}Q_{ij,p}V_p~dx =
\int_{\partial\Omega}\left|\frac{\partial Q}{\partial n}\right|^2
- \int_{\partial\Omega} \frac{|\nabla Q|^2}{2}~d\sigma +
O\left(s_+^2\right),\,\textrm{ as }L\to 0
\end{equation}

On the other hand,
\begin{equation}
\label{eq:bbh4} \frac{1}{L}\int_{\Omega}Q_{ij,p}V_p~\frac{\partial
\tilde{f}_B(Q)}{\partial Q_{ij}}~dx =- \frac{1}{L}\int_{\Omega}
\tilde{f}_B(Q) \nabla \cdot V~dx \leq O(s_+^2),\,\textrm{ as }L\to
0
\end{equation}
since $\frac{1}{L}\int_{\Omega}\tilde{f}_B(Q)~dV \leq C(\Omega)
s_+^2$ from energy minimality and $\tilde f_B(Q_b)=0$ by our
choice of $Q_b$.

Equating (\ref{eq:bbh3}) and (\ref{eq:bbh4}), we obtain
\begin{equation}
\label{eq:bbh5}
 \int_{\partial\Omega}\left|\frac{\partial Q}{\partial n}\right|^2 -
\int_{\partial\Omega} \frac{|\nabla Q|^2}{2}~d\sigma
 = \frac{1}{2}\int_{\partial\Omega}\left|\frac{\partial Q}{\partial
n}\right|^2 -
\left|\frac{\partial Q_b}{\partial \tau}\right|^2~d\sigma \leq C(\Omega)
s_+^2
\end{equation}
and (\ref{eq:bbh}) now follows. $\Box$

We now prove the uniform convergence   of the bulk energy density, $
\tilde{f}_B(Q^{(L)})$, to its minimal value, on compact subsets, $K\subset
\overline{\Omega}$, that do not contain defects of the limiting
harmonic map $Q^{(0)}$. This extends the  result in Proposition~\ref{uniformbulk} where the uniform convergence is proven only  for $K\subset\Omega$.

\begin{proposition}
\label{prop:uniformboundary} Let $Q^{(L)}$ denote a global minimizer of
$\tilde F_{LG}[Q]$ in the admissible space
$\mathcal{A}_Q$ defined in (\ref{eq:am41}). Consider a sequence  $\{Q^{(L_k)}\}_{k\in\mathbb{N}}$  which converges to 
a limiting harmonic map $Q^{(0)}$ strongly in $W^{1,2}(\Omega,S_0)$ as $L_k\to 0$. 
\par  Let $x^0\in \partial\Omega$ be a boundary point. We assume that
the region $\Omega_r$ in (\ref{eq:omegar}) contains no singularity
of the limiting harmonic map $Q^{(0)}$. Then
\begin{equation}
\label{eq:bulk2} \lim_{L_k\to 0} \tilde{f}_B(Q^{(L_k)}(x)) = 0 \quad\forall x \in
\Omega_r
\end{equation}
and the limit is uniform on $\Omega_r$.
\end{proposition}

\textbf{Proof.} We set $\alpha = \tilde{f}_B\left(Q^{(L_k)}(x^0)\right) \geq 0$.
Consider
the region $\Omega_\rho \subset \Omega_r$ where $\rho < r \leq r_0$ (here $r_0$ is the constant from Lemmas ~\ref{lemma:boundarysize} and \ref{boundarymonotonicity}). Then
the boundary monotonicity inequality (\ref{eq:bound2}) implies that
\begin{equation}
\label{eq:bulk4} \Ecal_\rho \leq \Ecal_r +
C(a^2,b^2,c^2,Q_b,,\Omega)\left(r - \rho\right)
\end{equation} for $\rho<r<r_0$.

\par Take an arbitrary $\varepsilon>0$. Recall that  $Q^{(L_k)} \rightarrow Q^{(0)}$ in $W^{1,2}$ as $L_k\to 0$ and
$\Omega_r$ contains no singularities of $Q^{(0)}$. Using the same
arguments as in Proposition~\ref{uniformbulk}, we have that there exists an $r_1<\min\{r_0,\varepsilon\}$ and $\bar L>0$ ( both depending on $\varepsilon$) so that
for  $L_k<\bar L$
$$ \frac{1}{r_1}\int_{\Omega_{r_1}}\frac{1}{2}|\nabla Q^{(L_k)}|^2~dx \leq \varepsilon
$$

Similarly, we have that there exists an $\tilde L>0$ (depending on $\varepsilon$) so that
$$ \frac{1}{r_1}\int_{\Omega_{r_1}}\frac{\tilde{f}_B\left(Q^{(L_k)}\right)}{L_k}~dV\leq
\varepsilon $$ for $L_k<\tilde L$ ( see the proof of Lemma ~\ref{lemma:W^{1,2}}).
Combining the above, we obtain
\begin{equation}\label{eq:bulk5}
\frac{1}{\rho}\int_{\Omega_\rho}
\frac{\tilde{f}_B\left(Q^{(L_k)}\right)}{L_k}~dx \leq \Ecal_\rho \leq
C'\varepsilon\end{equation} for any $\rho<r_1$ and for $L_k<\min\{\bar L,\tilde L\}$ where the constant $C^{'} >0$ is independent of $L_k$  .

\par Using arguments very close to those in \cite{bbh} ( Lemma $A.2$  and the way it is used in Step $B.1$ of the proof of Theorem $1$) together with Proposition ~\ref{prop:max}, one can  easily obtain:
\begin{equation}
\label{eq:bulk3} \|\nabla Q^{(L_k)}\|_{L^\infty(\Omega)} \leq
\frac{H(a^2,b^2,c^2,\Omega)}{\sqrt{L_k}}.
\end{equation}
On the other hand, $\tilde{f}_B(Q)$ is a Lipschitz function of the
$Q$-tensor and
one can infer the following from (\ref{eq:bulk3}) and Proposition ~\ref{prop:max}:
\begin{equation}
\label{eq:bulk6} \|\nabla \tilde{f}_B\left(Q^{(L_k)}\right)\|_{L^\infty} \leq \frac{
D(a^2,b^2,c^2,\Omega)}{\sqrt{L_k}}
\end{equation}
so that
\begin{equation}
\label{eq:bulk7} \tilde{f}_B\left(Q^{(L_k)}(x)\right) \geq \alpha -
\frac{ D(a^2,b^2,c^2,\Omega)}{\sqrt{L_k}}\rho \quad\forall x \in
\Omega_\rho.
\end{equation}

We take
$$\rho = \frac{\alpha \sqrt{L_k}}{2 D(a^2,b^2,c^2,\Omega)} . $$      

There exists a constant $\gamma(\Omega)$  so that
$$\left| \Omega_\rho\right|\geq \ga(\Omega) \rho^3 $$ (see  also \cite{bbh} for the 2D version of the above) Combining the above with (\ref{eq:bulk7}) , we obtain the following inequality
\begin{equation}
\label{eq:bulk8} \frac{1}{\rho}\int_{\Omega_\rho}
\frac{\tilde{f}_B\left(Q^{(L_k)}\right)}{L_k}~dx \geq \frac{\ga
\rho^2}{L_k}\left(\alpha - \frac{
D(a^2,b^2,c^2,\Omega)}{\sqrt{L_k}}\rho\right) = \ga
\frac{\alpha^3}{D^{'}(a^2,b^2,c^2,\Omega)}
\end{equation}
where the constant $D^{'}>0$ is independent of $L_k$. Combining
(\ref{eq:bulk5}) and (\ref{eq:bulk8}), we have that
\begin{equation}
\label{eq:bulk9} \alpha^3 \leq
D^{''}(a^2,b^2,c^2,Q_b,\Omega)\varepsilon
\end{equation} where $D^{''}>0$ is independent of $L_k$. The upper bound
(\ref{eq:bulk9}) is independent of $x^0$ and  $\varepsilon >0$ was chosen arbitrarily. Therefore,
Proposition~\ref{prop:uniformboundary} now follows. $\Box$

\section{Consequences of the convergence results}
\label{sec:consequence}

In this section, we discuss some consequences of the convergence
results in Propositions ~\ref{uniformbulk},
\ref{prop:unifconvinterior}, and \ref{prop:uniformboundary}. We consider a sequence of  global
minimizers $\{Q^{(L_k)}\}_{k\in\mathbb{N}}$  converging to a limiting harmonic map $Q^{(0)}$. From
Proposition~\ref{prop:unifconvinterior}, we have that for  a ball $B(x,r_0)\subset \Omega$, where
$B(x,r_0)$ does not contain any singularities of  $Q^{(0)}$
\begin{equation}
\label{eq:cons3} \left| Q^{(L_k)}(y) - Q^{(0)} (y) \right| \leq \eps (L_k)
\quad y \in B(x,r_0)
\end{equation}
where $\eps (L_k) \to 0^+$ as $L_k\to 0$. Further, the small energy
regularity in Lemma~\ref{lemma:uniformreg} implies that for $L_k$
sufficiently small,
\begin{equation}
\label{eq:cons1} e_{L_k}\left(Q^{(L_k)}(y)\right) =
\frac{\tilde{f}_B(Q^{(L_k)}(y))}{L_k} + \frac{1}{2}\left|\grad Q^{(L_k)}(y)
\right|^2 \leq C\left(a^2, b^2, c^2, \Omega\right) \quad y \in
B(x,r_0)
\end{equation}
where $C>0$ is a positive constant independent of $L_k$. Therefore,
for sufficiently small $L_k$, one has
\begin{eqnarray}
&& |\grad Q^{(L_k)}(y)|^2 \leq C\left(a^2, b^2, c^2, \Omega\right)
\nonumber \\ && \tilde{f}_B(Q^{(L_k)}(y)) \leq C\left(a^2, b^2, c^2,
\Omega\right) L_k \quad y \in B(x,r_0). \label{eq:cons2}
\end{eqnarray}
One immediate consequence of the uniform convergence in
(\ref{eq:cons3}) and the bounds in (\ref{eq:cons2}) is the
following

\begin{lemma}
\label{lem:cons1} Let $Q^{(L)}$ denote a global minimizer of
$\tilde F_{LG}[Q]$ in the admissible class $\Acal_Q$. Consider a sequence  $\{Q^{(L_k)}\}_{k\in\mathbb{N}}$  which converges to 
a limiting harmonic map $Q^{(0)}$ strongly in $W^{1,2}(\Omega,S_0)$ as $L_k\to 0$. Let
$x\in\overline{\Omega}$ be such that
$B(x,r_0)\cap\overline{\Omega}$ (for $r_0$ smaller than the $r_0$ used in Lemma  ~\ref{lemma:boundarysize}), does not contain
any singularities of the limiting map $Q^{(0)}$. Then
\begin{eqnarray}
\label{eq:cons4} && Q^{(L_k)}(y) = S^{(L_k)} \left( n^{(L_k)} \otimes n^{(L_k)} -
\frac{1}{3}Id\right) + R^{(L_k)} \left(m^{(L_k)}\otimes m^{(L_k)} - p^{(L_k)}\otimes
p^{(L_k)}\right) \quad \nonumber \\&& where~|S^{(L_k)} - s_+| \leq \eps_1(L_k),~|R^{(L_k)}|
\leq \eps_2(L_k)
\end{eqnarray}
with $n^{(L_k)}, m^{(L_k)}$ and $p^{(L_k)}$  unit eigenvectors of $Q^{(L_k)}$, and $ \eps_1(L_k),
\eps_2(L_k) \to 0^+$ as $L_k \to 0$. Secondly, if $x\in \Omega$ is an
interior point such that $B(x,r_0)\subset \Omega$ does not contain
any singularities of $Q^{(0)}$, then we also have that
\begin{equation}
\label{eq:cons4new}
(n^{(L_k)} \cdot n^{(0)})^2 \geq 1 - \eps_3(L_k)
\end{equation}
where $\eps_3(L_k)\to 0^+$ as $L_k\to 0$ and $n^{(0)}$ has been defined in
(\ref{eq:n0}).
\end{lemma}

\textbf{Proof} The representation (\ref{eq:cons4}) is a direct
consequence of Propositions ~\ref{prop:secondrep},~\ref{uniformbulk} and
~\ref{prop:uniformboundary}. In the following we drop the superscripts $L_k$ for convenience, but keep the superscript $0$ in $Q^{(0)}$ and $n^{(0)}$. From Proposition~\ref{uniformbulk} and
Proposition~\ref{prop:uniformboundary}, we have that
$$ \tilde{f}_B(Q(y)) \to 0 \quad as ~ L_k \to 0$$ for $y \in
B(x,r_0)\cap\overline{\Omega}$ where
$B(x,r_0)\cap\overline{\Omega}$ does not contain any singularities
of $Q^{(0)}$. The bulk energy density $\tilde{f}_B(Q)$ is a smooth
function of the order parameters $(S,R)$ in
Proposition~\ref{prop:secondrep}. Therefore, as
$\tilde{f}_B(Q(y)) \to 0$, the corresponding order parameters
$(S,R)$ approach the bulk energy minimum defined by $(S, R) =
(s_+,0)$ and the inequalities (\ref{eq:cons4}) follow. Further, if $B(x,r_0) \subset \Omega$, then the uniform
convergence (\ref{eq:cons3}) holds. A direct computation shows
that for
$$Q = S\left( n \otimes n -
\frac{1}{3}Id\right) + R\left(m\otimes m - p\otimes p\right),$$ we
have
\begin{equation}
\left|Q(y) - Q^{(0)}(y)\right|^2 = \frac{2}{3}\left(S^2 +
s_+^2\right) - 2S s_+\left[\left(n\cdot n^{(0)}\right)^2 -
\frac{1}{3}\right] + 2 R^2 - 2s_+ R\left(\left(m\cdot
n^{(0)}\right)^2 - \left(p\cdot n^{(0)}\right)^2 \right).
\nonumber
\end{equation} The lower bound on $\left(n\cdot
n^{(0)}\right)^2$ now follows from (\ref{eq:cons3}) and the fact that 
$|S-s_+|\leq \eps_1(L_k),|R|\leq \eps_2(L_k)$ for sufficiently small values of $L_k$. $\Box$

\begin{proposition}
\label{lem:cons2} Let $Q^{(L)}$ denote a global minimizer of
$\tilde F_{LG}[Q]$ in the admissible space $\Acal_Q$.  Consider a sequence  $\{Q^{(L_k)}\}_{k\in\mathbb{N}}$  which converges to 
a limiting harmonic map $Q^{(0)}$ strongly in $W^{1,2}(\Omega,S_0)$ as $L_k\to 0$. Then $Q^{(L_k)}$
converges uniformly to the limiting harmonic map $Q^{(0)}$, away from
the singular set of $Q^{(0)}$, in the interior of $\Omega$. Let
$K\subset \Omega$ be an interior subset that does not contain any
singularities of $Q^{(0)}$. Then (i)\begin{eqnarray}
\label{eq:cons6} && \beta( Q^{(L_k)}(y)) \leq C L_k \nonumber \\
&& \left| \left|Q^{(L_k)}(y)\right| - \sqrt{\frac{2}{3}}s_+\right| \leq
D \sqrt{L_k} \quad y \in K
\end{eqnarray} where $C$ and $D$ are positive constants
independent of $L_k$. (ii)(rate of convergence of eigenvalues) Let
$\left\{\la^{(L_k)}_i\right\}$ denote the set of eigenvalues of $Q^{(L_k)}$
and $\left\{\la_i\right\}$ denote the set of eigenvalues of $Q^{(0)}$.
Then
\begin{equation}
\label{eq:second7} \left| \la^{(L_k)}_i(y) - \la_i(y)\right|^2 \leq
\alpha(a^2,b^2,c^2) L_k \quad y \in K;~i=1\ldots 3
\end{equation} where $\alpha$ is a positive constant independent of $L_k$.
\end{proposition}
\textbf{Proof} (i) This follows directly from (\ref{eq:cons2}) and
Proposition~\ref{prop:beta}. In Proposition~\ref{prop:beta}, we
obtain a lower bound for $\tilde{f}_B(Q)$ in terms of $|Q|$ and
$\beta$ and in (\ref{eq:cons2}) we have an upper bound for $\tilde f_B(Q)$ in terms of $L_k$ as shown below
$$ \left[\frac{
 \frac{2}{3} c^2 s_+^2 + a^2}{2}\right]\left(|Q| - \sqrt{\frac{2}{3}}s_+\right)^2 +
\frac{b^2}{6\sqrt{6}}\beta(Q)|Q|^3\leq
\tilde{f}_B\left(Q(y)\right) \leq C\left(a^2,b^2,c^2,\Omega\right)
L_k \quad y\in K
$$

(ii)  In the following we drop the superscripts ``$(L_k)$'' for convenience. From (\ref{eq:cons2}), we have the following upper bound for
the bulk energy density on the set $K\subset\Omega$
$$
\tilde{f}_B(Q(y)) \leq C~L_k  \qquad y \in K
$$
where $C$ is a positive constant independent of $L_k$. Using the
representation formula (\ref{eq:cons4}), we have that
\begin{equation}
\label{eq:second9} Q = S\left(n\otimes n - \frac{1}{3}Id\right) +
R\left(m\otimes m - p \otimes p\right)
\end{equation}
where \begin{equation} \left| S - s_+\right| \leq \eps_4(L_k) =o(1)
\label{eq:second11}
\end{equation} and
\begin{equation}
\label{eq:second12} \left| R \right| \leq \eps_5(L_k) =o(1).
\end{equation}

A direct computation shows that
\begin{equation}
\label{eq:second10} Q_{ij}Q_{ij} = \frac{2}{3}S^2 + 2R^2 \quad\quad Q_{ip}Q_{pj}Q_{ij} = \frac{2S^3}{9} - 2SR^2.
\end{equation}
From (\ref{eq:second11}) and (\ref{eq:second12}), we represent
$Q$ on the subset $K \subset \Omega$ as follows:
\begin{equation}
\label{eq:second13} Q = \left(s_+ + \eps \right)\left(n \otimes
n - \frac{1}{3}Id\right) + \ga \left(m\otimes m - p\otimes p\right)
\end{equation}
where $|\eps|,|\ga| =o(1)$. Using (\ref{eq:second10}), we find that
$$ \left|Q\right|^2 = \frac{2}{3}\left(s_+^2 + \eps^2 + 2s_+ \eps
\right) + 2\ga^2 $$
 and from the maximum principle (Proposition~\ref{prop:max}),
$$|Q(x)|^2 \leq \frac{2}{3}s_+^2 \quad x \in K.$$ This necessarily
implies that
$\eps \leq 0$.

The bulk energy density $\tilde{f}_B$ is given by
\begin{equation}
\label{eq:second14} \tilde{f}_B(Q) = \frac{a^2}{3}\left(s_+^2 -
S^2\right) + \frac{2 b^2}{27}\left(s_+^3 - S^3\right) -
\frac{c^2}{9}\left(s_+^4 - S^4\right) - a^2 R^2 +
\frac{2b^2}{3}SR^2 + \frac{2c^2}{3}S^2 R^2 + c^2 R^4
\end{equation}
where we have merely expressed $\textrm{tr}Q^2$ and
$\textrm{tr}Q^3$ in terms of the order parameters $S$ and $R$. We
write the bulk energy density as the sum of two contributions -
\begin{equation}
\label{eq:second15} \tilde{f}_B(Q) = F\left(S \right) + G\left(
S,R \right)
\end{equation}
where
$$ F\left(S \right)  = \frac{a^2}{3}\left(s_+^2 - S^2\right) + \frac{2
b^2}{27}\left(s_+^3 - S^3\right) - \frac{c^2}{9}\left(s_+^4 - S^4\right)
$$ and
$$G\left( S,R \right) = - a^2 R^2 + \frac{2b^2}{3}SR^2 + \frac{2c^2}{3}S^2
R^2 + c^2
R^4.$$

The function $F(S)$ is analyzed in (\ref{eq:am35}); the function
$F(S)$ is bounded from below by
\begin{equation}
\label{eq:second16} F(S) \geq D(a^2,b^2,c^2)~\left( S -
s_+\right)^2,\quad D(a^2,b^2,c^2)\ge 0
\end{equation}
Similarly, since $2c^2 s_+^2 = b^2 s_+ + 3a^2$ and $0< s_+ - S =o(1)$ (for $L_k$ sufficiently small), we have the following inequality
\begin{equation}
\label{eq:second17} - a^2 R^2 + \frac{2b^2}{3}SR^2 +
\frac{2c^2}{3}S^2 R^2 \geq \frac{b^2s_+}{2} \ga^2.
\end{equation}
Combining (\ref{eq:second16}), (\ref{eq:second17}) and
(\ref{eq:cons2}), we obtain the following
\begin{equation}
\label{eq:second18} D(a^2,b^2,c^2)~\eps^2 + \frac{b^2}{2} s_+
\ga^2 + c^2 \ga^4 \leq \tilde{f}_B\left(Q(x)\right) \leq C~L_k
\end{equation}
from which we deduce
$$\eps^2 \leq C_1 ~L_k\quad and \quad \ga^2 \leq C_2 ~L_k,$$
where $C_1, C_2$ are positive constants independent of $L_k$. The
inequalities (\ref{eq:second7}) now follow. $\Box$

Next, we have a lemma about the leading eigenvector $n$ in the
representation (\ref{eq:cons4}).

\begin{lemma}\label{lem:eigenvector} Let
$Q=S\left(n\otimes n-\frac{1}{3}Id\right)+R\left(m\otimes
m-p\otimes p\right)$ with $S>8|R|$ and $n,m, p\in\mathbb{S}^2$, pairwise perpendicular. Then the minimum of $$\left|Q - s_+\left(a\otimes a-\frac{1}{3}Id\right)\right|^2 $$ with
$a\in\mathbb{S}^2$ is attained by $a=\pm n$.
\end{lemma}

\smallskip\par{\bf Proof.}  A direct computation shows that
\begin{eqnarray}
&& \left|Q - s_+\left(a\otimes a-\frac{1}{3}Id\right)\right|^2 = \frac{2}{3}\left(S^2 + s_+^2 + S s_+\right) + 2 R^2 - 2 S s_+\left(n\cdot a \right)^2 - 2 s_+R \left(\left(m\cdot a\right)^2 - \left(p\cdot a \right)^2 \right) \nonumber \\ &&  = \frac{2}{3}\left( S^2 + s_+^2 + Ss_+\right)+ 2 R^2 + 2s_+R - 2s_+\left(S + R\right) \left(n\cdot a \right)^2 - 4 s_+ R \left(m\cdot a\right)^2
\label{eq:leadingeigenvector}
\end{eqnarray} where in the last line of (\ref{eq:leadingeigenvector}), we use the equality $\left(n\cdot a \right)^2 +\left(m\cdot a\right)^2 + \left(p\cdot a \right)^2 = 1$. 
Since $S>8|R|$, one can immediately verify that (\ref{eq:leadingeigenvector}) is minimized for $\left(n\cdot a \right)^2=1$ or equivalently $a=\pm n$.
 $\Box$

\par We can now provide a result about the regularity the leading ``eigendirection'' $n\otimes n\in M^{3\times 3}$ where $n\in \mathbb{S}^2$ is the leading eigenvector. For a thorough discussion about the relationships between the regularity of the eigenvector $n\in\mathbb{S}^2$ and that of the eigendirection $n\otimes n\in M^{3\times 3}$ see \cite{bz}.

 \begin{corollary}
 \label{eigenvector}
 Let $Q^{(L)}$ denote a global minimizer of
$\tilde F_{LG}[Q]$ in the admissible class $\Acal_Q$. Consider a sequence  $\{Q^{(L_k)}\}_{k\in\mathbb{N}}$  which converges to 
a limiting harmonic map $Q^{(0)}$ strongly in $W^{1,2}(\Omega,S_0)$ as $L_k\to 0$. Let $K\subset \Omega$ be a compact
 subset of $\Omega$ that does not contain singularities of the
 limiting map $Q^{(0)}$. Then, for $L_k$ small enough (depending on $K$),  $Q^{(L_k)}$ can be represented as in
 (\ref{eq:cons4}) on the set $K\subset \Omega$ and the
 \emph{leading eigendirection} $n^{(L_k)}\otimes n^{(L_k)} \in C^\infty\left(K;M^{3\times 3}\right)$.
 \end{corollary} \textbf{Proof} From (\ref{eq:cons4}), we can
 represent $Q^{(L_k)}$ as
 $$ Q^{(L_k)}(x) = S^{(L_k)}\left(n^{(L_k)}\otimes n^{(L_k)}-\frac{1}{3}Id\right)+R^{(L_k)}\left(m^{(L_k)}\otimes
m^{(L_k)}-p^{(L_k)}\otimes p^{(L_k)}\right) $$
 where $|S^{(L_k)} - s_+|=o(1), |R^{(L_k)}| =o(1)$, and $n^{(L_k)}, m^{(L_k)}, p^{(L_k)} \in \mathbb{S}^2$ are the eigenvectors of $Q^{(L_k)}$. 

 \par Let $\pi(Q)$ be the nearest neighbor projection
onto the manifold of global minimizers of the bulk energy density,
denoted by $Q_{\min} = \{s_+\left(a\otimes
a-\frac{1}{3}Id\right),a\in\mathbb{S}^2\}$ as in (\ref{eq:Qmin}).
Namely, $\pi(Q)$ associates with each $Q'$, (in a neighborhood of
the manifold $Q_{\min}$) an element $Q^*\in Q_{\min}$ such that
$$\left| Q' - Q^* \right| = \min_{Q\in Q_{\min}}\left| Q' - Q \right|.$$
The projection $\pi$  is defined only in a neighborhood of the
manifold $Q_{\min}$ and moreover $\pi(Q') \in C^\infty\left(S_0, Q_{\min}\right)$ (see, for
instance, \cite{remarks}). The Lemma ~\ref{lem:cons1} and Lemma ~\ref{lem:eigenvector} show that in our
case $$\pi(Q^{(L_k)}) = s_+\left(n^{(L_k)}\otimes
n^{(L_k)}-\frac{1}{3}Id\right).$$ Therefore, the tensor $$(n^{(L_k)}\otimes n^{(L_k)} - \frac{1}{3} Id) \in C^\infty\left(K, S_0\right),$$ (since $s_+$ is a constant) and the conclusion of the lemma now follows. 
$\Box$

\section{ Biaxiality and uniaxiality}
\label{sec:6}

\subsection{The bulk energy density}

Our first proposition concerns the stationary
points of the bulk energy density.
\begin{proposition} \cite{bm}
\label{prop:bulk} Consider the bulk energy density
$\tilde{f}_B(Q)$ given by
\begin{equation}
\tilde{f}_B(Q) = -\frac{a^2}{2}\textrm{tr} Q^2 -
\frac{b^2}{3}\textrm{tr} Q^3 + \frac{c^2}{4}\left(\textrm{tr}
Q^2\right)^2 +\frac{a^2}{3}s_+^2 + \frac{2b^2}{27}s_+^3 -
\frac{c^2}{9}s_+^4. \label{eq:am15}
\end{equation}
Then $\tilde{f}_B(Q)$ attains its minimum for uniaxial $Q$-tensors
of the form
\begin{equation}
Q = s_+\left(n\otimes n - \frac{1}{3}\right), \label{eq:am16}
\end{equation}
where \begin{equation} s_+ = \frac{b^2 + \sqrt{b^4 + 24 a^2
c^2}}{4c^2} \label{eq:am17}
\end{equation} and $n:\Omega \rightarrow S^2$ is a unit eigenvector of
$Q$.
\end{proposition}

\textbf{Proof.} Proposition~\ref{prop:bulk} has been proven in
\cite{bm} and we reproduce the proof in the {\it Appendix} for
completeness. $\Box$

In the following proposition, we estimate $\tilde{f}_B(Q)$ in
terms of $\left|Q\right|$ and the biaxiality parameter $\beta(Q)$.

\begin{proposition}
\label{prop:beta} Let $Q\in S_0$. Then the bulk energy density
$\tilde{f}_B(Q)$ is bounded from below by
\begin{equation}
\label{eq:am18}\tilde{f}_B(Q) \geq \left[\frac{ \frac{2}{3}c^2s_+^2 +
a^2}{2}\right]\left(|Q| - \sqrt{\frac{2}{3}}s_+\right)^2 +
\frac{b^2}{6\sqrt{6}}\beta(Q)|Q|^3
\end{equation}
where $s_+$ has been defined in (\ref{eq:am17}).
\end{proposition}

\textbf{Proof.} From Lemma~\ref{betalemma}, we have the
inequality,
$$ \textrm{tr}Q^3 = |Q|^3\sqrt{\left(\frac{1 -
\beta}{6}\right)} \leq \frac{|Q|^3}{\sqrt{6}}\left(1 -
\frac{\beta}{2}\right)~for ~Q\in S_0. $$

From the definition of $\tilde{f}_B(Q)$ and $s_+$  in
(\ref{eq:am15}) and (\ref{eq:am17}), we can obtain a lower bound for $\tilde{f}_B(Q)$
in terms of $|Q|$ and $\beta(Q)$ as follows i.e.
\begin{eqnarray}
&&\tilde{f}_B(Q) = - \frac{a^2}{2}|Q|^2 -
\frac{b^2}{3\sqrt{6}}|Q|^3 \sqrt{1 - \beta} + \frac{c^2}{4}|Q|^4
+\frac{a^2}{2}\left(\sqrt{\frac{2}{3}}s_+\right)^2 + \frac{b^2}{3}\frac{2s^3_+}{9} - \frac{c^2}{4}\left(\sqrt{\frac{2}{3}}s_+\right)^4
\label{eq:am19} \\ && \qquad \geq  \left[- \frac{a^2}{2}|Q|^2 -
\frac{b^2}{3\sqrt{6}}|Q|^3 + \frac{c^2}{4}|Q|^4 +
\frac{a^2}{3}s_+^2 + \frac{2b^2}{27}s_+^3 -
\frac{c^2}{9}s_+^4 \right] +
\frac{b^2}{6\sqrt{6}}\beta(Q)|Q|^3. \label{eq:am23}
\end{eqnarray}

The bracketed term in (\ref{eq:am23}) can be further simplified by
carrying out a series of calculations. Consider the function
\begin{equation}
f(u) = - \frac{a^2}{2} u^2 - \frac{b^2}{3\sqrt{6}} u^3 +
\frac{c^2}{4} u^4.\label{eq:am20}
\end{equation} The stationary points of $f(u)$ are solutions of
the algebraic equation
\begin{equation}
\label{eq:am21} f^{'}(u) = u\left(c^2 u^2 - \frac{b^2}{\sqrt{6}}u
- a^2\right) = 0
\end{equation} and one can readily verify that $f(u)$ attains its
minimum for
\begin{equation}
u_{\min} = \sqrt{\frac{2}{3}}s_+ \label{eq:am22}
\end{equation}

The bracketed term in (\ref{eq:am23}) is non-negative by virtue of
(\ref{eq:am20})--(\ref{eq:am22}). Further, let $\delta = |Q| -
\sqrt{\frac{2}{3}}s_+$ where $c^2 \frac{2}{3}s_+^2 =
\frac{b^2}{\sqrt{6}}\sqrt{\frac{2}{3}}s_+ + a^2$ by the definition of
$s_+$. Then
\begin{eqnarray}
\label{eq:am24} && \left[- \frac{a^2}{2}|Q|^2 -
\frac{b^2}{3\sqrt{6}}|Q|^3 + \frac{c^2}{4}|Q|^4 +
\frac{a^2}{3}s_+^2 + \frac{2b^2}{27}s_+^3 -
\frac{c^2}{9}s_+^4  \right] =
\\ && = \delta\left[-a^2
\sqrt{\frac{2}{3}}s_+ - \frac{\sqrt{2}b^2}{3\sqrt{3}}s_+^2 + 
\frac{2\sqrt{2}}{3\sqrt{3}}c^2s_+^3\right] + \delta^2\left[-\frac{a^2}{2}
-\frac{b^2}{3}s_+ +s_+^2 c^2
\right] + \nonumber \\ && \qquad + \delta^3\left[c^2 \sqrt{\frac{2}{3}}s_+ -
\frac{b^2}{3\sqrt{6}}\right] + \frac{c^2}{4}\delta^4.
\end{eqnarray} The coefficient of $\delta$ vanishes from the
definition of $s_+$ in(\ref{eq:am17}). The coefficients of $\delta^2$ and
$\delta^3$ are positive since
\begin{eqnarray}
&& -\frac{a^2}{2} -\frac{b^2}{3}s_+ +
s_+^2 c^2 \geq \left[\frac{ \frac{2}{3}c^2s_+^2 +
a^2}{2}\right] \nonumber \\ &&  \sqrt{\frac{2}{3}} c^2s_+ - \frac{b^2}{3
\sqrt{6}} \geq \frac{2b^2}{3\sqrt{6}} ~\label{eq:am25}.
\end{eqnarray}
We substitute (\ref{eq:am25}) into (\ref{eq:am24}) to obtain
\begin{eqnarray}\label{eq:am26}
&& \left[- \frac{a^2}{2}|Q|^2 - \frac{b^2}{3\sqrt{6}}|Q|^3 +
\frac{c^2}{4}|Q|^4 +\frac{a^2}{3}s_+^2 + \frac{2b^2}{27}s_+^3 -
\frac{c^2}{9}s_+^4 
\right] \geq \nonumber \\ && \qquad \geq \delta^2 \left[\frac{\frac{2}{3}c^2
s_+^2 + a^2}{2}\right]
\end{eqnarray}
and on combining (\ref{eq:am26}) with (\ref{eq:am23}), the lower
bound (\ref{eq:am18}) follows. $\Box$

The bulk energy density, $\tilde{f}_B(Q)$, can be equivalently
expressed in terms of the order parameters $s$ and $r$ in
Proposition~\ref{prop:rep1}, as shown below

\begin{proposition}\label{prop:bulk2}
Let $Q\in S_0$ be represented as in Proposition~\ref{prop:rep1}
$$ Q = s\left( n\otimes n - \frac{1}{3}Id\right) + r\left(m\otimes m -
\frac{1}{3}Id\right)$$
with either $0\leq r \leq \frac{s}{2}$ or $\frac{s}{2}\leq r \leq
0$. Case (i) Non-negative order parameters, $0\leq r \leq
\frac{s}{2}$ with $0\leq s \leq s_+$, where $s_+$ is defined in
(\ref{eq:am17}). Then the bulk energy density, $\tilde{f}_B(Q)$,
is bounded from below by
\begin{equation}
\tilde{f}_B(Q) \geq \left(s_+ - s\right)^2\ga(a^2, b^2,c^2) +
\frac{r(s-r)}{9}\left(3a^2 + b^2 s - 2c^2 s^2 \right) +
\frac{5b^2}{27} r^2 s \quad 0\leq s \leq s_+
\label{eq:am27}
\end{equation}
where $\ga\left(a^2, b^2, c^2\right)$ is an explicitly computable
positive constant.

Case (ii) Non-negative order parameters, $0\leq r \leq
\frac{s}{2}$ and $s \geq s_+$. Then
\begin{equation}
\tilde{f}_B(Q) \geq \left[\frac{\frac{2}{3}c^2 s_+^2 +
a^2}{2}\right]\min \left\{\frac{2}{3}\left(s -
s_+\right)^2,~\frac{1}{6}\left(\sqrt{3}s - 2 s_+\right)^2\right\}
+ \tau b^2 s_+^{3}\left(\frac{r^2(s - r)^2}{s^4}\right)
\label{eq:am28}
\end{equation} where $\tau$ is an explicitly computable positive constant,
independent of $a^2,b^2,c^2$.

Case (iii) If $\frac{s}{2}\leq r \leq 0$, then
\begin{equation}
\label{eq:negative} \tilde{f}_B(Q) =  \tilde{f}_B(-Q)+
\frac{2b^2}{27}\left(2|s|^3 + 2|r|^3 - 3s^2 |r| - 3|s|r^2 \right),
\end{equation} where $-Q\in S_0$ has positive order parameters $0 \leq -r
\leq
-\frac{s}{2}$ and $ \tilde{f}_B(-Q)$ can be estimated using
(\ref{eq:am27}) and
(\ref{eq:am28}).
In particular,
\begin{equation}
\label{eq:rsneg} \tilde{f}_B(Q) \geq -\frac{a^4}{4 c^2} -
\frac{s_+^3}{3}\left(\frac{b^2}{9} - \frac{c^2}{3} s_+\right) > 0
\end{equation} for $Q$-tensors with $\frac{s}{2}\leq r \leq 0$.
\end{proposition}

\textbf{Proof.} From Proposition~\ref{prop:rep1}, it suffices to
consider the two cases $0\leq r \leq \frac{s}{2}$ and
$\frac{s}{2}\leq r \leq 0$.

Case (i): We can explicitly express the bulk energy density,
$\tilde{f}_B(Q)$, in terms of $s$ and $r$ as follows -
\begin{eqnarray}
\label{eq:am29} && \tilde{f}_B(Q) = -\frac{a^2}{3}\left(s^2 + r^2
- sr \right) - \frac{b^2}{27}\left(2s^3 + 2r^3 - 3s^2 r - 3s r^2
\right) + \nonumber \\ && \qquad + \frac{c^2}{9}\left(s^4 + r^4 +
3s^2 r^2 - 2sr^3 - 2s^3 r \right)  + \frac{a^2}{3}s_+^2 +
\frac{2b^2}{27} s_+^3 - \frac{c^2}{9}s_+^4,
\end{eqnarray} where we have expressed $\textrm{tr}Q^2$ and $\textrm{tr}Q^3$
in terms of $s$ and $r$ $$ \textrm{tr}Q^2 = \frac{2}{3}\left( s^2
+ r^2 - sr \right)$$ and
$$ \textrm{tr}Q^3 = \frac{1}{9}\left(2s^3 +
2r^3 - 3s^2 r - 3sr^2\right).$$
 The function
$\tilde{f}_B(Q)$ consists of two components -
\begin{eqnarray}
\label{eq:am30} && \tilde{f}_B(Q) = F(s) + G(s,r) \quad
where\nonumber
\\ && F(s) = -\frac{a^2}{3}\left(s^2 - s_+^2\right) -
\frac{2b^2}{27}\left(s^3 - s_+^3\right) + \frac{c^2}{9}\left(s^4 -
s_+^4\right) \nonumber \\ && G(s,r) = \frac{a^2}{3}\left(sr -
r^2\right) + \frac{b^2}{27}\left(3s^2 r + 3sr^2 - 2r^3\right) +
\frac{c^2}{9}\left(-2s^3 r + 3s^2 r^2 - 2sr^3 + r^4 \right).
\end{eqnarray} Recalling that $2c^2 s_+^2 = b^2 s_+ + 3a^2$ (from
the definition of $s_+$ in (\ref{eq:am17})), the function $F(s)$
can be expressed in terms of $\delta = s_+ - s \geq 0$ as follows
-
\begin{eqnarray}
&&  F(s) = \frac{\delta}{27}\left(18a^2 s_+ + 6b^2 s_+^2 - 12c^2
s_+^3\right)+ \nonumber \\ && \qquad + \delta^2\left(
\frac{3b^2}{27}s_+ + \frac{18a^2}{27} +
\delta\left(\frac{2b^2}{27} - \frac{4c^2}{9}s_+ +
\frac{c^2}{9}\delta \right)\right). \label{eq:am31}
\end{eqnarray} The coefficient of $\delta$ vanishes by virtue of
the definition of $s_+$ in (\ref{eq:am17}). We note that the
function
\begin{equation}
G\left(\delta\right) = \delta\left(\frac{2b^2}{27} -
\frac{4c^2}{9}s_+ + \frac{c^2}{9}\delta \right) \label{eq:am32}
\end{equation}
attains a minimum for
\begin{equation}
\delta_{\min} = 2s_+ - \frac{b^2}{3c^2} > s_+ \label{eq:am33}
\end{equation}
and, therefore,
\begin{equation}
\label{eq:am34} G(\delta) \geq G(s_+) = \frac{1}{27}\left(2b^2 s_+
- 9c^2 s_+^2 \right).
\end{equation} We substitute (\ref{eq:am34}) into (\ref{eq:am31})
to obtain the following lower bound for $F(s)$ -
\begin{equation}
\label{eq:am35} F(s) \geq \frac{c^2 s_+^2 + 3a^2}{27}\left( s_+ -
s \right)^2.
\end{equation}

We can analyze the function $G(s,r)$, in (\ref{eq:am30}), in an
analogous manner. Let $\ga = \frac{r}{s} \in \left[0,
\frac{1}{2}\right]$. Then
\begin{eqnarray}
&& G(s,r) = \ga s^2 \left[\frac{a^2}{3} + \frac{3b^2}{27} s -
\frac{2c^2}{9} s^2 \right] + \ga^2 s^2 \left[-\frac{a^2}{3} +
\frac{3b^2}{27} s + \frac{3c^2}{9} s^2 \right] + \ga^3 s^3
\left[-\frac{2b^2}{27} - \frac{2c^2 s}{9}  + \ga \frac{c^2 s}{9}
\right]. \label{eq:am36}
\end{eqnarray} The coefficient of $\ga $ is non-negative for all
$s\leq s_+$. Using the inequality $\ga \leq \frac{1}{2}$, one
readily obtains the following lower bound for $G(s,r)$ -
\begin{eqnarray}
\label{eq:am37}&&  G(s,r) \geq \ga s^2 \left[\frac{a^2}{3} +
\frac{3b^2}{27} s - \frac{2c^2}{9} s^2 \right] + \ga^2 s^2
\left[-\frac{a^2}{3} + \frac{2b^2}{27} s + \frac{2c^2}{9}
s^2\right] \geq \nonumber \\ && \qquad \geq
\frac{r(s-r)}{9}\left(3a^2 + b^2 s - 2c^2 s^2\right) +
\frac{5b^2}{27}r^2 s.
\end{eqnarray} Combining (\ref{eq:am35}) and (\ref{eq:am37}), the
lower bound for $0\leq s\leq s_+$ in (\ref{eq:am27}) follows.

Case(ii) The case $s \geq s_+$ can be dealt with similarly. For
any $Q\in S_0$ with $0\leq r \leq \frac{s}{2}$, we have that
\begin{equation}
\label{eq:amrevised}
 \frac{s}{\sqrt{2}} \leq \left| Q
\right|=\sqrt{\frac{2}{3}}\sqrt{\left(s^2 + r^2 -
sr \right)} \leq \sqrt{\frac{2}{3}} s.\end{equation}
For $s \geq s_+$, $|Q|^3 \geq \frac{s_+^3}{2\sqrt{2}}$ and
\begin{equation}
\beta(Q) \geq \eta \left(\frac{r^2(s - r)^2}{s^4}\right)
\label{eq:am38}
\end{equation}
where $\beta(Q)$ is the biaxiality parameter defined in
(\ref{eq:am2}) and $\eta$ is a positive constant independent of
$a^2, b^2$ or $c^2$ or $L$. Combining (\ref{eq:amrevised}),
(\ref{eq:am38}) and (\ref{eq:am18}), we readily obtain the lower
bound
\begin{eqnarray}
&& \tilde{f}_B(Q) \geq \left[\frac{ \frac{2}{3}c^2s_+^2 +
a^2}{2}\right]\left(|Q| - \sqrt{\frac{2}{3}}s_+\right)^2 +
\frac{b^2}{6\sqrt{6}}\beta(Q)|Q|^3 \geq \nonumber \\ && \qquad
\geq \left[\frac{ \frac{2}{3}c^2s_+^2 + a^2}{2}\right]\min
\left\{\frac{2}{3}\left(s -
s_+\right)^2,~\frac{1}{6}\left(\sqrt{3}s - 2 s_+\right)^2\right\}
+ \tau b^2 s_+^{3}\left(\frac{r^2(s - r)^2}{s^4}\right)
\label{eq:am39}
\end{eqnarray}
where $\tau$ is an explicitly computable positive constant.

Case (iii) Finally, we consider $Q\in S_0$ with negative order
parameters $\frac{s}{2}\leq r \leq 0$. In this case, one can
directly check that
$$ \textrm{tr} Q^3 = \frac{1}{9}\left( 2s^3 + 2r^3 - 3s^2 r - 3sr^2
\right) \leq 0$$
and therefore,
\begin{eqnarray}
\label{eq:negative1} && \tilde{f}_B(Q) = -\frac{a^2}{2}|Q|^2 -
\frac{b^2}{3}\textrm{tr} Q^3 + \frac{c^2}{4}|Q|^4 +
\frac{a^2}{3}s_+^2 + \frac{2b^2}{27} s_+^3 - \frac{c^2}{9}s_+^4
\nonumber \\ && = -\frac{a^2}{2}|Q|^2 -\frac{b^2}{3}\textrm{tr}
\left(-Q\right)^3 + \frac{c^2}{4}|Q|^4 + \frac{a^2}{3}s_+^2 +
\frac{2b^2}{27} s_+^3 - \frac{c^2}{9}s_+^4 +
\frac{2b^2}{3}\left|\textrm{tr} Q^3 \right|,
\end{eqnarray} since $\frac{b^2}{3}\textrm{tr} \left(-Q\right)^3 = -
\frac{b^2}{3}\textrm{tr} Q^3$ and $- \frac{b^2}{3}\textrm{tr} Q^3 =
\frac{b^2}{3}\left|\textrm{tr} Q^3\right|$.
The inequality (\ref{eq:negative}) follows from
(\ref{eq:negative1}) upon expressing $\textrm{tr} Q^3$ in terms of
$s$ and $r$.

For (\ref{eq:rsneg}), it suffices to note that for $s,r \leq 0$,
$\textrm{tr} Q^3 \leq 0$ and therefore,

\begin{eqnarray} &&
\tilde{f}_B(Q) \geq -\frac{a^2}{2}|Q|^2 + \frac{c^2}{4}|Q|^4
+\frac{a^2}{3}s_+^2 + \frac{2b^2}{27}s_+^3 - \frac{c^2}{9}s_+^4 =
\nonumber \\ &&  = -\frac{a^2}{3}\left(s^2 + r^2 - sr\right) +
\frac{c^2}{9}\left(s^2 + r^2 - sr\right)^2 -
\frac{s_+^3}{3}\left(\frac{b^2}{9}-\frac{c^2}{3}s_+\right).\nonumber
\end{eqnarray}
 A straightforward computation shows that the function
$$ -\frac{a^2}{3}\left(s^2 + r^2 - sr\right) + \frac{c^2}{9}\left(s^2
+ r^2 - sr\right)^2 \geq -\frac{a^4}{4c^2}$$ and
$$\frac{s_+^3}{3}\left(\frac{b^2}{9}-\frac{c^2}{3}s_+\right) <
-\frac{a^4}{4c^2}.$$ The inequality (\ref{eq:rsneg}) now follows.
 $\Box$

 \begin{remark} One can readily obtain lower bounds for $\tilde{f}_B(Q)$
in terms
 of the order parameters $\left(S,~R\right)$ in
 Proposition~\ref{prop:secondrep}, following the methods outlined
 in Proposition~\ref{prop:bulk2}. The details are omitted here for
 brevity.
\end{remark}
\begin{remark} Relation (\ref{eq:rsneg}) shows that if $f_B(Q^{(L_k)}(x))\to 0$ as $L_k\to 0$ then 
$Q^{(L_k)}$ cannot have an $(s,r)$ representation with $\frac{s}{2}<r<0$, if $L_k$ is sufficiently small.
\end{remark}

In view of Propositions~\ref{prop:max} and \ref{prop:beta}, we can
make qualitative predictions about the size of regions where a
global Landau-De Gennes minimizer $Q^*$ can have
$|Q^*|<<\sqrt{\frac{2}{3}}s_+$ and the size of regions where $Q^*$ can be
strongly biaxial.

\begin{proposition}
\label{prop:isotropic} Let $Q^*$ be a global minimizer of
$\tilde F_{LG}[Q]$ in (\ref{LDGfunctional}), in the admissible class
$\Acal_Q$ defined in (\ref{eq:am41}). Let $\Omega^* = \left\{x \in
\Omega;~|Q^*(x)| \leq \frac{1}{2}\sqrt{\frac{2}{3}}s_+\right\}$. Then
\begin{equation}
|\Omega^*| \leq \alpha \frac{L}{\left(c^2 s_+^2 +
a^2\right)}~\int_{\Omega}|\nabla n^{(0)}(x)|^2~dx,
\label{eq:am52}\end{equation} where $n^{(0)}$ is defined in (\ref{eq:n0}) and $\alpha$ is an explicitly
computable positive constant independent of $a^2, b^2, c^2$ or
$L$.
\end{proposition}

\textbf{Proof.} From Proposition~\ref{prop:beta}, we have that
\begin{equation}
\tilde{f}_B(Q^*(x)) \geq \frac{1}{\alpha}\left(c^2 s_+^2 + a^2
\right) s_+^2, \quad \textrm{$x\in \Omega^*$} \label{eq:am53}
\end{equation} for some explicitly computable positive constant $\alpha$,
since
$|Q^*| \leq  \frac{1}{2}\sqrt{\frac{2}{3}}s_+= \frac{1}{\sqrt{6}}s_+$ on
$\Omega^*$. On the other hand, recalling the definition of $Q^{(0)}$ in (\ref{eq:Q0}) and since $Q^*$ is a global minimizer
of $\tilde F_{LG}[Q]$, we have that
\begin{equation}
\label{eq:am54} \int_{\Omega^*} \tilde{f}_B(Q^*(x))~dx \leq
\mathcal{F}_{LG}[Q^{(0)}] = \int_{\Omega}\tilde{f}_B(Q^{(0)}) +
\frac{L}{2}|\nabla Q^{(0)}|^2~dx = L s_+^2 \int_{\Omega}|\nabla
n^{(0)}|^2~dx,
\end{equation} since $\tilde{f}_B(Q^{(0)})=0 $ everywhere in
$\Omega$. Substituting (\ref{eq:am53}) into (\ref{eq:am54}), we
obtain
\begin{equation}
\label{eq:am55} \frac{1}{\alpha}\left(c^2 s_+^2 + a^2 \right)
s_+^2 |\Omega^*| \leq L s_+^2 \int_{\Omega}|\nabla n^{(0)}|^2~dx,
\end{equation}
from which the inequality (\ref{eq:am52}) follows. $\Box$

\begin{proposition}
\label{prop:nondefect} Let $Q^*$ be a global minimizer of
$\tilde F_{LG}[Q]$ in (\ref{LDGfunctional}), in the admissible class
$\Acal_Q$ defined in (\ref{eq:am41}). Let $\Omega^{\la} = \left\{ x
\in \Omega;~|Q^*(x)| \geq \frac{1}{2}\sqrt{\frac{2}{3}}s_+,~ \beta(Q(x))> \la
\right\}$ for some positive constant $\la$. Then,
\begin{equation}
\label{eq:am56} |\Omega^{\la}| \leq \alpha
\frac{L}{\la s_+ b^2} \int_{\Omega}|\nabla
n^{(0)}|^2\,dx
\end{equation}where $n^{(0)}$ is defined in (\ref{eq:n0}) and $\alpha$ is an explicitly
computable positive constant independent of $a^2, b^2, c^2$ or
$L$.
\end{proposition}

\textbf{Proof.} From Proposition~\ref{prop:beta}, we have that
\begin{equation}
\tilde{f}_B(Q^*(x)) \geq \frac{b^2}{6\sqrt{6}}\beta(Q^*(x))|Q^*(x)|^3 \geq
\frac{1}{\alpha}b^2 \la s_+^3 \quad \textrm{$x\in \Omega^{\la}$}
\label{eq:am57}
\end{equation} for some explicitly computable positive constant $\alpha$,
since
$|Q^*| \geq  \frac{1}{2}\sqrt{\frac{2}{3}}s_+= \frac{1}{\sqrt{6}}s_+$ on
$\Omega^{\la}$. On the other hand, recalling the definition of $Q^{(0)}$, (\ref{eq:Q0}), and since $Q^*$ is a global
minimizer of $\tilde F_{LG}[Q]$, we have that
\begin{equation}
\label{eq:am58} \int_{\Omega^{\la}} \tilde{f}_B(Q^*(x))~dx\leq
\int_{\Omega}\tilde{f}_B(Q^{(0)}) + \frac{L}{2}|\nabla Q^{(0)}|^2~dx = L
s_+^2 \int_{\Omega}|\nabla n^{(0)}|^2~dx
\end{equation} since $\tilde{f}_B(Q^{(0)})=0 $ everywhere in
$\Omega$. Substituting (\ref{eq:am57}) into (\ref{eq:am58}), we
obtain
\begin{equation}
\label{eq:am59} \frac{1}{\alpha}b^2 \la s_+^3 |\Omega^{\la}| \leq
L s_+^2 \int_{\Omega}|\nabla n^{(0)}|^2~dx,
\end{equation}
from which the inequality (\ref{eq:am56}) follows. $\Box$

Proposition~\ref{prop:isotropic} is relevant to the size of defect
cores in global energy minimizers whereas
Proposition~\ref{prop:nondefect} is relevant to the equilibrium
behaviour far away from the defect cores.

\subsection{Analyticity and uniaxiality}
\label{section:analyticitybiaxiality}

\par We define a new biaxiality parameter $\tilde\beta(Q)$ as follows:

$$\tilde\beta(Q)\stackrel{def}{=}(\textrm{tr}(Q^2))^3-6(\textrm{tr}(Q^3))^2.$$
Then $\tilde\beta(Q)\ge 0$ with $\tilde\beta(Q)=0$ if and only if
$Q$ is uniaxial i.e. $Q=s\left(n\otimes n-\frac{1}{3}Id\right)$
for some $s\in\mathbb{R}\setminus\{0\},n\in\mathbb{S}^2$ or $Q=0$. The function
$\tilde\beta(Q)$ is a real analytic function of $Q$ and this is
particularly important given that global energy minimizers of the
functional $\mathcal{F}$ (subject to smooth boundary conditions)
are real analytic:

\begin{proposition} Let $\Omega$ be a simply-connected bounded open set. Let $Q^{(L)}$ be a global energy minimizer of $\tilde F_{LG}[Q]$ in (\ref{LDGfunctional}) in the admissible space $\Acal_Q$. 
Then $Q^{(L)}$ is real analytic in $\Omega$. \label{analyticity}
\end{proposition}

\smallskip\par{\bf Proof.} \par We drop the superscript $L$ from $Q^{(L)}$ for convenience. As
$-\frac{a^2}{2}\textrm{tr}\left(Q^2\right) -
\frac{b^2}{3}\textrm{tr}\left(Q^3\right) +
\frac{c^2}{4}\left(\textrm{tr}Q^2\right)^2$  is bounded from below
(see also the {\it Appendix}) we have that there exists an $H^1$
global energy minimizer satisfying the Euler-Lagrange system:
$$L\Delta Q_{ij}=-a^2Q_{ij}-b^2\left(Q_{ik}Q_{kj}-\frac{\delta_{ij}}{3}\textrm{tr}(Q^2)\right)+c^2Q_{ij}\textrm{tr}(Q^2)$$
\par For $Q$ an $H^1$ solution of the equation one uses
$H^1\hookrightarrow L^6$ (in $\mathbb{R}^3$) and H\"{o}lder's
inequality to obtain that the right hand side of each  equation is
in $L^2$. Elliptic regularity gives that $Q\in H^2\hookrightarrow
W^{1,6}\hookrightarrow{L^\infty}$ hence the right hand side of the
equation is in $H^1$. Elliptic regularity gives $Q\in H^3$ and one
can continue bootstrapping to obtain the full regularity allowed
by the regularity of boundary data and that of the domain.
\par In order to prove the analyticity we use a general abstract result
due to A.Friedman,\cite{friedman}.  We define growth classes as follows:
let $M_n$ be a sequence of
positive numbers. Then a function $F:C^\infty(D)\to\mathbb{C}$,
with $D\subset\mathbb{R}^d$ an open set, belongs to the class
$C\{M_n;D\}$ if for  any closed subset $D_0\subset D$ there exist
constants $H_0,H$ with

$$|\partial^j F(x)|\le H_0H^j M_j,x\in D_0$$ where we have used
multiindex-notation $(\partial^j
F=\partial^{j_1}_1\dots\partial^{j_d}_d F;j=\Sigma_{i=1}^d j_i)$.
Let us observe that $C\{n!;D\}$ is the class of functions analytic
in $D$.
\par In \cite{friedman}  the following theorem is proved for general
elliptic systems:

\begin{proposition}(\cite{friedman},p.45) Let $u(x)$ be a real
solution of the elliptic system
$$\Phi_l(x;u,\nabla u,\nabla^2 u,\dots, \nabla
^{2m}u)=0,x\in\Omega\subset\mathbb{R}^d;u\in\mathbb{R}^N,\,l=1,\dots,N$$
in $\Omega\subset\mathbb{R}^d$. Let $E$ be some open set
containing $E_1\stackrel{def}{=}\{u(x),\nabla
u(x),\dots,\nabla^{2m} u(x);x\in\Omega\}$. Assume that:
\par {\it (i)} $\Phi_l\in C\{M_n;\Omega\times E\}$ and that the $M_n$
satisfy the monotonicity conditions

\par {\it (ii)} ${n \choose i} M_i M_{n-i}\le AM_n;\,\, 0\le i\le
n,\,n\in\mathbb{N}$ for some $A>0$.
\par If $u\in C^{2m+\alpha}(\Omega),0<\alpha<1$ then $u\in
C\{M_{n-2m+1};\Omega\}$ (where $M_{-i}=1$ for $i\in\mathbb{N}$)
\end{proposition}

\smallskip\par In our case, for the system (\ref{ELeq}) we have
$m=1$ and $\Phi_l$ is analytic hence of class $C\{n!;\Omega\}$.
The constants $M_n=n!$ satisfy the monotonicity conditions {\it
(ii)} in the theorem, with $A=1$. We have that $Q\in
C^\infty(\Omega)$ and hence by the theorem $Q$ is in the class
$C\{(n-1)!;\Omega\}$ therefore real analytic. $\Box$

\smallskip
\begin{proposition} Let $Q$ be a real analytic function
$Q:\Omega\subset\mathbb{R}^3\to S_0$. Then the set where
$Q$ is uniaxial or isotropic is either the whole of $\Omega$ or has zero Lebesgue
measure. \label{zeroL}
\end{proposition}

\smallskip\par{\bf Proof.} If there is no $x\in\Omega$ such that
$\tilde\beta(Q(x))\not =0$ then $Q$ is uniaxial or isotropic everywhere.  If there
exists a $P\in\Omega$ such that $\tilde\beta(Q(P))\not=0$ then let us
consider the lines passing through $P$.  The restriction of
$Q$ to any such line is real analytic and then so is
$\tilde\beta(Q)$. Thus $\tilde\beta(Q)$  has at most countably
many zeroes on such a line.  We claim that this implies that the set of
zeroes of
$\tilde\beta(Q)$  in $\Omega$ is of measure zero.
\par We assume, without loss of generality, that $P=0$. We denote
$N^*=N\setminus\{0\}$ and decompose
$\mathbb{R}^n\cap\Omega=\cup_{n\in\mathbb{N}^*}
\left(\overline{B_{\frac{1}{n}}\setminus
B_{\frac{1}{n+1}}}\cap\Omega\right)\cup\left(\cup_{n\in\mathbb{N}^*}
\left(\overline{B_{n+1}\setminus
B_{n}}\cap\Omega\right)\right)$. We
claim that for any $n,\frac{1}{n}\in\mathbb{N}^*$  the set
$\left(\tilde\beta(Q)\right)^{-1}(0)\cap\left(\Omega\cap
\overline{B_{\frac{1}{n}}\setminus B_{\frac{1}{n+1}}} \right)$ is a set of
measure
zero. This implies that $\tilde\beta(Q)^{-1}(0)\cap\Omega$, which is a
countable
union of sets as before, is also
a set of measure zero.

\par We consider the bi-Lipschitz functions
$$f_n:[\frac{1}{n+1},\frac{1}{n}]\times\underbrace{[0,\pi]\times\dots\times
[0,\pi]}_{n-2\, times}\times[0,2\pi)\to
\overline{B_{\frac{1}{n}}\setminus B_{\frac{1}{n+1}}},\forall
n,\frac{1}{n}\in\mathbb{N}$$ that realize the
change of coordinates from polar to usual cartesian  coordinates.
\par We have that $f_n^{-1}\left(\tilde\beta(Q)^{-1}(0)\cap\Omega\cap
\overline{B_{\frac{1}{n}}\setminus B_{\frac{1}{n+1}}} \right)\subset
[\frac{1}{n+1},\frac{1}{n}]\times\underbrace{[0,\pi]\times\dots\times
[0,\pi]}_{n-2\, times}\times[0,2\pi)$. We recall that the Lebesgue
measure $\mu$ on
the $n$-dimensional product space
$[\frac{1}{n+1},\frac{1}{n}]\times\underbrace{[0,\pi]\times\dots\times
[0,\pi]}_{n-2\, times}\times[0,2\pi)$  is the completion  of the product
measure
$\mu_1\times \mu_2$ where $\mu_1$ is the $1$ dimensional Lebesgue  measure on
$[\frac{1}{n+1},\frac{1}{n}]$ and $\mu_2$ is the $n-1$ dimensional
Lebesgue measure on
$\underbrace{[0,\pi]\times\dots\times
[0,\pi]}_{n-2\, times}\times[0,2\pi)$. Then for any set $E\subset
[\frac{1}{n+1},\frac{1}{n}]\times\underbrace{[0,\pi]\times\dots\times
[0,\pi]}_{n-2\, times}\times[0,2\pi)$ we have $$(\mu_1\times
\mu_2)(E)=\int_{\underbrace{[0,\pi]\times\dots\times
[0,\pi]}_{n-2\, times}\times[0,2\pi)}\mu_1(E^y)\mu_2(dy)$$ where $E^y=\{x\in
[\frac{1}{n+1},\frac{1}{n}], (x,y_1,\dots,y_{n-1})\in E\}\subset
[\frac{1}{n+1},\frac{1}{n}]$. In our case, letting
$$E\stackrel{def}{=}f_n^{-1}\left(\tilde\beta(Q)^{-1}(0)\cap\Omega\cap
\overline{B_{\frac{1}{n}}\setminus B_{\frac{1}{n+1}}}\right)$$ we have
that $E^y$ is
made of finitely many points  for almost all $y\in
\underbrace{[0,\pi]\times\dots\times
[0,\pi]}_{n-2\, times}\times[0,2\pi)$ (as  a consequence of the first
paragraph in
this proof; because $E^y$ is just the set of the distances to $P$ of the
uniaxial or isotropic
points that are in
 $\Omega\cap\overline{B_{\frac{1}{n}}\setminus B_{\frac{1}{n+1}}}$,  on a
a segment
through $P$, segment that has in
polar coordinates  the direction $y\in \underbrace{[0,\pi]\times\dots\times
[0,\pi]}_{n-2\, times}\times[0,2\pi)$ ). Thus $\mu_1(E^y)=0, \mu_2-\,a.e.
\,y$ hence $\mu_1\times \mu_2(E)=0$ thus $\mu(E)=0$.

\par  As bi-Lipschitz functions carry sets of measure zero into
sets of measure zero we have that  $\tilde\beta(Q)^{-1}(0)\cap\Omega\cap
\overline{B_{\frac{1}{n}}\setminus B_{\frac{1}{n+1}}}$ is a set of measure
zero. On the other hand $\tilde\beta(Q)^{-1}(0)\cap \Omega$ is a countable
union of  sets as before, hence it has measure zero. $\Box$

\begin{corollary}
 \label{eigenvectornew} Let $Q^{(L)}$ be a global minimizer of
 $\tilde F_{LG}[Q]$ in the admissible class $\Acal_Q$. Then there exists a set of measure zero, possibly empty, $\Omega_0$ in $\Omega$ such that 
 the  eigenvectors of $Q^{(L)}$ are smooth at all points $x\in\Omega\setminus\Omega_0$. The uniaxial-biaxial interfaces,
isotropic-uniaxial or isotropic-biaxial interfaces are contained in $\Omega_0$.
\end{corollary}
\textbf{Proof.} The global minimizer $Q^{(L)}\in
C^{\infty}\left(\Omega;\Acal\right)$. The eigenvectors of $Q^{(L)}$
 have the same degree of regularity as $Q^{(L)}$ on sets
$K\subset\Omega$, where $Q^{(L)}$ has the same number of distinct
eigenvalues i.e. where $Q^{(L)}$ is either  biaxial or
 uniaxial or isotropic, \cite{nomizu}, but not necessarily otherwise \cite{kato}. If $Q^{(L)}$ is uniaxial everywhere then $\Omega_0=\emptyset$. 
If $Q^{(L)}$ is either uniaxial or isotropic on the whole of $\Omega$ (i.e. $\tilde\beta(Q^{(L)})=0$ in $\Omega$), with  $Q^{(L)}\not=0$ at some point in $\Omega$, then let 
$\tilde\Omega=\{x\in\Omega, Q^{(L)}(x)=0\}$ denote the zero-set of $Q^{(L)}$. Let us observe that $\tilde\Omega=\left(|Q|^2\right)^{-1}(0)$ and $|Q|^2$ is an analytic function. By an argument similar to the proof of Proposition ~\ref{zeroL} and since $Q(x)\not=0$  for at least one point $x\in\Omega$, we have that 
$\tilde\Omega$ has measure zero and we take $\Omega_0\stackrel{def}{=}\tilde \Omega$.

 If $Q^{(L)}$ is biaxial somewhere then Proposition ~\ref{zeroL} shows that the set of points where $\tilde \beta(Q)=0$ has measure zero. We denote this set
 by $\Omega_0$ and observe that $\Omega\setminus\Omega_0$ is an open set and the eigenvectors have the same regularity as $Q^{(L)}$ on $\Omega\setminus\Omega_0$, see \cite{nomizu}. $\Box$

\bigskip
\section{Acknowledgements}
A.~Majumdar was supported by a Royal Commission for the Exhibition
of 1851 Research Fellowship till October 2008. She is now
supported by Award No. KUK-C1-013-04 , made by King Abdullah
University of Science and Technology (KAUST). A.~Zarnescu is
supported by the EPSRC Grant EP/E010288/1  - Equilibrium Liquid
Crystal Configurations: Energetics, Singularities and
Applications. We thank John Ball and Christof Melcher for stimulating discussions.

\vskip 20ex
\begin{center}
{\bf Appendix}
\end{center}

\begin{proposition} \cite{bm}
\label{prop:bulk1}
Consider the bulk energy density $f_B(Q)$ given by
\begin{equation}
f_B(Q) = -\frac{a^2}{2}\textrm{tr} Q^2 - \frac{b^2}{3}\textrm{tr}
Q^3 + \frac{c^2}{4}\left(\textrm{tr} Q^2\right)^2. \label{eq:am15new}
\end{equation}
Then $f_B(Q)$ attains its minimum for uniaxial $Q$-tensors of the
form
\begin{equation}
Q = s_+\left(n\otimes n - \frac{1}{3}\right), \label{eq:am16new}
\end{equation}
where $n:\Omega \rightarrow S^2$ is a unit eigenvector of $Q$ and
\begin{equation}
s_+ = \frac{b^2 + \sqrt{b^4 + 24 a^2 c^2}}{4c^2}. \label{eq:am17new}
\end{equation}
\end{proposition}

\textbf{Proof.}
Proposition~\ref{prop:bulk1} has been proven in \cite{bm}. We reproduce
the proof
here for completeness.

We recall that for a symmetric, traceless matrix Q of the form
$$Q = \sum_{i=1}^3 \la_i e_i \otimes e_i, $$ $\textrm{tr}Q^n =
\sum_{i=1}^{3}\la_i^n$
subject to the tracelessness condition so that the bulk energy
density $f_B$ in (\ref{eq:am15new}) only depends on the eigenvalues
$\lambda_1,
\lambda_2$ and $\lambda_3$. Then the stationary points of the bulk
energy density $f_B$ are given by the stationary points of the
function $f:\Rr^3 \rightarrow \Rr$ defined by
\begin{equation}
f\left(\la_1,\la_2,\la_3\right) = -\frac{a^2}{2}\sum_{i=1}^{3}\la_i^2
- \frac{b^2}{3}\sum_{i=1}^{3}\la_i^3 +
\frac{c^2}{4}\left(\sum_{i=1}^{3}\la_i^2 \right)^2 -
2\delta\sum_{i=1}^{3}\la_i .\label{eq:16new}
\end{equation} where we have recast $f_B$ in terms of the eigenvalues and
introduced
a Lagrange multiplier $\delta$ for the tracelessness condition.

The equilibrium equations are given by a system of three algebraic
equations
\begin{equation}
\frac{\partial f}{\partial \lambda_i} = 0 \Leftrightarrow -a^2 \la_i
- b^2\la_i^2 + c^2\left(\sum_{k=1}^{3}\la_k^2 \right)\la_i = 2\delta
\quad \textrm{for $i=1 \ldots 3$,} \label{eq:17new}
\end{equation}
or equivalently
\begin{equation}
\left( \la_i - \la_j \right)\left[ -a^2 - b^2\left(\la_i + \la_j
\right) + c^2 \sum_{k=1}^{3}\la_k^2 \right] = 0 \quad 1\leq i< j
\leq 3.\label{eq:18new}
\end{equation}
Let $\left\{\lambda_i\right\}$ be a solution of the system
(\ref{eq:17new}) with three distinct eigenvalues $\la_i \neq \la_2
\neq \la_3$.  We consider equation (\ref{eq:18new}) for the pairs
$\left(\la_1,\la_2\right)$ and $\left(\la_1,\la_3\right)$. This
yields two equations
\begin{eqnarray}
&& -a^2 - b^2\left(\la_1 + \la_2\right) + c^2\sum_{k=1}^{3}\la_k^2 = 0
\nonumber \\&& -a^2 - b^2\left(\la_1 + \la_3\right) +
c^2\sum_{k=1}^{3}\la_k^2 = 0 \label{eq:19new}
\end{eqnarray}
from which we obtain
\begin{equation}
-b^2\left(\la_2 - \la_3 \right) = 0, \label{eq:20new}
\end{equation} contradicting our initial hypothesis $\la_2 \neq
\la_3$. We, thus, conclude that a stationary point of the bulk
energy density must have at least two equal eigenvalues and
therefore correspond to either a uniaxial or isotropic liquid
crystal state.

We consider an arbitrary uniaxial state given by
$\left(\la_1,\la_2,\la_3\right) =
\left(\frac{2s}{3},-\frac{s}{3},-\frac{s}{3}\right)$ and the
corresponding Q-tensor is $Q = s\left(e_1\otimes e_1 -
\frac{1}{3}Id\right)$. The function $f_B$ is then a quartic
polynomial in the order parameter $s$ ie.
\begin{equation}
f_B(s) = \frac{s^2}{27}\left( -9a^2 - 2b^2s + 3c^2s^2\right)
\label{eq:f1new}
\end{equation}
and the stationary points are solutions of the algebraic equation
$\frac{d f_B}{d s} = 0$,
\begin{equation}
\frac{d f_B}{d s} = \frac{1}{27}\left(-18a^2s - 6b^2s^2 + 12c^2s^3\right)
= 0. \label{eq:f2new}
\end{equation}
The cubic equation (\ref{eq:f2new}) admits three solutions;
\begin{equation}
s = 0 \quad \textrm{and $s_{\pm} = \frac{b^2 \pm \sqrt{b^4 + 24 a^2
c^2}}{4c^2}$}\label{eq:f3}
\end{equation}
where
\begin{equation}
f_B(0) = 0 \quad \textrm{and} \quad f_B(s_+) < f_B(s_-) < 0. \label{eq:f4}
\end{equation} Symmetry considerations show that we obtain the same set of
stationary points for the remaining two uniaxial choices. The global
minimizer is,
therefore, a uniaxial $Q$-tensor of the form
 \begin{equation}
Q = s_+\left(n\otimes n - \frac{1}{3}Id\right), \,n\in\mathbb{S}^2  \label{eq:uniaxial2}
\end{equation} where $s_+$ has been defined in
(\ref{eq:am17new}).$\Box$

\begin{lemma}
\label{lem:1} Let $Q \in S_0$. We define the biaxiality parameter $\beta(Q)$
to be
\begin{equation}
\label{eq:am2new} \beta(Q) = 1 -
6\frac{\left(\textrm{tr}Q^3\right)^2}{\left(\textrm{tr}Q^2\right)^3}.
\end{equation} (i) The biaxiality parameter $\beta(Q)\in\left[0,1\right]$ and
$\beta(Q)=0$ if and only if $Q$ is purely uniaxial i.e. if $Q$ is
of the form, $Q=s\left(n\otimes n-\frac{1}{3}Id\right)$ for some
$s\in\mathbb{R},n\in\mathbb{S}^2$. (ii) The biaxiality parameter,
$\beta(Q)$, can be bounded in terms of the ratio $\frac{r}{s}$,
where $(s,~r)$ are the scalar order parameters in
Proposition~\ref{prop:rep1} . These bounds are given by
\begin{equation}
\label{eq:am3version} \frac{1}{2}\left(1 - \sqrt{1 -
\sqrt{\beta}}\right) \leq \frac{r}{s} \leq \frac{1}{2}\left(1 +
\sqrt{1 - \sqrt{\beta}}\right).
\end{equation} Equivalently,
\begin{equation}
\label{eq:am3newversion} \frac{1 -\sqrt{1 - \sqrt{\beta}}}{3 +
\sqrt{1 - \sqrt{\beta}}} \leq \frac{R}{S}\leq \frac{1 +\sqrt{1 -
\sqrt{\beta}}}{3 - \sqrt{1 - \sqrt{\beta}}} \end{equation} where
$(S,~R)$ are the order parameters in
Proposition~\ref{prop:secondrep}. Further $\beta(Q)=1$ if and only
if $r=\frac{s}{2}$ or if and only if $\frac{R}{S}=\frac{1}{3}$.
(iii)For an arbitrary $Q\in S_0$, we have that
\begin{equation}
-\frac{|Q|^3}{\sqrt{6}}\left(1 - \frac{\beta}{2}\right) \leq
\textrm{tr} Q^3 \leq \frac{|Q|^3}{\sqrt{6}}\left(1 -
\frac{\beta}{2}\right).\label{eq:am9new}
\end{equation}
\end{lemma}
\textbf{Proof:} (i) The quantity $\beta(Q)$ is known as the
biaxiality parameter in the liquid crystal literature
\cite{mg} and it is well-known that
$\beta(Q)\in\left[0,1\right]$. We present a simple proof here for
completeness.

Following Proposition~\ref{prop:rep1}, we represent an arbitrary
$Q\in S_0$ as
\begin{equation} \label{eq:am1}Q =
s\left(n\otimes n - \frac{1}{3}Id\right) + r \left(m\otimes m -
\frac{1}{3}Id\right) \quad 0\leq r\leq
\frac{s}{2}~or~\frac{s}{2}\leq r \leq 0.
\end{equation}Since
$6\frac{\left(\textrm{tr}Q^3\right)^2}{\left(\textrm{tr}Q^2\right)^3}
\geq 0$, the inequality $\beta(Q) \leq 1$ is trivial. To show
$\beta(Q)\geq 0$, we use the representation (\ref{eq:am1}) to
express $\textrm{tr}Q^3$ and $\textrm{tr}Q^2$ in terms of the
order parameters $s$ and $r$.
\begin{eqnarray}
\label{eq:biaxiality3} && \textrm{tr}Q^3 =
\frac{1}{9}\left(2s^3 + 2r^3 - 3s^2 r - 3sr^2\right)\nonumber \\
&& \textrm{tr}Q^2 = \frac{2}{3}\left( s^2 + r^2 - sr \right)
\end{eqnarray} A straightforward calculation shows that
$$\left(\textrm{tr}Q^3\right)^2 = \frac{1}{81}\left(4s^6 +
4r^6 - 12 s^5 r - 12 s r^5 + 26 s^3 r^3 - 3 s^4 r^2 - 3s^2 r^4
\right)$$ and
$$\left(\textrm{tr}Q^2\right)^3 = \frac{8}{27}\left(s^6 + r^6
- 3 s^5 r - 3 s r^5 - 7 s^3 r^3 + 6 s^2 r^4 + 6 s^4 r^2 \right).$$
One can then directly verify that
\begin{equation}
\label{eq:biaxiality4}\left(\textrm{tr}Q^2\right)^3 -
6\left(\textrm{tr}Q^3\right)^2 = 2 s^2 r^2\left( s- r
\right)^2 \geq 0
\end{equation}
as required. It follows immediately from (\ref{eq:biaxiality4})
that $\beta(Q)=0$ if and only if either $s=0, r=0$ or $s=r$. From
(\ref{eq:am1}), the three cases, $s=0, r=0$ and $s=r$, correspond
to uniaxial nematic states (in fact all uniaxial states can be
described by one of these three conditions) and therefore,
$\beta(Q)=0$ if and only if $Q$ is uniaxial.

(ii) From Proposition~\ref{prop:rep1}, it suffices to consider
$Q$-tensors with either $0\leq r\leq \frac{s}{2}$ or
$\frac{s}{2}\leq r \leq 0$. Let $\gamma = \frac{r}{s}$, then
$\gamma \in \left[0,\frac{1}{2}\right]$ for the two cases under
consideration. The biaxiality parameter, $\beta(Q)$, can be
expressed in terms of the ratio $\gamma$ as follows
\begin{equation}
\label{eq:am5} \frac{\left(2 - 3\gamma - 3\gamma^2 +
2\gamma^3\right)^2}{\left(1 - \gamma + \gamma^2\right)^3} = 4\left(1
- \beta \right).
\end{equation}
From (\ref{eq:biaxiality4}), we have that
\begin{equation}
\label{eq:am6} \left(2 - 3\gamma - 3\gamma^2 + 2\gamma^3\right)^2 =4
\left(1 - \gamma + \gamma^2\right)^3 - 27\gamma^2\left(1 -
\gamma\right)^2, \end{equation} which in turn, yields the
following equality
\begin{equation}
\frac{27\gamma^2 \left(1 - \gamma\right)^2}{\left(1 - \gamma +
\gamma^2\right)^3} = 4\beta. \label{eq:am7}
\end{equation} Noting that for
$\gamma=\left[0,\frac{1}{2}\right]$, the polynomial $1-\gamma +
\gamma^2 \geq \frac{3}{4}$, we obtain the following upper bound
\begin{equation}
\label{eq:am8} \beta \leq 16\gamma^2\left(1 - \gamma \right)^2
\end{equation}
and the bounds (\ref{eq:am3version}) readily follow from
(\ref{eq:am8}).

One can readily see from (\ref{eq:am7}) that $\beta(Q)=1$ if and
only if $\frac{r}{s}=\frac{1}{2}$ . The bounds
(\ref{eq:am3newversion}) follow directly from
(\ref{eq:am3version}) on noting that
$$r = 2R~and~s = S + R.$$ One can see directly from
(\ref{eq:am3newversion}) that if $\beta=1$, then
$\frac{R}{S}=\frac{1}{3}$. On the other hand, if
$\frac{R}{S}=\frac{1}{3}$, then $\frac{r}{s}=\frac{1}{2}$ and
(\ref{eq:am7}) implies that $\beta(Q)=1$. The claims in (ii) now
follow.

(iii) From the definition of the biaxiality parameter in
(\ref{eq:am2new}), we necessarily have that
\begin{equation}
\label{eq:am13}\textrm{tr} Q^3 = \pm
\frac{|Q|^3}{\sqrt{6}}\sqrt{1-\beta(Q)}.
\end{equation}It is
easily checked that
\begin{equation}
\sqrt{1 - \beta} \leq 1 - \frac{\beta}{2} \label{eq:am14}
\end{equation}The bounds (\ref{eq:am9new}) follow from combining
(\ref{eq:am13}) and
(\ref{eq:am14}) .$\Box$


\begin{thebibliography}
\smallskip
\par
\bibitem{lieb} F.J. Almgren and E.H. Lieb,  Singularities of energy minimizing maps 
from the ball to the sphere: examples, counterexamples, and bounds.  
Ann. of Math. (2)  128  (1988),  no. 3, 483--530.
\bibitem{bz}  J.M. Ball and A. Zarnescu, Orientability and energy
minimization for liquid crystals, in preparation
\bibitem{bbh} F. Bethuel, H. Brezis  and F.H\'{e}lein,
Asymptotics for the minimization of a Ginzburg-Landau functional.
Calc. Var. Partial Differential Equations  1  (1993), no. 2,
123--148
\bibitem{bc}  F. Bethuel and D.Chiron,  Some questions related to the lifting problem in Sobolev spaces. 
Perspectives in nonlinear partial differential equations, 125--152, Contemp. Math., 446, Amer. Math. Soc., Providence, RI, 2007.
\bibitem{brezis}H. Brezis, The interplay between analysis and topology in
some nonlinear PDE problems. Bull. Amer. Math. Soc. (N.S.) 40 (2003), no.
2, 179--201 (electronic).
\bibitem{remarks} Y. Chen, and F.Lin, Remarks on approximate harmonic maps.
 Comment. Math. Helv.  70  (1995),  no. 1, 161--169
 \bibitem{gartland} T.~Davis and E.~Gartland, Finite element analysis of
the Landau--De Gennes minimization problem for liquid crystals.
SIAM Journal of Numerical Analysis, 35, 336-362 (1998).
 \bibitem{ericksen} J. L. Ericksen,  Liquid crystals with variable degree
of orientation.
 Arch. Rational Mech. Anal.  113  (1990),  no. 2, 97--120
 \bibitem{evans} L.~Evans, Partial Differential Equations. American
Mathematical Society, Providence, 1998.
 \bibitem{fraenkel} L.E. Fraenkel,  On regularity of the boundary in the
theory of Sobolev spaces.  Proc. London Math. Soc. (3)  39  (1979), no. 3,
385--427
\bibitem{oseenfrank} F.C. Frank,  On the theory of liquid crystals.
 Disc. Faraday Soc., 25(1958)1
 \bibitem{friedman} A. Friedman, On the regularity of the solutions of
nonlinear
 elliptic and parabolic systems of partial differential equations.
 J. Math. Mech. 7, 43-59 (1958)
\bibitem{dg} P. G. De Gennes,  The physics of liquid crystals.
Oxford, Clarendon Press. 1974
\bibitem{giaquinta} M. Giaquinta, {\it Multiple integrals in the calculus
of variations and nonlinear elliptic systems.}
 Annals of Mathematics Studies, 105. Princeton University Press,
Princeton, NJ, 1983.
\bibitem{gilbarg} D.Gilbarg and N.Trudinger, Elliptic Partial Differential
Equations of Second Order.
 Springer, 224, 2, 1977
\bibitem{partialcrystal} R. Hardt, D. Kinderlehrer and F. H. Lin,
Existence and partial regularity of static liquid crystals
configurations, Comm. Math. Phys., 105 (1986), 547-570
\bibitem{kato} T. Kato, Perturbation theory for linear operators, rundlehren der 
Mathematischen Wissenschaften, Band 132. Springer-Verlag, Berlin-New York, 1976
 \bibitem{lin} F.~H.~Lin and C.~Liu, Static and Dynamic Theories of Liquid
Crystals. Journal of Partial Differential Equations, 14, no. 4,
289--330  (2001).
\bibitem{linpoon} F.  Lin and C. Poon,  On Ericksen's model for liquid
crystals. J. Geom. Anal. 4 (1994), no. 3, 379--392
\bibitem{linriviere} F. Lin and T.Rivi\'ere, Complex Ginzburg-Landau
equations in high dimensions and
codimension two area minimizing currents.  J. Eur. Math. Soc. (JEMS)  1
(1999),  no. 3, 237--311.
\bibitem{bm} A. Majumdar, Equilibrium order parameters of liquid crystals
in the Landau--De Gennes theory, preprint.
\bibitem{virgadematteis} G. De Matteis and E.G. Virga,  Tricritical points in biaxial liquid crystal phases,
Phys. Rev. E 71, 061703 (2005)
\bibitem{mg} S. Mkaddem and E. C. Gartland,
 Fine structure of defects in radial nematic droplets.  Phys. Rev. E, 62, 6694 --
 6705, 2000
 \bibitem{moser}  R. Moser,  Partial regularity for harmonic maps and
related problems. World Scientific Publishing , Hackensack, NJ, 2005.
\bibitem{newtonmottram} N.J.Mottram and C.Newton,  Introduction to
\textbf{Q}-tensor Theory.
 University of Strathclyde, Department of Mathematics, Research Report,
10, 2004
\bibitem{nomizu} K. Nomizu, Characteristic roots and vectors of a differentiable family of symmetric matrices, Linear and
Multilinear Algebra 1, 159-162 (1973)
\bibitem{smallelastic}  E. B. Priestley,  P. J Wojtowicz and  P.
Sheng, Intorduction to Liquid Crystals, Plenum, New York, 1975
\bibitem{rosso&virga} R. Rosso and  E.Virga, Metastable nematic hedgehogs.
J. Phys. A: Math. Gen. 29, 4247 -- 4264,  1996
 \bibitem{schoen} Schoen, R. {\it Analytic Aspects of the Harmonic Map
Problem.} Seminar
on Nonlinear Partial Differential Equations. Chern, S.S., Ed.; MSRI
Publications 2, Springer-Verlag, 1984.
\bibitem{harmreg} R. Schoen and K. Uhlenbeck,  A Regularity Theory for Harmonic 
Mappings. J. Diff. Geom. 1982, 17, 307-335. 
\bibitem{taylor} M.E. Taylor, Partial differential equations. III.
Nonlinear equations.
Applied Mathematical Sciences, 117. Springer-Verlag, New York, 1997
\bibitem{virga}E. G. Virga, Variational theories for liquid crystals.
Chapman and Hall, London 1994
\end{thebibliography}
\end{document}